\documentclass[notitlepage]{amsart}
\usepackage{amssymb}

\usepackage{graphicx}
\usepackage{amscd}

\newtheorem{theorem}{Theorem}
\theoremstyle{plain}

\newtheorem{corollary}{Corollary}

\newtheorem{definition}{Definition}

\newtheorem{lemma}{Lemma}

\newtheorem{proposition}{Proposition}
\newtheorem{remark}{Remark}

\numberwithin{equation}{section}

\input{diagrams.tex}
\newarrow{Onto}{-}{-}{-}{-}{->>}
\newarrow{Into}{C}{-}{-}{-}{>}
\newarrow{Corr}{<}{-}{-}{-}{>}
\newarrow{Ev}{|}{-}{-}{-}{>}
\input{tcilatex}
\begin{document}
\title{\emph{Twisted} Internal coHom Objects in the Category of Quantum Linear
Spaces}
\author{S. Grillo}
\author{H. Montani}
\address{\textit{Centro At\'{o}mico Bariloche and Instituto Balseiro}\\
\textit{\ 8400-S. C. de Bariloche}\\
\textit{\ Argentina}}
\date{December 11, 2001}

\begin{abstract}
Adapting the idea of twisted tensor products \cite{cap} to the category of
finitely generated algebras, we define on its opposite, the category $%
\mathsf{QLS}$ of quantum linear spaces \cite{man0}\cite{man1}, a family of
objects $\underline{hom}^{\Upsilon }\left[ \mathcal{B},\mathcal{A}\right]
^{op}\in \mathsf{QLS}$, one for each pair $\mathcal{A}^{op},\mathcal{B}%
^{op}\in \mathsf{QLS}$, with analogous properties to $\underline{hom}\left[ 
\mathcal{B},\mathcal{A}\right] ^{op}=\underline{Hom}\left[ \mathcal{B}^{op},%
\mathcal{A}^{op}\right] $ (its internal Hom objects), but representing \emph{%
spaces of transformations }whose \emph{coordinate rings} $\underline{hom}%
^{\Upsilon }\left[ \mathcal{B},\mathcal{A}\right] $ and the ones of their
respective domains $\mathcal{B}^{op}$ do not commute among themselves. (The
work is mainly developed in a full subcategory of $\mathsf{QLS}$ whose
opposite objects are given by finitely generated graded algebras.) They give
rise to a $\mathsf{QLS}$-based category different from the one defined by
the function $\left( \mathcal{A},\mathcal{B}\right) \mapsto \underline{hom}%
\left[ \mathcal{B},\mathcal{A}\right] $. The mentioned non commutativity is
controlled by a collection of twisting maps $\tau _{\mathcal{A},\mathcal{B}}$
which define \emph{product spaces }$\left[ \underline{hom}^{\Upsilon }\left[ 
\mathcal{B},\mathcal{A}\right] \circ _{\tau }\mathcal{B}\right] ^{op}$
satisfying the inclusions $\underline{hom}^{\Upsilon }\left[ \mathcal{B},%
\mathcal{A}\right] \circ _{\tau }\mathcal{B}\subset \underline{hom}%
^{\Upsilon }\left[ \mathcal{B},\mathcal{A}\right] \otimes _{\tau }\mathcal{B}
$. We show that the (bi)algebras $\underline{end}^{\Upsilon }\left[ \mathcal{%
A}\right] \doteq \underline{hom}^{\Upsilon }\left[ \mathcal{A},\mathcal{A}%
\right] $, under certain circumstances, are 2-cocycle twisting of the
quantum semigroups\emph{\ }$\underline{end}\left[ \mathcal{A}\right] $. This
fact generalizes the twist equivalence (at a semigroup level) between, for
instance, the quantum groups $GL_{q}(n)$ and their multiparametric versions $%
GL_{q,\phi }(n)$.
\end{abstract}

\maketitle

\section*{Introduction}

Quantum linear spaces, or simply \emph{quantum spaces}, were defined by Y.
Manin as opposite (or dual) objects to finitely generated algebras ($\mathrm{%
FGA}$). He described the latter by pairs $\left( \mathbf{A}_{1},\mathbf{A}%
\right) $ where $\mathbf{A}$ is a unital algebra generated by some finite
linear subspace $\mathbf{A}_{1}\subset \mathbf{A}$; arrows between them are
unital algebra homomorphisms $\mathbf{A}\rightarrow \mathbf{B}$ which
restricted to $\mathbf{A}_{1}\subset \mathbf{A}$ give linear morphisms $%
\mathbf{A}_{1}\rightarrow \mathbf{B}_{1}$ ($\subset \mathbf{B}$). Every pair 
$\mathcal{A}\doteq \left( \mathbf{A}_{1},\mathbf{A}\right) $ is interpreted
as the (generically) \emph{noncommutative coordinate ring} of the quantum
space $\mathcal{A}^{op}$. In other words, if we call $\mathsf{QLS}$ the
category of quantum spaces and $\mathrm{FGA}$\textrm{\ }the category of
their coordinate rings (the pairs), then $\mathsf{QLS}\mathrm{\ }\doteq 
\mathrm{FGA}^{op}$. Since the duality between $\mathsf{QLS}$\textrm{\ }and $%
\mathrm{FGA}$, for quantum spaces\emph{\ }we understand the objects of each
one of these categories.

Initially, Manin \cite{man0} defined quantum spaces in terms of quadratic
algebras, which constitute a full subcategory \textrm{QA} of $\mathrm{FGA}$.
Later \cite{man1}, he extended the concept to arbitrary finitely generated
algebras. Then we refer to the objects of \textrm{QA} and\textrm{\ }$\mathrm{%
FGA}$ as\textrm{\ }\emph{quadratic} and \emph{general} quantum spaces,
respectively.

$\mathrm{FGA}$\textrm{\ }is a monoidal category $\left( \mathrm{FGA},\circ ,%
\mathcal{I}\right) $. Thus $\mathsf{QLS}$\textrm{\ }has, by duality, an
associated direct product of quantum spaces (of course, the term `direct' is
not used in the sense of `cartesian').

Manin's work is based on the existence of internal coHom objects in $\mathrm{%
FGA}$, which define the\textbf{\ }\emph{quantum matrix spaces }$\underline{%
hom}\left[ \mathcal{B},\mathcal{A}\right] $ with related coevaluation
morphisms $\mathcal{A}\rightarrow \underline{hom}\left[ \mathcal{B},\mathcal{%
A}\right] \circ \mathcal{B}$; or in their dual version, the internal Hom
objects 
\begin{equation*}
\underline{Hom}\left[ \mathcal{B}^{op},\mathcal{A}^{op}\right] =\underline{%
hom}\left[ \mathcal{B},\mathcal{A}\right] ^{op}
\end{equation*}
of the monoidal category $\left( \mathsf{QLS},\times ,\mathcal{I}%
^{op}\right) $, giving rise to the (noncommutative) algebraic geometric
version of quantum groups. More precisely, he defined quantum groups as the
universal Hopf envelope $aut\left[ \mathcal{A}\right] =GL\left[ \mathcal{A}%
\right] $ of the bialgebras $\underline{end}\left[ \mathcal{A}\right] =%
\underline{hom}\left[ \mathcal{A},\mathcal{A}\right] $. Restricted to 
\textrm{QA}, this approach is close to the FRT construction of quantum
matrix bialgebras \cite{frt}\cite{kassel}.

However the product $\circ $ is in essence the usual tensor product of
algebras. Following the work of \v{C}ap, Schichl and Van\v{z}ura \cite{cap},
we may ask why such a product is used in the construction above. Using it,
one is assuming that the coordinates of the space $\underline{Hom}\left[ 
\mathcal{B}^{op},\mathcal{A}^{op}\right] $ commute with the ones of $%
\mathcal{B}^{op}$, although the coordinates on the individual factors do not
commute among themselves. The main aim of this work is to carry out a twisted%
\emph{\ }version of the quantum matrix spaces, by means of replacing on the
coevaluation arrows the product $\circ $ by a twisted tensor product, which
we shall indicate $\circ _{\tau }$. It is worth remarking that we are not
going to change the monoidal structure on the category $\mathrm{FGA}$, which
would be to make the Manin construction in a Yang-Baxter or braided category 
\cite{joyal-street}. We just change the products between certain objects. In
other words, we study certain subclasses of maps $\mathcal{A}\rightarrow 
\mathcal{H}\circ _{\tau }\mathcal{B}$ for fixed $\mathcal{A}$ and $\mathcal{B%
}$ (each subclass representing a particular kind of non commutativity
between factors $\mathcal{H}^{op}$ and $\mathcal{B}^{op}$ of a quantum space
product $\left[ \mathcal{H}\circ _{\tau }\mathcal{B}\right] ^{op}$) from
which we construct universal objects, the \emph{twisted internal coHom
objects} $\underline{hom}^{\Upsilon }\left[ \mathcal{B},\mathcal{A}\right] $%
, with analogous properties to the proper coHom ones. In this way we extend,
and endow with an universal character, the results given in a previous work 
\cite{hr}.

We shall see that, in some cases, the associated twisted coEnd objects,
given by bialgebras $\underline{end}^{\Upsilon }\left[ \mathcal{A}\right] $,
correspond to twisting by counital 2-cocycles of the bialgebras $\underline{%
end}\left[ \mathcal{A}\right] $ in the sense of Drinfeld twisting process 
\cite{drin}. That is to say, $\underline{end}\left[ \mathcal{A}\right] $ and 
$\underline{end}^{\Upsilon }\left[ \mathcal{A}\right] $ are twist equivalent
as bialgebras.

\bigskip

The paper is organized as follows. In \S \textbf{1 }we analyze an
interesting full subcategory of $\mathrm{FGA}$ (formed out by finitely
generated graded algebras) on which our main work is based. This class of
objects will be called conic\emph{\ }quantum spaces, and they can be
regarded as an intermediate version between general and quadratic ones.
Before that, we review some well-known properties of the latter classes such
as monoidal structures and related internal (co)Hom objects.

We point out in \S \textbf{2 }a related semigroupoid structure associated to
certain comma categories of which the internal coHom objects are initial
objects. This structure is the main guide towards our construction. We show%
\textbf{\ }how to arrive at the notions of (co)evaluation, (co)identity
morphisms and (co)composition\emph{\ }of morphisms (in a compatible way)
using the universal character of the objects $\underline{hom}\left[ \mathcal{%
B},\mathcal{A}\right] $ and the above mentioned semigroupoid structure, and
then relate those notions to the concept of based categories. At this point
we introduce the necessary ingredients to pose in precise algebraic terms
the plan of our work and the main results.

In \S \textbf{3} we define a twisted tensor product $\circ _{\tau }$ between
objects of $\left( \mathrm{FGA},\circ \right) $ and morphisms between the
resulting objects. They are such that, as in the algebra case, their
isomorphism classes are in bijection to linear transformations with the same
properties as twisting maps\emph{\ }\cite{cap}. We introduce this concept in
the algebraic geometric terminology used by Manin, which allows us to
associate the (trivial or) non twisted product isomorphism class (i.e. the
class related to the standard product $\circ $, or to the flipping map) to 
\emph{commuting points} of factor spaces. Conversely, we are connecting the
idea of (generically) non commuting points with the concepts of (non
trivial) twisted tensor products and twisting maps.

Finally, we build up in \S \textbf{4} the twisted quantum matrix spaces $%
\underline{hom}^{\Upsilon }\left[ \mathcal{B},\mathcal{A}\right] $. They
give rise to a $\mathsf{QLS}$ or $\mathrm{FGA}^{op}$-based category
different from the one defined by the proper Hom objects, and with an
evaluation notion coming from the arrow (dual to) $\mathcal{A}\rightarrow 
\underline{hom}^{\Upsilon }\left[ \mathcal{B},\mathcal{A}\right] \circ
_{\tau }\mathcal{B}$. We also show the twist equivalence $\underline{end}%
^{\Upsilon }\left[ \mathcal{A}\right] \backsim \underline{end}\left[ 
\mathcal{A}\right] $, as bialgebras.

\bigskip

We often adopt definitions and notation extracted form Mac Lane's book \cite
{mac}. By $\Bbbk $ we indicate some of the numerics fields, $\mathbb{R}$ or $%
\mathbb{C}$. The usual tensor product on $\Bbbk \mathrm{-Alg}$\textrm{\ }and 
\textrm{Vct}$_{\Bbbk }$ (the categories of unital associative $\Bbbk $%
-algebras and of $\Bbbk $-vector spaces, respectively) is denoted by $%
\otimes $.

\section{Quantum linear spaces}

From its very definition, any pair $\left( \mathbf{A}_{1},\mathbf{A}\right)
\ $in $\mathrm{FGA}$\textrm{\ }has associated a canonical epimorphism of
unital algebras $\Pi :\mathbf{A}_{1}^{\otimes }\twoheadrightarrow \mathbf{A}$%
, where $\mathbf{A}_{1}^{\otimes }=\bigoplus_{n\in \mathbb{N}_{0}}\mathbf{A}%
_{1}^{\otimes n}$ is the tensor algebra of $\mathbf{A}_{1}$ ($\mathbb{N}_{0}$
denoting the non negative integers), such that restricted to $\mathbf{A}_{1}$
gives the inclusion $\mathbf{A}_{1}\hookrightarrow \mathbf{A}$, and it
defines a canonical isomorphism $\left. \mathbf{A}_{1}^{\otimes }\right/
\ker \Pi \backsimeq \mathbf{A}$. Moreover, the gradation in $\mathbf{A}%
_{1}^{\otimes }$ induces a filtration in $\mathbf{A}$ and $\ker \Pi $, with 
\begin{equation*}
\mathbf{A}=\bigcup\nolimits_{n\in \mathbb{N}_{0}}\mathbf{F}_{n};\qquad \;%
\mathbf{F}_{n}=\Pi \left( \bigoplus\nolimits_{i=0}^{n}\mathbf{A}%
_{1}^{\otimes i}\right)
\end{equation*}
and $\ker \Pi =\cup _{n\in \mathbb{N}_{0}}\mathbf{f}_{n};\;\mathbf{f}%
_{n}\subset \bigoplus_{i=0}^{n}\,\mathbf{A}_{1}^{\otimes i};\;\mathbf{f}%
_{0,1}=\left\{ 0\right\} $. On the other hand, for any pair of quantum
spaces $\left( \mathbf{A}_{1},\mathbf{A}\right) $ and $\left( \mathbf{B}_{1},%
\mathbf{B}\right) $, the arrows $\left( \mathbf{A}_{1},\mathbf{A}\right)
\rightarrow \left( \mathbf{B}_{1},\mathbf{B}\right) $ are given by algebra
homomorphisms $\alpha :\mathbf{A}\rightarrow \mathbf{B}$ completely defined
by the linear maps $\left. \alpha \right| _{\mathbf{A}_{1}}:\mathbf{A}%
_{1}\rightarrow \mathbf{B}_{1}$. In fact, they can be characterized by
arrows $\alpha _{1}:\mathbf{A}_{1}\rightarrow \mathbf{B}_{1}$ in \textrm{Vct}%
$_{\Bbbk }$ such that 
\begin{equation*}
\alpha _{1}^{\otimes }\left( \ker \Pi _{\mathbf{A}}\right) \subset \ker \Pi
_{\mathbf{B}},\;\;being\;\;\alpha _{1}^{\otimes }=\bigoplus\nolimits_{n\in 
\mathbb{N}_{0}}\alpha _{1}^{\otimes n}:\mathbf{A}_{1}^{\otimes }\rightarrow 
\mathbf{B}_{1}^{\otimes }
\end{equation*}
the unique extension of $\alpha _{1}$ to $\mathbf{A}_{1}^{\otimes }$ as a
morphism of algebras. We will say that $\alpha $ is the quotient map of $%
\alpha _{1}^{\otimes }$.

In these terms, the category \textrm{QA }of quadratic algebras, which
constitute a full subcategory of $\mathrm{FGA}$, is formed by pairs $\left( 
\mathbf{A}_{1},\mathbf{A}\right) $ such that the kernel of $\Pi :\mathbf{A}%
_{1}^{\otimes }\twoheadrightarrow \mathbf{A}$ is (a bilateral ideal)
algebraically generated by a subspace $\mathbf{I}\subset \mathbf{A}%
_{1}\otimes \mathbf{A}_{1}$, i.e. 
\begin{equation*}
\ker \Pi =\,\bigoplus\nolimits_{n\geq 2}\mathbf{I}_{n};\;\;\mathbf{I}%
_{n}=\,\sum\nolimits_{k=0}^{n-2}\mathbf{A}_{1}^{\otimes k}\otimes \mathbf{I}%
\otimes \mathbf{A}_{1}^{\otimes n-k-2}.
\end{equation*}
These objects are called quadratic quantum spaces. Its arrows can be
characterized by linear maps $\alpha _{1}$ such that $\alpha _{1}^{\otimes
2}\left( \mathbf{I}\right) \subset \mathbf{J}$, being $\mathbf{J}\subset 
\mathbf{B}_{1}^{\otimes 2}$ the generator subspace of $\ker \Pi _{\mathbf{B}%
} $.

\bigskip

Now let us present an intermediate version between quadratic and general
quantum spaces. We shall name conic algebras or conic quantum spaces those
pairs $\left( \mathbf{A}_{1},\mathbf{A}\right) $ of $\mathrm{FGA}$\textrm{\ }%
such that $\mathbf{A}$ is a graded algebra 
\begin{equation*}
\mathbf{A}=\,\bigoplus\nolimits_{n\in \mathbb{N}_{0}}\mathbf{A}_{n};\;%
\mathbf{A}_{n}=\Pi \left( \mathbf{A}_{1}^{\otimes n}\right) ,
\end{equation*}
or equivalently, $\ker \Pi =\bigoplus_{n\in \mathbb{N}_{0}}\mathbf{I}_{n}$, $%
\mathbf{I}_{n}\subset \mathbf{A}_{1}^{\otimes n}$ and $\mathbf{I}_{n}\otimes 
\mathbf{I}_{m}\subset \mathbf{I}_{n+m}$. The arrows between conic quantum
spaces are morphisms $\alpha _{1}:\mathbf{A}_{1}\rightarrow \mathbf{B}_{1}$
in \textrm{Vct}$_{\Bbbk }$\textrm{\ }such that $\alpha _{1}^{\otimes
n}\left( \mathbf{I}_{n}\right) \subset \mathbf{J}_{n}$. If we name $\mathrm{%
CA}$ the category of conic quantum spaces, the inclusion of full
subcategories $\mathrm{QA}\subset \mathrm{CA}\subset \mathrm{FGA}$ follows.

Examples of conic quantum spaces are, beside quadratics ones, those with
related ideals of the form $\ker \Pi =\mathbf{A}_{1}^{\otimes }\otimes 
\mathbf{I}\otimes \mathbf{A}_{1}^{\otimes }$, with $\mathbf{I}\subset 
\mathbf{A}_{1}^{\otimes m}$ for some $m\in \mathbb{N}$. We call them $m$%
\textbf{-}\emph{th quantum spaces}, and denote by $\mathrm{CA}^{m}$ the full
subcategory of $\mathrm{FGA}$ which has these pairs as objects. Thus \textrm{%
QA\ }$=\mathrm{CA}^{2}$.

Because our main constructions will be made on $\mathrm{CA}$, we are going
to study its associated internal coHom objects and some functorial
structures in terms of which they can be written. To this end, we first
refresh the monoidal (and other relevant) structures on $\mathrm{FGA}$%
\textrm{\ }and \textrm{QA}.

\subsection{Monoidal structure on $\mathrm{FGA}$ and its internal coHom
objects}

On $\mathrm{FGA}$ a bifunctor $\circ :\mathrm{FGA}\times \mathrm{FGA}%
\rightarrow \mathrm{FGA}$ can be defined, such that to every pair of quantum
spaces $\mathcal{A}=\left( \mathbf{A}_{1},\mathbf{A}\right) $ and $\mathcal{B%
}=\left( \mathbf{B}_{1},\mathbf{B}\right) $ it assigns the quantum space 
\begin{equation*}
\mathcal{A}\circ \mathcal{B}\doteq \left( \mathbf{A}_{1}\otimes \mathbf{B}%
_{1},\mathbf{A}\circ \mathbf{B}\right) ,
\end{equation*}
where $\mathbf{A}\circ \mathbf{B}$ is the subalgebra of $\mathbf{A}\otimes 
\mathbf{B}$ generated by $\mathbf{A}_{1}\otimes \mathbf{B}_{1}$. That means,
from the (canonical) isomorphism 
\begin{equation*}
\left[ \mathbf{A}_{1}\otimes \mathbf{B}_{1}\right] ^{\otimes }\backsimeq %
\left[ \mathbf{A}_{1}\otimes \mathbf{B}_{1}\right] ^{\widehat{\otimes }%
}\doteq \bigoplus\nolimits_{n\in \mathbb{N}_{0}}\mathbf{A}_{1}^{\otimes
n}\otimes \mathbf{B}_{1}^{\otimes n}
\end{equation*}
where $\left[ \mathbf{A}_{1}\otimes \mathbf{B}_{1}\right] ^{\widehat{\otimes 
}}$ is the subalgebra of $\mathbf{A}_{1}^{\otimes }\otimes \mathbf{B}%
_{1}^{\otimes }$ generated by $\mathbf{A}_{1}\otimes \mathbf{B}_{1}$,%
\footnote{%
The above isomorphism is built up from permutations; for instance,
restricted to $\left[ \mathbf{A}_{1}\otimes \mathbf{B}_{1}\right] ^{\otimes
2}$ is given by the map $S_{23}$ which acts as the flipping map on second
and thirst factors, and trivially on the others.} the kernel of the
canonical projection $\left[ \mathbf{A}_{1}\otimes \mathbf{B}_{1}\right]
^{\otimes }\twoheadrightarrow \mathbf{A}\circ \mathbf{B}$ is isomorphic to 
\begin{equation*}
\left( \mathbf{A}_{1}^{\otimes }\otimes \ker \Pi _{\mathbf{B}}+\ker \Pi _{%
\mathbf{A}}\otimes \mathbf{B}_{1}^{\otimes }\right) \cap \left[ \mathbf{A}%
_{1}\otimes \mathbf{B}_{1}\right] ^{\widehat{\otimes }}.
\end{equation*}
Note that $\mathbf{A}_{1}^{\otimes }\circ \mathbf{B}_{1}^{\otimes }=\left[ 
\mathbf{A}_{1}\otimes \mathbf{B}_{1}\right] ^{\widehat{\otimes }}$. We
frequently identify the latter algebra and $\left[ \mathbf{A}_{1}\otimes 
\mathbf{B}_{1}\right] ^{\otimes }$. Returning to the functor $\circ $, on
morphisms $\alpha :\mathcal{A}\rightarrow \mathcal{A}^{\prime }$ and $\beta :%
\mathcal{B}\rightarrow \mathcal{B}^{\prime }$, this functor gives 
\begin{equation*}
\alpha \circ \beta \doteq \left. \alpha \otimes \beta \right| _{\mathbf{A}%
\circ \mathbf{B}}:\mathcal{A}\circ \mathcal{B}\rightarrow \mathcal{A}%
^{\prime }\circ \mathcal{B}^{\prime },
\end{equation*}
the restriction of $\alpha \otimes \beta $ to the subalgebra $\mathbf{A}%
\circ \mathbf{B}\subset \mathbf{A}\otimes \mathbf{B}$. It is easy to see
that $\circ $ supplies $\mathrm{FGA}$ with an structure of symmetric
monoidal category, with unit object $\mathcal{I}\doteq \left( \Bbbk ,\Bbbk
\right) $. The functor $\circ $ also defines a monoid on the $m$-th's and
conic quantum spaces, but the unit object in these cases is $\mathcal{K}%
=\left( \Bbbk ,\Bbbk ^{\otimes }\right) =\left( \Bbbk ,\Bbbk \left[ e\right]
\right) $, i.e. the free algebra generated by the indeterminate $e$.

For every couple of objects $\mathcal{A}$ and $\mathcal{B}$ in $\left( 
\mathrm{FGA},\circ ,\mathcal{I}\right) \mathsf{\ }$(c.f. \cite{man1}) we
have a (left) internal coHom object $\underline{hom}\left[ \mathcal{B},%
\mathcal{A}\right] $, i.e. an initial object of the comma category $\left( 
\mathcal{A}\downarrow \mathrm{FGA}\circ \mathcal{B}\right) $. The objects of
each $\left( \mathcal{A}\downarrow \mathrm{FGA}\circ \mathcal{B}\right) $, 
\emph{diagrams }in the Manin terminology, are pairs 
\begin{equation*}
\left\langle \varphi ,\mathcal{H}\right\rangle _{_{\mathcal{A}\mathbf{,}%
\mathcal{B}}}=\left\langle \varphi ,\mathcal{H}\right\rangle ,\;\;with\;%
\mathcal{H}\in \mathrm{FGA},
\end{equation*}
and $\varphi $ a morphism $\mathcal{A}\rightarrow $ $\mathcal{H}\circ 
\mathcal{B}$; and its arrows $\left\langle \varphi ,\mathcal{H}\right\rangle
\rightarrow \left\langle \varphi ^{\prime },\mathcal{H}^{\prime
}\right\rangle $ are morphisms $\alpha :\mathcal{H}\rightarrow \mathcal{H}%
^{\prime }$ satisfying $\varphi ^{\prime }=\left( \alpha \circ I_{B}\right)
\,\varphi $. For every such comma category there exist an embedding 
\begin{equation*}
\frak{P}:\left( \mathcal{A}\downarrow \mathrm{FGA}\circ \mathcal{B}\right)
\hookrightarrow \mathrm{FGA}\;\;/\;\;\frak{P}\left\langle \varphi ,\mathcal{H%
}\right\rangle =\mathcal{H}.
\end{equation*}

If $\mathcal{A}$ and $\mathcal{B}$ are generated by $\mathbf{A}_{1},\mathbf{B%
}_{1}\in \mathrm{Vct}_{\Bbbk }$ with $\dim \mathbf{A}_{1}=\frak{n}$ and $%
\dim \mathbf{B}_{1}=\frak{m}$, the initial object of $\left( \mathcal{A}%
\downarrow \mathrm{FGA}\circ \mathcal{B}\right) $ is the pair $\left\langle
\delta _{_{\mathcal{A},\mathcal{B}}},\underline{hom}\left[ \mathcal{B},%
\mathcal{A}\right] \right\rangle $ where $\underline{hom}\left[ \mathcal{B},%
\mathcal{A}\right] $ is given by an algebra $\left. \mathbf{E}_{1}^{\otimes
}\right/ \mathbf{K}$ being 
\begin{equation}
\mathbf{E}_{1}\doteq span\left[ z_{i}^{j}\right] _{i,j=1}^{\frak{n},\frak{m}}
\label{E}
\end{equation}
and $\mathbf{K}$ a bilateral ideal depending on $\ker \Pi _{\mathbf{A}}$ and 
$\ker \Pi _{\mathbf{B}}$, and $\delta _{_{\mathcal{A},\mathcal{B}}}=\delta $
is the coevaluation map 
\begin{equation}
\delta :\mathcal{A}\rightarrow \underline{hom}\left[ \mathcal{B},\mathcal{A}%
\right] \circ \mathcal{B}\;\;\;/\;\;\;a_{i}\mapsto z_{i}^{j}\otimes b_{j},
\label{eva}
\end{equation}
with $\left\{ a_{i}\right\} $ and $\left\{ b_{j}\right\} $ basis of $\mathbf{%
A}_{1}$ and $\mathbf{B}_{1}$. Generically each object $\underline{hom}\left[ 
\mathcal{B},\mathcal{A}\right] $ is equal to $\left( \mathbf{E}_{1},\left. 
\mathbf{E}_{1}^{\otimes }\right/ \mathbf{K}\right) $, in the sense that
symbols $z_{i}^{j}$ are linearly independent.\footnote{%
Nevertheless, suppose a relation $c^{i}a_{i}+\lambda 1$ ($\lambda \neq 0$)
is present in a quantum space $\mathcal{A}$, such that $c^{k}=0$, for some $%
k $. Then $\underline{end}\left[ \mathcal{A}\right] $ is generated by
elements $z_{i}^{j}$ related, among other possible constraints, by the
linear form $c^{i}z_{i}^{k}$.} From formal properties of internal (co)Hom
objects we have the cocomposition and coidentity maps 
\begin{eqnarray}
\underline{hom}\left[ \mathcal{B},\mathcal{A}\right] &\rightarrow &%
\underline{hom}\left[ \mathcal{C},\mathcal{A}\right] \circ \underline{hom}%
\left[ \mathcal{B},\mathcal{C}\right] ,  \notag \\
&&  \label{coun} \\
\underline{end}\left[ \mathcal{A}\right] &\doteq &\underline{hom}\left[ 
\mathcal{A},\mathcal{A}\right] \rightarrow \mathcal{I},  \notag
\end{eqnarray}
which on each coEnd objects $\underline{end}\left[ \mathcal{A}\right] \in 
\mathrm{FGA}$ define the (injective) comultiplications 
\begin{equation*}
\Delta :\underline{end}\left[ \mathcal{A}\right] \hookrightarrow \underline{%
end}\left[ \mathcal{A}\right] \circ \underline{end}\left[ \mathcal{A}\right]
\;\;\;/\;\;\;z_{i}^{j}\mapsto z_{i}^{k}\otimes z_{k}^{j}
\end{equation*}
and (in our case surjective) counits 
\begin{equation*}
\varepsilon :\underline{end}\left[ \mathcal{A}\right] \twoheadrightarrow 
\mathcal{I}\;\;/\;\;z_{i}^{j}\mapsto \delta _{i}^{j}\in \Bbbk ,
\end{equation*}
giving $\underline{end}\left[ \mathcal{A}\right] $ a bialgebra structure $%
\left( \Delta ,\varepsilon \right) $. In addition, the (injective)
coevaluation $\mathcal{A}\hookrightarrow \underline{end}\left[ \mathcal{A}%
\right] \circ \mathcal{A}$ defines a coaction which supplies $\mathcal{A}$
with a structure of left $\underline{end}\left[ \mathcal{A}\right] $%
-comodule algebra.

\bigskip

\textbf{The functor }$\frak{F}$\textbf{:} The following observation will be
crucial to construct our twisted quantum matrix spaces. For any $%
\left\langle \varphi ,\mathcal{H}\right\rangle $ in $\left( \mathcal{A}%
\downarrow \mathrm{FGA}\circ \mathcal{B}\right) $, with $\mathcal{H}=\left( 
\mathbf{H}_{1},\mathbf{H}\right) $, there exist an associated linear space 
\begin{equation*}
\mathbf{H}_{1}^{\varphi }\doteq \,\pounds \left\langle
h_{i}^{j}\right\rangle _{i,j=1}^{\frak{n},\frak{m}}\subset \mathbf{H}_{1}
\end{equation*}
defined by the restriction of $\varphi $ to $\mathbf{A}_{1}$, i.e. by a
linear map $a_{i}\mapsto h_{i}^{j}\otimes b_{j}$. In particular, there is a
linear surjection 
\begin{equation*}
\pi ^{\varphi }:\mathbf{B}_{1}^{\ast }\otimes \mathbf{A}_{1}%
\twoheadrightarrow \mathbf{H}_{1}^{\varphi }\;\;:\;\;b^{j}\otimes
a_{i}\mapsto h_{i}^{j},
\end{equation*}
being $\left\{ b^{i}\right\} \subset \mathbf{B}_{1}^{\ast }$ the dual basis
of $\left\{ b_{i}\right\} $. In what follows we identify the linear spaces $%
Lin\left[ \mathbf{B}_{1},\mathbf{A}_{1}\right] $ and $\mathbf{B}_{1}^{\ast
}\otimes \mathbf{A}_{1}$. Moreover, $\left\langle \varphi ,\mathcal{H}%
\right\rangle $ can be related to a subalgebra $\mathbf{H}^{\varphi }\ $of $%
\mathbf{H}$ generated by $\mathbf{H}_{1}^{\varphi }$, and the corresponding
quantum space $\mathcal{H}^{\varphi }=\left( \mathbf{H}_{1}^{\varphi },%
\mathbf{H}^{\varphi }\right) $. For the initial object, we have generically
that $\mathbf{E}_{1}=\mathbf{E}_{1}^{\delta }$ and $\pi ^{\delta }$ gives an
isomorphism $Lin\left[ \mathbf{B}_{1},\mathbf{A}_{1}\right] \backsimeq 
\mathbf{E}_{1}$.

In this way, for every couple $\mathcal{A},\mathcal{B}\in \mathrm{FGA}$, the
map 
\begin{equation*}
\left\langle \varphi ,\mathcal{H}\right\rangle \in \left( \mathcal{A}%
\downarrow \mathrm{FGA}\circ \mathcal{B}\right) \longmapsto \mathcal{H}%
^{\varphi }\in \mathrm{FGA},
\end{equation*}
extended to any arrow $\left\langle \varphi ,\mathcal{H}\right\rangle 
\overset{\alpha }{\rightarrow }\left\langle \psi ,\mathcal{G}\right\rangle $
as $\alpha \mapsto \left. \alpha \right| _{\mathbf{H}^{\varphi }}$, provides
a collection of functors 
\begin{equation*}
\frak{F}:\left( \mathcal{A}\downarrow \mathrm{FGA}\circ \mathcal{B}\right)
\rightarrow \mathrm{FGA}.
\end{equation*}
For any $\alpha :\left\langle \varphi ,\mathcal{H}\right\rangle \rightarrow
\left\langle \psi ,\mathcal{G}\right\rangle $, $\frak{F}\alpha =\left.
\alpha \right| _{\mathbf{H}^{\varphi }}:\mathbf{H}^{\varphi }\rightarrow 
\mathbf{G}^{\psi }$ defines the epi arrow 
\begin{equation}
\frak{F}\alpha :\frak{F}\left\langle \varphi ,\mathcal{H}\right\rangle
\twoheadrightarrow \frak{F}\left\langle \psi ,\mathcal{G}\right\rangle .
\label{epig}
\end{equation}
Obviously, $\frak{F}=\frak{P}$ for all $\left\langle \varphi ,\mathcal{H}%
\right\rangle $ such that $\mathcal{H}^{\varphi }=\mathcal{H}$; in
particular, $\frak{F}=\frak{P}$ for generic initial objects.

\bigskip

\textbf{The geometric role of initial objects: }We will make a brief comment
on the algebraic geometric interpretation of the relationship between
diagrams $\left\langle \varphi ,\mathcal{H}\right\rangle $ and internal
coHom objects in $\mathrm{FGA}$.

The initial property of any $\underline{hom}\left[ \mathcal{B},\mathcal{A}%
\right] $ means that for a given object $\mathcal{H}$ and a morphism $%
\varphi :\mathcal{A}\rightarrow \mathcal{H}\circ \mathcal{B}$, there exist a
unique morphism of quantum spaces $\alpha :\underline{hom}\left[ \mathcal{B},%
\mathcal{A}\right] \rightarrow \mathcal{H}$ making commutative the diagram 
\begin{equation}
\begin{diagram} & & \mathcal{A} & & \\ & \ldTo^{\delta
_{\mathcal{A},\mathcal{B}}} & & \rdTo^{\varphi } & \\ \underline{hom}\left[
\mathcal{B},\mathcal{A}\right] \circ \mathcal{B}& & \rTo^{\alpha \circ
I_{\mathcal{B}}}& &\mathcal{H}\circ \mathcal{B} \\ \end{diagram}
\label{diaun}
\end{equation}
From Eq. $\left( \ref{epig}\right) $, we have the epi\textbf{\ }(or
epimorphism, for the underlying algebras) 
\begin{equation}
\underline{hom}\left[ \mathcal{B},\mathcal{A}\right] \twoheadrightarrow 
\frak{F}\left\langle \varphi ,\mathcal{H}\right\rangle  \label{epi}
\end{equation}
for every $\left\langle \varphi ,\mathcal{H}\right\rangle \in \left( 
\mathcal{A}\downarrow \mathrm{FGA}\circ \mathcal{B}\right) $. Thus, the
opposite of any $\frak{F}\left\langle \varphi ,\mathcal{H}\right\rangle $
can be regarded as a \emph{subspace} 
\begin{equation*}
\frak{F}\left\langle \varphi ,\mathcal{H}\right\rangle ^{op}\hookrightarrow 
\underline{hom}\left[ \mathcal{B},\mathcal{A}\right] ^{op}=\underline{Hom}%
\left[ \mathcal{B}^{op},\mathcal{A}^{op}\right]
\end{equation*}
of the space of quantum linear maps from $\mathcal{B}$ to $\mathcal{A}$.
More precisely, the monic\textbf{\ }$\frak{F}\left\langle \varphi ,\mathcal{H%
}\right\rangle ^{op}\hookrightarrow \underline{hom}\left[ \mathcal{B},%
\mathcal{A}\right] ^{op}$ gives the representative of an equivalence class
of monics defining a subobject of\textbf{\ }$\underline{hom}\left[ \mathcal{B%
},\mathcal{A}\right] ^{op}$.

We could say $\underline{hom}\left[ \mathcal{B},\mathcal{A}\right] $ is a 
\textbf{`}\emph{generic point}' of the `\emph{noncommutative algebraic
variety}' $\underline{Hom}\left[ \mathcal{B}^{op},\mathcal{A}^{op}\right] $
(or that the quantum space $\underline{Hom}\left[ \mathcal{B}^{op},\mathcal{A%
}^{op}\right] $ is the `\emph{locus of\ }$\underline{hom}\left[ \mathcal{B},%
\mathcal{A}\right] $'), and regard the object $\frak{F}\left\langle \varphi ,%
\mathcal{H}\right\rangle =\mathcal{H}^{\varphi }$ as a `\emph{specialization}%
' of it.

\subsection{Quadratic quantum spaces}

There is another symmetric monoidal structure on $\mathrm{QA}$ given by the
bifunctor $\bullet :\mathrm{QA}\times \mathrm{QA}\longrightarrow \mathrm{QA}$
such that to every pair of quadratic quantum spaces $\mathcal{A}=\left( 
\mathbf{A}_{1},\mathbf{A}\right) $ and $\mathcal{B}=\left( \mathbf{B}_{1},%
\mathbf{B}\right) $, it assigns the quantum space 
\begin{equation*}
\mathcal{A}\bullet \mathcal{B}\doteq \left( \mathbf{A}_{1}\otimes \mathbf{B}%
_{1},\mathbf{A}\bullet \mathbf{B}\right) ;\;\;\mathbf{A}\bullet \mathbf{B}%
\doteq \left. \left[ \mathbf{A}_{1}\otimes \mathbf{B}_{1}\right] ^{\widehat{%
\otimes }}\right/ I\left[ \mathbf{I}\otimes \mathbf{J}\right] ,
\end{equation*}
where $I\left[ \mathbf{X}\right] $ means the bilateral ideal generated by
the set $\mathbf{X}$.\footnote{%
Manin uses $\left[ \mathbf{A}_{1}\otimes \mathbf{B}_{1}\right] ^{\otimes }$
instead of $\left[ \mathbf{A}_{1}\otimes \mathbf{B}_{1}\right] ^{\widehat{%
\otimes }}$. We are working with the latter in order to deal with the
homogeneous case bellow.} For the morphisms $\mathcal{A}\overset{\alpha }{%
\rightarrow }\mathcal{A}^{\prime }$ and $\mathcal{B}\overset{\beta }{%
\rightarrow }\mathcal{B}^{\prime }$, the morphism $\alpha \bullet \beta :%
\mathcal{A}\bullet \mathcal{B}\rightarrow $ $\mathcal{A}^{\prime }\bullet 
\mathcal{B}^{\prime }$ is defined as the quotient map of 
\begin{equation*}
\left. \alpha _{1}^{\otimes }\otimes \beta _{1}^{\otimes }\right| _{\left[ 
\mathbf{A}_{1}\otimes \mathbf{B}_{1}\right] ^{\widehat{\otimes }}}=\left[
\alpha _{1}\otimes \beta _{1}\right] ^{\widehat{\otimes }}=\bigoplus%
\nolimits_{n\in \mathbb{N}_{0}}\alpha _{1}^{\otimes n}\otimes \beta
_{1}^{\otimes n},
\end{equation*}
which is well-defined because $\alpha _{1}^{\otimes 2}\otimes \beta
_{1}^{\otimes 2}\left( \mathbf{I}\otimes \mathbf{J}\right) \subset \mathbf{I}%
^{\prime }\otimes \mathbf{J}^{\prime }$. The unit object is $\mathcal{U}%
=\left( \Bbbk ,\mathbf{U}\right) $, being $\mathbf{U}=\left. \Bbbk \left[ e%
\right] \right/ I\left[ e^{2}\right] $.

Now, consider the covariant functor $!:\mathrm{QA}^{op}\rightarrow \mathrm{QA%
}$ such that 
\begin{equation*}
\mathcal{A}^{!}=\left( \mathbf{A}_{1}^{\ast },\mathbf{A}^{!}\right) ,\qquad
\;\mathbf{A}^{!}\doteq \left. \mathbf{A}_{1}^{\ast \otimes }\right/ I\left[ 
\mathbf{I}^{\perp }\right] ;
\end{equation*}
where $\mathbf{I}^{\perp }\doteq \left\{ r\in \mathbf{A}_{1}^{\ast \otimes
2}:\left\langle r,q\right\rangle =0,\;\forall q\in \mathbf{I}\right\} $ is
the annihilator of $\mathbf{I}$ in relation with the standard pairing. For a
morphism $\alpha :\mathcal{A}\rightarrow \mathcal{B}$ it assigns $\alpha
^{!}:\mathcal{B}^{!}\rightarrow \mathcal{A}^{!}$, being $\alpha ^{!}$ the
quotient map\footnote{%
Remember that $\alpha _{1}^{\ast }:\mathbf{B}_{1}^{\ast }\rightarrow \mathbf{%
A}_{1}^{\ast }$ is given by $\left\langle \alpha _{1}^{\ast }\left( b\right)
,a\right\rangle =\left\langle b,\alpha _{1}\left( a\right) \right\rangle $
(calculated on the respective pairings) for all $a\in \mathbf{A}_{1}^{\ast }$
and $b\in \mathbf{B}_{1}^{\ast }$. Then $\alpha _{1}^{\ast \otimes 2}\left( 
\mathbf{J}^{\perp }\right) \subset \mathbf{I}^{\perp }$ follows immediately.}
of $\alpha _{1}^{\ast \otimes }:\mathbf{B}_{1}^{\ast \otimes }\rightarrow 
\mathbf{A}_{1}^{\ast \otimes }$. The main relationships among the functors $%
\circ $, $\bullet $ and $!$ can be summarized in the following equations 
\begin{equation}
\mathcal{A}^{!!}\backsimeq \mathcal{A},\;\;\;\left( \mathcal{A}\circ 
\mathcal{B}\right) ^{!}\backsimeq \mathcal{A}^{!}\bullet \mathcal{B}%
^{!},\;\;\;\left( \mathcal{A}\bullet \mathcal{B}\right) ^{!}\backsimeq 
\mathcal{A}^{!}\circ \mathcal{B}^{!},\;\;\;\mathcal{K}^{!}\backsimeq 
\mathcal{U}.  \label{p}
\end{equation}
Manin showed that $\left( \mathrm{QA},\circ ,\mathcal{K}\right) $ has
internal coHom objects given by $\underline{hom}\left[ \mathcal{B},\mathcal{A%
}\right] =\mathcal{B}^{!}\bullet \mathcal{A}$, with \textbf{(}co\textbf{)}%
evaluations, compositions and (co)identities (replacing $\mathcal{I}$ by $%
\mathcal{K}$) given by $\left( \ref{eva}\right) $ and $\left( \ref{coun}%
\right) $.\ $\mathrm{QA}$ has also internal Hom objects in relation with the
monoidal structure $\bullet $, that is, for $\left( \mathrm{QA},\bullet ,%
\mathcal{U}\right) $. They can be defined as 
\begin{equation*}
\underline{Hom}\left[ \mathcal{B},\mathcal{A}\right] =\mathcal{B}^{!}\circ 
\mathcal{A}\backsimeq \underline{hom}\left[ \mathcal{B}^{!},\mathcal{A}^{!}%
\right] ^{!}.
\end{equation*}

We shall be mainly concerned with coHom objects, because they define on $%
\mathsf{QLS}=\mathrm{FGA}^{op}$ the spaces of quantum linear maps 
\begin{equation*}
\underline{Hom}\left[ \mathcal{B}^{op},\mathcal{A}^{op}\right] =\underline{%
hom}\left[ \mathcal{B},\mathcal{A}\right] ^{op}\in \mathsf{QLS},
\end{equation*}
or the \emph{noncommutative algebraic varieties} given by the \emph{locus of 
}$\underline{hom}\left[ \mathcal{B},\mathcal{A}\right] $.

The structures described in this section can be defined also on $m$-th
quantum spaces (we just have to change $2$ by $m$ in all definitions). For
instance, a unit element for $\left( \mathrm{CA}^{m},\bullet \right) $ would
be $\mathcal{U}_{m}=\left( \Bbbk ,\mathbf{U}_{m}\right) $, being $\mathbf{U}%
_{m}=\Bbbk \left[ e\right] /I\left[ e^{m}\right] $. Note that $\mathcal{U}%
_{2}=\mathcal{U}$, and that we can define $\mathcal{U}_{\infty }\doteq 
\mathcal{K}$.

\subsection{The conic case}

From now on, we understand by quantum spaces the conic ones, unless we say
the contrary. If $\left( \mathbf{A}_{1},\mathbf{A}\right) $ and $\left( 
\mathbf{B}_{1},\mathbf{B}\right) $ are objects of this class, then 
\begin{equation*}
\ker \Pi _{\mathbf{A}}=\,\bigoplus\nolimits_{n\in \mathbb{N}_{0}}\mathbf{I}%
_{n},\;\;\;\;\ker \Pi _{\mathbf{B}}=\bigoplus\nolimits_{n\in \mathbb{N}%
_{0}}\,\mathbf{J}_{n},
\end{equation*}
Define the functors $\odot $ and $!$ as 
\begin{equation*}
\mathbf{A}\odot \mathbf{B}\doteq \left. \left[ \mathbf{A}_{1}\otimes \mathbf{%
B}_{1}\right] ^{\widehat{\otimes }}\right/ \bigoplus\nolimits_{n\in \mathbb{N%
}_{0}}\mathbf{I}_{n}\otimes \mathbf{J}_{n}
\end{equation*}
and $\mathbf{A}^{!}\doteq \left. \mathbf{A}_{1}^{\ast \otimes }\right/ I%
\left[ \bigoplus\nolimits_{n\geq 2}\mathbf{I}_{n}^{\perp }\right] $, where
now 
\begin{equation*}
\mathbf{I}_{n}^{\perp }=\left\{ x\in \mathbf{A}_{1}^{\ast \otimes
n}:\left\langle x,y\right\rangle =0,\;\forall y\in \mathbf{I}_{n}\right\}
,\;\;n\geq 2,\;\;\;\mathbf{I}_{0,1}^{\perp }=\left\{ 0\right\} .
\end{equation*}
A unit element for $\left( \mathrm{CA},\odot \right) \ $is $\mathcal{U}$ as
for $\mathrm{QA}$, however the restriction of $\odot $ to every $\mathrm{CA}%
^{m}$ (in particular $\mathrm{QA}$) does not coincide with the analogous
functor $\bullet $ defined on these subcategories. On the other hand, the
functor $!$ does coincide with the corresponding to the $m$-th cases because 
\begin{equation}
\left( \mathbf{A}_{1}^{\otimes r}\otimes \mathbf{I}\otimes \mathbf{A}%
_{1}^{\otimes s}\right) ^{\perp }=\mathbf{A}_{1}^{\ast \otimes r}\otimes 
\mathbf{I}^{\perp }\otimes \mathbf{A}_{1}^{\ast \otimes s},\;\forall r,s\in 
\mathbb{N}_{0},  \label{k}
\end{equation}
where $\mathbf{I}\subset \mathbf{A}_{1}^{\otimes m}$ is the generator of $%
\ker \Pi $; but $!^{2}\ncong id_{\mathrm{CA}}$. The preserved properties are
(w.r.t. Eq. $\left( \ref{p}\right) $), 
\begin{equation*}
\left( \mathcal{A}\odot \mathcal{B}\right) ^{!}\backsimeq \mathcal{A}%
^{!}\circ \mathcal{B}^{!},\qquad \mathcal{K}^{!}\backsimeq \mathcal{U}%
,\qquad \;\mathcal{K}\backsimeq \mathcal{U}^{!},
\end{equation*}
while $\left( \mathcal{A}\circ \mathcal{B}\right) ^{!}\ncong \mathcal{A}%
^{!}\odot \mathcal{B}^{!}$ and $\mathcal{A}^{!!}\ncong \mathcal{A}$.

Although the functors $\odot $ and $!$ do not have the nice properties of $%
\bullet $ and $!$, they lead us to define on $\mathrm{CA}$ some kind of
covariant mixing\emph{\ }of them, namely 
\begin{equation*}
\begin{array}{c}
\triangleright :\mathrm{CA}^{op}\times \mathrm{CA}\rightarrow \mathrm{CA}%
,\;\triangleleft :\mathrm{CA}\times \mathrm{CA}^{op}\rightarrow \mathrm{CA},
\\ 
\\ 
\diamond \,:\mathrm{CA}^{op}\times \mathrm{CA}^{op}\rightarrow \mathrm{CA},
\end{array}
\end{equation*}
given on objects by (identifying each object with its opposite)\footnote{%
The maps on morphisms are immediate; for instance, we have for $%
\triangleright $ that $\alpha ^{op}\times \beta $ is sent to the quotient
map of $\left[ \alpha _{1}^{\ast }\otimes \beta _{1}\right] ^{\widehat{%
\otimes }}$.} 
\begin{eqnarray*}
\mathbf{A}\triangleright \mathbf{B} &\doteq &\left. \left[ \mathbf{A}%
_{1}^{\ast }\otimes \mathbf{B}_{1}\right] ^{\widehat{\otimes }}\right/ I%
\left[ \bigoplus\nolimits_{n\in \mathbb{N}_{0}}\mathbf{I}_{n}^{\perp
}\otimes \mathbf{J}_{n}\right] \qquad \\
\mathbf{A}\triangleleft \mathbf{B} &\doteq &\left. \left[ \mathbf{A}%
_{1}\otimes \mathbf{B}_{1}^{\ast }\right] ^{\widehat{\otimes }}\right/ I%
\left[ \bigoplus\nolimits_{n\in \mathbb{N}_{0}}\mathbf{I}_{n}\otimes \mathbf{%
J}_{n}^{\perp }\right] \\
\mathbf{A}\diamond \mathbf{B} &\doteq &\left. \left[ \mathbf{A}_{1}^{\ast
}\otimes \mathbf{B}_{1}^{\ast }\right] ^{\widehat{\otimes }}\right/ I\left[
\bigoplus\nolimits_{n\in \mathbb{N}_{0}}\mathbf{I}_{n}^{\perp }\otimes 
\mathbf{J}_{n}^{\perp }\right] .
\end{eqnarray*}
Using $\left( \ref{k}\right) $, it can be shown that restricted to every $%
\mathrm{CA}^{m}$, 
\begin{equation}
\triangleright \,=\bullet \,\left( !\times id\right) ,\;\triangleleft
\,=\bullet \,\left( id\times !\right) ,\;\diamond \,=\bullet \,\left(
!\times !\right) .  \label{compa}
\end{equation}

Associated to $\triangleright $ and $\triangleleft $ we have $\mathcal{K}$
as left and as a right unit, respectively, in the sense that $\mathcal{K}%
\triangleright \mathcal{A}\backsimeq \mathcal{A}$ and $\mathcal{A}%
\triangleleft \mathcal{K}\backsimeq \mathcal{A}$ for any $\mathcal{A}$ in $%
\mathrm{CA}$ (while $\mathcal{A}\triangleright \mathcal{U}\backsimeq 
\mathcal{A}^{!}$, $\mathcal{U}\triangleleft \mathcal{A}\backsimeq \mathcal{A}%
^{!}$ and $\mathcal{K}\diamond \mathcal{A}\backsimeq \mathcal{A}\diamond 
\mathcal{K}\backsimeq \mathcal{A}^{!}$).

\begin{theorem}
The monoidal category $\left( \mathrm{CA},\circ ,\mathcal{K}\right) $ has 
\textbf{internal coHom objects} given by 
\begin{equation*}
\underline{hom}\left[ \mathcal{B},\mathcal{A}\right] =\mathcal{B}%
\triangleright \mathcal{A}.
\end{equation*}
In other words, the functor $\triangleright $ is an internal coHom functor.\
\ \ $\blacksquare $
\end{theorem}

Before passing to the proof, we introduce some notation and make some
observations useful in next sections. Consider $\mathcal{A}=\left( \mathbf{A}%
_{1},\mathbf{A}\right) $ and $\mathcal{B}=\left( \mathbf{B}_{1},\mathbf{B}%
\right) $, with associated graded ideals linearly generated by the relations 
\begin{equation}
\left\{ R_{\lambda _{n}}^{k_{1}...k_{n}}\;a_{k_{1}}...a_{k_{n}}\right\}
_{\lambda _{n}\in \Lambda _{n}}\subset \mathbf{I}_{n},\;\;\;\left\{ S_{\mu
_{n}}^{k_{1}...k_{n}}\;b_{k_{1}}...b_{k_{n}}\right\} _{\mu _{n}\in \Phi
_{n}}\subset \mathbf{J}_{n},  \label{ideals}
\end{equation}
(sum over repeated indices are assumed) being $\left\{ a_{i}\right\} $ and $%
\left\{ b_{i}\right\} $ basis of $\mathbf{A}_{1}$ and $\mathbf{B}_{1}$. For
each $\mathbf{J}_{n}$ consider some complement $\mathbf{J}_{n}^{c}$, such
that $\mathbf{J}_{n}\oplus \mathbf{J}_{n}^{c}=\mathbf{B}_{1}^{\otimes n}$,
and indicate by $\left\{ J_{\omega _{n}}^{c}\right\} _{\omega _{n}\in \Omega
_{n}}$ a basis for $\mathbf{J}_{n}^{c}$. Then we can write 
\begin{equation}
b_{k_{1}}...b_{k_{n}}=\left( S^{\bot }\right) _{k_{1}...k_{n}}^{\omega
_{n}}\;J_{\omega _{n}}^{c}+J_{k_{1}...k_{n}};\;J_{k_{1}...k_{n}}\in \mathbf{J%
}_{n}.  \label{deco}
\end{equation}
It is important to note that each $\mathbf{J}_{n}^{\perp }\subset \mathbf{B}%
_{1}^{\ast \otimes n}$ is precisely the linear space spanned by 
\begin{equation}
\left\{ b^{k_{1}}...b^{k_{n}}\;\left( S^{\bot }\right)
_{k_{1}...k_{n}}^{\omega _{n}}\right\} _{\omega _{n}\in \Omega _{n}},
\label{duu}
\end{equation}
being $\left\{ b^{i}\right\} \subset \mathbf{B}_{1}^{\ast }$ the dual basis
of $\left\{ b_{i}\right\} $. Those subspaces define the ideal associated to $%
\mathcal{B}^{!}$.

\begin{proof}
An internal coHom object $\underline{hom}\left[ \mathcal{B},\mathcal{A}%
\right] $ of $\mathrm{CA}$, if there exists, is an initial object of the
comma category $\left( \mathcal{A}\downarrow \mathrm{CA}\circ \mathcal{B}%
\right) $. So, let us show that there exists an initial diagram $%
\left\langle \delta ,\mathcal{E}\right\rangle $ in every such a category.

Let $\mathbf{E}_{1}$ be as in Eq. $\left( \ref{E}\right) $, and consider the
map $\delta _{1}:a_{i}\mapsto z_{i}^{j}\otimes b_{j}$. Then, $\delta
_{1}^{\otimes }$ evaluated on an element $R_{\lambda
_{n}}^{k_{1}...k_{n}}\;a_{k_{1}}...a_{k_{n}}\in \mathbf{I}_{n}$ gives (using 
$\left( \ref{deco}\right) $) 
\begin{equation*}
R_{\lambda _{n}}^{k_{1}...k_{n}}\;z_{k_{1}}^{j_{1}}\cdot
z_{k_{2}}^{j_{2}}\;...\;z_{k_{n}}^{j_{n}}\otimes \left( \left( S^{\bot
}\right) _{j_{1}...j_{n}}^{\omega _{n}}J_{\omega
_{n}}^{c}+J_{k_{1}...k_{n}}\right) .
\end{equation*}
Defining $\mathbf{E}\doteq \left. \mathbf{E}_{1}^{\otimes }\right/ \mathbf{K}
$, with $\mathbf{K}=\bigoplus_{n}\mathbf{K}_{n}$ the ideal algebraically
generated by 
\begin{equation}
\left\{ R_{\lambda _{n}}^{k_{1}...k_{n}}\;z_{k_{1}}^{j_{1}}\cdot
z_{k_{2}}^{j_{2}}\;...\;z_{k_{n}}^{j_{n}}\;\left( S^{\bot }\right)
_{j_{1}...j_{n}}^{\omega _{n}}\right\} _{\omega _{n}\in \Omega
_{n}}^{\lambda _{n}\in \Lambda _{n}}\subset \mathbf{K}_{n},  \label{ideal}
\end{equation}
$\delta _{1}$ can be extended to an homomorphism $\delta :\mathbf{A}%
\rightarrow \mathbf{E}\otimes \mathbf{B}$, resulting the latter a quotient
map of $\delta _{1}^{\otimes }$. Thus, the pair $\mathcal{E}=\left( \mathbf{E%
}_{1},\mathbf{E}\right) $ defines the diagram $\left\langle \delta ,\mathcal{%
E}\right\rangle $. Its initial character is obvious. Moreover, from the
comment leading to Eq. $\left( \ref{duu}\right) $, it is immediate that the
linear map $z_{i}^{j}\mapsto b^{j}\otimes a_{i}$ extends to an algebra
isomorphism $\mathbf{E}\backsimeq \mathbf{B}\triangleright \mathbf{A}$.
Then, we can define $\underline{hom}\left[ \mathcal{B},\mathcal{A}\right] =%
\mathcal{B}\triangleright \mathcal{A}$ ($\backsimeq \mathcal{E}$).
\end{proof}

Nevertheless, the existence of internal Hom objects for $\left( \mathrm{CA}%
,\odot ,\mathcal{U}\right) $ can not be established, as we have made for
each $\left( \mathrm{CA}^{m},\bullet ,\mathcal{U}_{m}\right) $, with $m\geq
2 $.

\section{Semigroupoids and (co)based categories}

The main consequences related to the existence of internal (co)Hom objects
in a given category are the notions of (co)evaluation or (co)action,
(co)composition or (co)multiplication, and (co)identity or (co)unit, which
enable us to think of them as the `(co)spaces of morphisms' between the
objects in the corresponding category. We will first see how such notions
arise in $\mathrm{FGA}$. Then, we will encode them in the framework of
(co)based categories, since in this scheme our twisted version of quantum
matrices is developed.

\subsection{A semigroupoid structure on the comma categories}

Let $\mathrm{FGA}^{\circ }$ the disjoint union of the family $\left\{ \left( 
\mathcal{A}\downarrow \mathrm{FGA}\circ \mathcal{B}\right) \right\} _{%
\mathcal{A},\mathcal{B}\in \mathrm{FGA}}$, 
\begin{equation*}
\mathrm{FGA}^{\circ }\doteq \bigvee\nolimits_{\mathcal{A},\mathcal{B}\in 
\mathrm{FGA}}\left( \mathcal{A}\downarrow \mathrm{FGA}\circ \mathcal{B}%
\right) .
\end{equation*}
The objects of $\mathrm{FGA}^{\circ }$ are the disjoint union of objects of $%
\left( \mathcal{A}\downarrow \mathrm{FGA}\circ \mathcal{B}\right) $, for all
pairs $\mathcal{A}\mathbf{,}\mathcal{B}\in \mathrm{FGA}$, and its morphisms
are the union of morphisms of each comma category (in particular, the set of
morphisms from objects of $\left( \mathcal{A}\downarrow \mathrm{FGA}\circ 
\mathcal{B}\right) $ to the ones of $\left( \mathcal{C}\downarrow \mathrm{FGA%
}\circ \mathcal{D}\right) $ are empty, unless $\mathcal{A}=\mathcal{C}$ and $%
\mathcal{B}=\mathcal{D}$). The functors $\frak{P}$ of \S \textbf{1.1 }extend
to an obvious embedding $\mathrm{FGA}^{\circ }\hookrightarrow \mathrm{FGA}$
which we shall also call $\frak{P}$.

From the formal properties of the monoidal product in $\left( \mathrm{FGA}%
,\circ ,\mathcal{I}\right) $, $\mathrm{FGA}^{\circ }$ inherits a
semigroupoid structure given by a partial product functor 
\begin{equation}
\begin{array}{r}
\widehat{\circ }:\,\bigvee\nolimits_{\mathcal{A},\mathcal{B}\in \mathrm{FGA}%
}\left( \mathcal{A}\downarrow \mathrm{FGA}\circ \mathcal{B}\right) \times
\left( \mathcal{B}\downarrow \mathrm{FGA}\circ \mathcal{C}\right)
\longrightarrow \mathrm{FGA}^{\circ }, \\ 
\\ 
\left( \mathcal{A}\downarrow \mathrm{FGA}\circ \mathcal{B}\right) \times
\left( \mathcal{B}\downarrow \mathrm{FGA}\circ \mathcal{C}\right)
\longrightarrow \left( \mathcal{A}\downarrow \mathrm{FGA}\circ \mathcal{C}%
\right) ,
\end{array}
\label{dok}
\end{equation}
such that 
\begin{equation}
\left\langle \varphi ,\mathcal{H}\right\rangle _{_{\mathcal{A}\mathbf{,}%
\mathcal{B}}}\widehat{\circ }\,\left\langle \chi ,\mathcal{G}\right\rangle
_{_{\mathcal{B}\mathbf{,}\mathcal{C}}}=\left\langle \left( id\circ \chi
\right) \,\varphi ,\mathcal{H}\circ \mathcal{G}\right\rangle _{_{\mathcal{A}%
\mathbf{,}\mathcal{C}}},\qquad f\,\widehat{\circ }\,g=f\circ g.  \notag
\end{equation}
If $f$ and $g$ are arrows $\left\langle \varphi ,\mathcal{H}\right\rangle
_{_{\mathcal{A}\mathbf{,}\mathcal{B}}}\rightarrow \left\langle \varphi
^{\prime },\mathcal{H}^{\prime }\right\rangle _{_{\mathcal{A}\mathbf{,}%
\mathcal{B}}}$ and $\left\langle \chi ,\mathcal{G}\right\rangle _{_{\mathcal{%
B}\mathbf{,}\mathcal{C}}}\rightarrow \left\langle \chi ^{\prime },\mathcal{G}%
^{\prime }\right\rangle _{_{\mathcal{B}\mathbf{,}\mathcal{C}}}$, then $%
f\circ g$ will be an arrow $\mathcal{H}\circ \mathcal{G}\rightarrow \mathcal{%
H}^{\prime }\circ \mathcal{G}^{\prime }$ in $\mathrm{FGA}$ satisfying 
\begin{equation*}
\left( id\circ \chi ^{\prime }\right) \,\varphi ^{\prime }=\left( f\circ
g\circ id\right) \,\left( id\circ \chi \right) \,\varphi ,
\end{equation*}
because equations $\varphi ^{\prime }=\left( f\circ id\right) \,\varphi $
and $\chi ^{\prime }=\left( g\circ id\right) \,\chi $ hold. Hence $f\circ g$
is a morphism 
\begin{equation*}
\left\langle \left( id\circ \chi \right) \,\varphi ,\mathcal{H}\circ 
\mathcal{G}\right\rangle _{_{\mathcal{A}\mathbf{,}\mathcal{C}}}\rightarrow
\left\langle \left( id\circ \chi ^{\prime }\right) \,\varphi ^{\prime },%
\mathcal{H}^{\prime }\circ \mathcal{G}^{\prime }\right\rangle _{_{\mathcal{A}%
\mathbf{,}\mathcal{C}}}.
\end{equation*}
The associativity of $\widehat{\circ }$ follows from that of $\circ $, and
it is immediate that the objects $\left\langle \ell _{_{\mathcal{A}}},%
\mathcal{I}\right\rangle \in \left( \mathcal{A}\downarrow \mathrm{FGA}\circ 
\mathcal{A}\right) $ are units for $\widehat{\circ }$. In particular, $%
\widehat{\circ }$ supplies $\left( \mathcal{A}\downarrow \mathrm{FGA}\circ 
\mathcal{A}\right) $ with a monoidal structure, for any $\mathcal{A}\in 
\mathrm{FGA}$. Moreover, $\frak{P}\,\,\widehat{\circ }=\circ \,\left( \frak{P%
}\times \frak{P}\right) $ and $\frak{P}\,\left\langle \ell _{\mathcal{A}},%
\mathcal{I}\right\rangle =\mathcal{I}$, so $\frak{P}:\mathrm{FGA}^{\circ
}\hookrightarrow \mathrm{FGA}$ is a functor of categories with unital
associative partial products.

This groupoid structure on $\mathrm{FGA}^{\circ }$ is independent of the
existence of initial objects on each $\left( \mathcal{A}\downarrow \mathrm{%
FGA}\circ \mathcal{B}\right) $, nevertheless, every object $\left( \mathcal{C%
}\downarrow \mathrm{FGA}\circ \mathcal{I}\right) $, $\mathcal{C}\in \mathrm{%
FGA}$ has an initial object given by the pair $\left\langle r_{_{\mathcal{C}%
}},\mathcal{C}\right\rangle $, being $r_{_{\mathcal{C}}}$ the functorial
isomorphism $\mathcal{C}\backsimeq \mathcal{C}\circ \mathcal{I}$. (That is
true for all monoidal categories.)

Since 
\begin{equation*}
\left\langle \left( id\otimes \delta _{_{\mathcal{C},\mathcal{B}}}\right)
\delta _{_{\mathcal{B},\mathcal{C}}},\underline{hom}\left[ \mathcal{C},%
\mathcal{A}\right] \circ \underline{hom}\left[ \mathcal{B},\mathcal{C}\right]
\right\rangle \in \left( \mathcal{A}\downarrow \mathrm{FGA}\circ \mathcal{B}%
\right)
\end{equation*}
and $\left\langle \ell _{\mathcal{A}},\mathcal{I}\right\rangle \in \left( 
\mathcal{A}\downarrow \mathrm{FGA}\circ \mathcal{A}\right) $, there exist
unique arrows, cocomposition and coidentity (in $\mathrm{FGA}^{\circ }$) 
\begin{equation}
\begin{array}{l}
\Delta _{_{\mathcal{B},\mathcal{C},\mathcal{A}}}:\underline{hom}\left[ 
\mathcal{B},\mathcal{A}\right] \rightarrow \underline{hom}\left[ \mathcal{C},%
\mathcal{A}\right] \circ \underline{hom}\left[ \mathcal{B},\mathcal{C}\right]
, \\ 
\\ 
\varepsilon _{_{\mathcal{A}}}:\underline{hom}\left[ \mathcal{A},\mathcal{A}%
\right] =\underline{end}\left[ \mathcal{A}\right] \rightarrow \mathcal{I},
\end{array}
\label{coi}
\end{equation}
for every $\mathcal{A},\mathcal{B},\mathcal{C}\in \mathrm{FGA}$. There is
also the notion of (left) codual object through $\mathcal{A}^{\ast }\doteq 
\underline{hom}\left[ \mathcal{A},\mathcal{I}\right] $, with coevaluation $%
\mathcal{I}\rightarrow \mathcal{A}^{\ast }\circ \mathcal{A}$. As remarked
above, $\underline{hom}\left[ \mathcal{I},\mathcal{C}\right] \doteq \mathcal{%
C}$, $\forall \mathcal{C}\in \mathrm{FGA}$; hence any coevaluation $\mathcal{%
A}\rightarrow \underline{hom}\left[ \mathcal{B},\mathcal{A}\right] \circ 
\mathcal{B}$ can be regarded as the cocomposition 
\begin{equation}
\underline{hom}\left[ \mathcal{I},\mathcal{A}\right] \rightarrow \underline{%
hom}\left[ \mathcal{B},\mathcal{A}\right] \circ \underline{hom}\left[ 
\mathcal{I},\mathcal{B}\right] ,  \label{ev}
\end{equation}
i.e. $\delta _{_{\mathcal{A}\mathbf{,}\mathcal{B}}}=\Delta _{_{\mathcal{I},%
\mathcal{B},\mathcal{A}}}$. In particular, for duality we have 
\begin{equation}
\underline{hom}\left[ \mathcal{I},\mathcal{I}\right] \rightarrow \underline{%
hom}\left[ \mathcal{A},\mathcal{I}\right] \circ \underline{hom}\left[ 
\mathcal{I},\mathcal{A}\right] .  \label{dual}
\end{equation}
Due to associativity and existence of unit element for $\widehat{\circ }$,
the following diagrams (associative composition\emph{\ }and identity
morphisms) are commutative\footnote{%
Of course, we have a similar commutative diagram where $\mathcal{I}$ is on
the right.} 
\begin{equation}
\begin{diagram} \underline{hom}\left[ \mathcal{C},\mathcal{A}\right] \circ
\underline{hom}\left[ \mathcal{B},\mathcal{C}\right] &\rTo&
\underline{hom}\left[ \mathcal{D},\mathcal{A}\right] \circ
\underline{hom}\left[ \mathcal{C},\mathcal{D}\right] \circ
\underline{hom}\left[ \mathcal{B},\mathcal{C}\right]\\ \uTo& &\uTo \\
\underline{hom}\left[ \mathcal{B},\mathcal{A}\right] &\rTo&
\underline{hom}\left[ \mathcal{D},\mathcal{A}\right] \circ
\underline{hom}\left[ \mathcal{B},\mathcal{D}\right]\\ \end{diagram}
\label{evco}
\end{equation}
\medskip 
\begin{equation}
\begin{diagram} \mathcal{I}\circ \underline{hom}\left[
\mathcal{B},\mathcal{A}\right] & &\lTo & & \underline{end}\left[
\mathcal{A}\right] \circ \underline{hom}\left[
\mathcal{B},\mathcal{A}\right] \\ & \luCorr~{=} & &\ruTo & \\ & &
\underline{hom}\left[ \mathcal{B},\mathcal{A}\right] & &\\ \end{diagram}
\label{idco}
\end{equation}
These diagrams together with Eq. $\left( \ref{ev}\right) $ reflects the
compatibility between the coevaluation and the notions of cocomposition and
coidentity. In dual terms (and for concrete categories) they give equations
like 
\begin{equation}
\left[ f\,g\right] \left( v\right) =f\left( g\left( v\right) \right)
\;\;\;and\;\;\;id\left( v\right) =v.  \label{interp}
\end{equation}
On the other hand, diagram $\left( \ref{idco}\right) $ states that, for
every $\mathcal{A},\mathcal{B}$ in \textrm{FGA}, arrows $\underline{hom}%
\left[ \mathcal{B},\mathcal{A}\right] \rightarrow \underline{end}\left[ 
\mathcal{A}\right] \circ \underline{hom}\left[ \mathcal{B},\mathcal{A}\right]
$ are left invertible.

In particular, there exist injections $\underline{end}\left[ \mathcal{A}%
\right] \hookrightarrow \underline{end}\left[ \mathcal{A}\right] \circ 
\underline{end}\left[ \mathcal{A}\right] $, $\mathcal{A}\hookrightarrow 
\underline{end}\left[ \mathcal{A}\right] \circ \mathcal{A}$ as it was
mentioned in \S \textbf{1.1}.

\subsection{The (co)based categories}

Given a monoidal category $\left( \mathrm{C},\circ ,\mathbf{I}\right) $, a $%
\mathrm{C}$\textbf{-}\emph{based category}\textbf{\ }is (c.f. \cite{mac}):

a) A set of objects $\mathbf{a},\mathbf{b},\mathbf{c,...}$;

b) A function which assigns to each ordered pair of objects $\left( \mathbf{a%
},\mathbf{b}\right) $ an object $Hom\left[ \mathbf{a},\mathbf{b}\right] \in 
\mathrm{C}$, called the $Hom$ objects (to resemble the $Hom$\ $sets$ that
define a proper category);

c) For each ordered triple of objects $\left( \mathbf{a},\mathbf{b},\mathbf{c%
}\right) $ a morphism 
\begin{equation*}
Hom\left[ \mathbf{c},\mathbf{b}\right] \circ Hom\left[ \mathbf{a},\mathbf{c}%
\right] \rightarrow Hom\left[ \mathbf{a},\mathbf{b}\right] ;
\end{equation*}

d) For each object a morphism $\mathbf{I}\rightarrow Hom\left[ \mathbf{a},%
\mathbf{a}\right] \doteq End\left[ \mathbf{a}\right] $.

These morphisms must satisfy obvious associativity and unit constraints
(dual to diagrams $\left( \ref{evco}\right) $ and $\left( \ref{idco}\right) $%
).\footnote{%
Any (small) category can be regarded as a $\mathrm{Set}$-based category.}

In a (co)based category there is not, \emph{a priori}, a notion of
(co)evaluation.

It is clear that every category $\left( \mathrm{C},\circ ,\mathbf{I}\right) $
with internal Hom objects gives rise to a $\mathrm{C}$-based category with
objects in $\mathrm{C}$.\footnote{%
In the latter case and for concrete categories (in which the objects $Hom%
\left[ \mathbf{a},\mathbf{b}\right] $ define Hom-sets), it says that $%
\mathrm{C}$ has been \emph{enriched} by the objects $Hom\left[ \mathbf{a},%
\mathbf{b}\right] \in \mathrm{C}$. This is the case of $\left( \mathrm{QA}%
,\bullet \right) $, which is enriched by its internal Hom objects $\mathcal{B%
}^{!}\circ \mathcal{A}$.} In the dual case, of internal\ coHom objects,
there is a related $\mathrm{C}^{op}$-based category, which we call $\mathrm{C%
}$\textbf{-}\emph{cobased} category. The function $\left( \mathcal{A},%
\mathcal{B}\right) \mapsto \underline{hom}\left[ \mathcal{B},\mathcal{A}%
\right] $, together with arrows $\left( \ref{coi}\right) $ and diagrams $%
\left( \ref{evco}\right) $ and $\left( \ref{idco}\right) $, form an $\mathrm{%
FGA}$-cobased category. Dually, the function 
\begin{equation*}
\left( \mathcal{B}^{op},\mathcal{A}^{op}\right) \mapsto \underline{Hom}\left[
\mathcal{B}^{op},\mathcal{A}^{op}\right] =\underline{hom}\left[ \mathcal{B},%
\mathcal{A}\right] ^{op},
\end{equation*}
with reversed arrows and diagrams, defines a $\mathsf{QLS}$-based category%
\textbf{\ }with objects in $\mathsf{QLS}$.

This can be settled in more general terms: suppose a family of categories $%
\left\{ \digamma ^{\mathcal{A},\mathcal{B}}/\mathcal{A},\mathcal{B}\in 
\mathrm{FGA}\right\} $ (or equivalently, a function $\left( \mathcal{A},%
\mathcal{B}\right) \mapsto \digamma ^{\mathcal{A},\mathcal{B}}$ to
categories) with initial objects is given together a related family of
functors $\frak{P}^{\digamma }:\digamma ^{\mathcal{A},\mathcal{B}%
}\rightarrow \mathrm{FGA}$ and a semigroupoid structure $\circ _{\digamma }$
on its disjoint union $\digamma ^{\cdot }$. Suppose in addition that 
\begin{equation*}
\frak{P}^{\digamma }\,\circ _{\digamma }=\circ \,\left( \frak{P}^{\digamma
}\times \frak{P}^{\digamma }\right)
\end{equation*}
and the respective units are preserved. Then, the function $\left( \mathcal{A%
},\mathcal{B}\right) \mapsto \underline{hom}^{\digamma }\left[ \mathcal{B},%
\mathcal{A}\right] $, being $\underline{hom}^{\digamma }\left[ \mathcal{B},%
\mathcal{A}\right] $ the image under $\frak{P}^{\digamma }$ of the initial
object of $\digamma ^{\mathcal{A},\mathcal{B}}$, defines an $\mathrm{FGA}$%
-cobased category (for the dual case we just must change initial by terminal
objects).\smallskip\ That is a direct consequence of the discussion we have
made in the previous subsection.

\bigskip

\textbf{\textbf{Example.} }Consider the category $\mathrm{CA}_{\star }$ of
conic quantum spaces equipped with a marked basis (and with the same
morphisms of $\mathrm{CA}$). For each $\mathcal{A}_{\star }=\left( \mathcal{A%
},g\right) ,\mathcal{B}_{\star }=\left( \mathcal{B},f\right) \in \mathrm{CA}%
_{\star }$, where $g=\left\{ a_{i}\right\} \subset \mathbf{A}_{1}$ and $%
f=\left\{ b_{i}\right\} \subset \mathbf{B}_{1}$ are the marked basis of $%
\mathcal{A}$ and $\mathcal{B}$, we define the categories $\digamma $'s as
the ones given by triples $\left( \mathcal{H}_{\star },\varphi ,\varphi
^{\dagger }\right) $, being $\varphi :\mathcal{A}\rightarrow \mathcal{H}%
\circ \mathcal{B}$ and $\varphi ^{\dagger }:\mathcal{A}^{!}\rightarrow 
\mathcal{H}\circ \mathcal{B}^{!}$ arrows in $\mathrm{CA}$ such that 
\begin{equation}
\left\langle \varphi ^{\dagger }\left( a^{i}\right) ,b_{j}\right\rangle
=\left\langle b^{j},\varphi \left( a_{i}\right) \right\rangle ,
\label{duren}
\end{equation}
where $\left\langle h\otimes x,y\right\rangle =\left\langle x,h\otimes
y\right\rangle =\left\langle x,y\right\rangle \,h\in \mathbf{H}$, and $%
\left\{ a^{i}\right\} ,\left\{ b^{i}\right\} $ the corresponding duals to $g$
and $f$. Arrows 
\begin{equation*}
\left( \mathcal{H}_{\star },\varphi ,\varphi ^{\dagger }\right) \rightarrow
\left( \mathcal{G}_{\star },\psi ,\psi ^{\dagger }\right)
\end{equation*}
between triples are morphisms $\mathcal{H}\overset{\alpha }{\rightarrow }%
\mathcal{G}$ in $\mathrm{CA}$ such that $\psi =\left( \alpha \circ I\right)
\,\varphi $ (i.e. they are morphisms in $\left( \mathcal{A}\downarrow 
\mathrm{CA}\circ \mathcal{B}\right) $).\footnote{%
Note that $\psi =\left( \alpha \circ I\right) \,\varphi $ and Eq. $\left( 
\ref{duren}\right) $ implies $\psi ^{\dagger }=\left( \alpha \circ I\right)
\,\varphi ^{\dagger }$.} Of course, the functors $\frak{P}^{\digamma }$ are
the embedding $\left( \mathcal{H}_{\star },\varphi ,\varphi ^{\dagger
}\right) \mapsto \mathcal{H}$. The initial objects are quotients of the
corresponding to $\left( \mathcal{A}\downarrow \mathrm{CA}\circ \mathcal{B}%
\right) $. In particular, for each $\digamma ^{\mathcal{A}_{\star }}$ we can
define $\underline{end}^{\digamma }\left[ \mathcal{A}_{\star }\right] \doteq
e\left[ \mathcal{A},g\right] $, the one defined by Manin (for $\mathrm{QA}$)
in \cite{man0} to add the so-called \emph{missing relations}.\footnote{%
Such relations turns $\underline{end}\left[ \mathcal{A}\right] $ into a
commutative algebra $e\left[ \mathcal{A}\right] $ whenever $\mathcal{A}$ is
commutative.} The semigroupoid structure on $\digamma ^{\cdot }$ is given by 
\begin{equation*}
\left( \mathcal{H}_{\star },\varphi ,\varphi ^{\dagger }\right) \times
\left( \mathcal{G}_{\star },\psi ,\psi ^{\dagger }\right) \mapsto \left( 
\left[ \mathcal{H}\circ \mathcal{G}\right] _{\star },\left( I_{H}\circ \psi
\right) \,\varphi ,\left( I_{H}\circ \psi ^{\dagger }\right) \,\varphi
^{\dagger }\right)
\end{equation*}
and $\alpha \times \beta \mapsto \alpha \circ \beta $, where the marked
basis of $\mathcal{H}\circ \mathcal{G}$ is the direct product between the
ones of $\mathcal{H}$ and $\mathcal{G}$.\ \ \ $\blacksquare \smallskip
\smallskip $

\bigskip

In \S \textbf{4} we construct a family of categories $\left\{ \Upsilon ^{%
\mathcal{A},\mathcal{B}}:\mathcal{A},\mathcal{B}\in \mathrm{FGA}\right\} $
whose objects are essentially pairs 
\begin{equation*}
\left\langle \mathcal{A}\rightarrow \mathcal{H}\circ _{\tau }\mathcal{B},%
\mathcal{H}\right\rangle ,\;\;\mathcal{H}\in \mathrm{FGA},
\end{equation*}
with $\circ _{\tau }$ some given twisted tensor product between $\mathcal{H}$
and $\mathcal{B}$ (as we shall define in the next chapter). Each category $%
\Upsilon ^{\mathcal{A},\mathcal{B}}$ has an initial object, $\underline{hom}%
^{\Upsilon }\left[ \mathcal{B},\mathcal{A}\right] $, and a map $\mathcal{A}%
\rightarrow \underline{hom}^{\Upsilon }\left[ \mathcal{B},\mathcal{A}\right]
\circ _{\tau }\mathcal{B}$. They also have a related faithful functor $\frak{%
P}^{\Upsilon }:\Upsilon ^{\mathcal{A},\mathcal{B}}\hookrightarrow \mathrm{FGA%
}$. Moreover, on the subfamily $\left\{ \Upsilon ^{\mathcal{A},\mathcal{B}}:%
\mathcal{A},\mathcal{B}\in \mathrm{CA}\right\} $ corresponding to the conic
case, a semigroupoid structure $\circ _{\Upsilon }$ can be defined in such a
way that 
\begin{equation*}
\frak{P}^{\Upsilon }\,\circ _{\Upsilon }=\circ \,\left( \frak{P}^{\Upsilon
}\times \frak{P}^{\Upsilon }\right) .
\end{equation*}
As a consequence, the function $\left( \mathcal{A},\mathcal{B}\right)
\rightarrow \underline{hom}^{\Upsilon }\left[ \mathcal{B},\mathcal{A}\right] 
$ defines an $\mathrm{CA}$-cobased category and, in addition, the maps $%
\mathcal{A}\rightarrow \underline{hom}^{\Upsilon }\left[ \mathcal{B},%
\mathcal{A}\right] \circ _{\tau }\mathcal{B}$ provide a coevaluation notion.
(We must mention that the last notion is not compatible with cocomposition
and counit -see the comments before Eq. $\left( \ref{interp}\right) $-,
because in general, $\underline{hom}^{\Upsilon }\left[ \mathcal{K},\mathcal{A%
}\right] \neq \mathcal{A}$).

\section{Twisted tensor product of quantum spaces}

In this chapter we define a twisted tensor product structure on quantum
spaces in a way intimately related to the corresponding algebra case,
following the developments accomplished in ref. \cite{cap} . In doing so, we
firstly redefine the twisted tensor products on the category $\Bbbk \mathrm{%
-Alg}$ (with respect to its usual monoid $\otimes $) in geometric terms.

\subsection{The twisted tensor products of algebras and the geometric
language}

Given an algebra $\left( \mathbf{A},m,\eta \right) $, $m$ denotes its
associative product and $\eta :\Bbbk \rightarrow \mathbf{A}$ its unit map.
The unit element of any algebra (including $\Bbbk $) is denoted by $``1"$.

Given $\mathbf{A},\mathbf{B}\in \Bbbk \mathrm{-Alg}$, the set $Hom_{\Bbbk 
\mathrm{-Alg}}\left[ \mathbf{A},\mathbf{B}\right] $ of algebra homomorphisms
is also called the set of $\mathbf{B}$\emph{-points} of $\mathbf{A}$. This
terminology is adequate when we are thinking about $\mathbf{A}\ $and $%
\mathbf{B}$ as the coordinate (or function) rings of certain
-noncommutative- spaces. It says that a $\mathbf{B}$-point of $\mathbf{A}$
is \emph{generic} if it is defined by a monomorphism $\mathbf{%
A\hookrightarrow B}$; and two $\mathbf{B}$-points $\alpha :\mathbf{X}%
\rightarrow \mathbf{B}$ of $\mathbf{X}$ and $\beta :\mathbf{Y}$ $\rightarrow 
\mathbf{B}$ of $\mathbf{Y}$ are \emph{commuting} if 
\begin{equation*}
\alpha \left( x\right) \cdot \beta \left( y\right) =\beta \left( y\right)
\cdot \alpha \left( x\right) ,\;\forall x\in \mathbf{X},\;\forall y\in 
\mathbf{Y},
\end{equation*}
or 
\begin{equation*}
\left[ \alpha ,\beta \right] _{B}\doteq \left( m_{B}-m_{B}^{op}\right)
\,\left( \alpha \otimes \beta \right) =0.
\end{equation*}

Let us express in these terms the definition of twisted tensor products of
algebras given in \cite{cap}.

\begin{definition}
Consider a pair of algebras $\mathbf{A}\ $and $\mathbf{B}.$ The \textbf{%
twisted tensor product \emph{(}TTP}\emph{)}\textbf{\ of} $\mathbf{A}\ $%
\textbf{and} $\mathbf{B}$ is a \textbf{triple} $\left( \mathbf{C}%
,i_{A},i_{B}\right) $ where $\mathbf{C}$ is an algebra, $i_{A}$ is a generic 
$\mathbf{C}$-point of $\mathbf{A}$ and $i_{B}$ is a generic $\mathbf{C}$%
-point of $\mathbf{B}$, such that the linear map 
\begin{equation*}
\varphi \doteq m_{C}\,\left( i_{A}\otimes i_{B}\right) :\mathbf{A}\otimes 
\mathbf{B\rightarrow C}
\end{equation*}
is an isomorphism. A morphism of TTP's $\left( \mathbf{C},i_{A},i_{B}\right)
\rightarrow \left( \mathbf{C}^{\prime },i_{A}^{\prime },i_{B}^{\prime
}\right) $ is given by an homomorphism of algebras $\rho :\mathbf{C}%
\rightarrow \mathbf{C}^{\prime }$, with $\rho \,i_{A}=i_{A}^{\prime }$ and $%
\rho \,i_{B}=i_{B}^{\prime }$.\ \ \ $\blacksquare $
\end{definition}

Given algebras $\left( \mathbf{A},m_{A},\eta _{A}\right) $ and $\left( 
\mathbf{B},m_{B},\eta _{B}\right) $, examples of these objects are algebras $%
\mathbf{A}\otimes _{\tau }\mathbf{B}$ built over the vector space $\mathbf{A}%
\otimes \mathbf{B}$ with product and unit map 
\begin{equation*}
m_{A\otimes B}^{\tau }=\left( m_{A}\otimes m_{B}\right) \,\left(
I_{A}\otimes \tau \otimes I_{B}\right) ;\;\;\eta _{A\otimes B}=\eta
_{A}\otimes \eta _{B},
\end{equation*}
being $\tau :\mathbf{B}\otimes \mathbf{A\rightarrow A}\otimes \mathbf{B}$ a
linear map satisfying 
\begin{equation}
\begin{array}{l}
\tau \left( m_{B}\otimes I_{A}\right) =\left( I_{A}\otimes m_{B}\right)
\,\left( \tau \otimes I_{B}\right) \,\left( I_{B}\otimes \tau \right) , \\ 
\tau \left( I_{B}\otimes m_{A}\right) =\left( m_{A}\otimes I_{B}\right)
\,\left( I_{A}\otimes \tau \right) \,\left( \tau \otimes I_{A}\right) ,
\end{array}
\label{mu}
\end{equation}
and 
\begin{equation}
\tau \left( I_{B}\otimes \eta _{A}\right) =\eta _{A}\otimes I_{B};\;\;\tau
\left( \eta _{B}\otimes I_{A}\right) =I_{A}\otimes \eta _{B}.  \label{un}
\end{equation}
The generic $\mathbf{A}\otimes _{\tau }\mathbf{B}$\textbf{-}points
correspond to the usual inclusion of algebras 
\begin{equation*}
i_{A}\doteq \mathbf{i}_{A}:a\mapsto a\otimes 1,\;\;i_{B}\doteq \mathbf{i}%
_{B}:b\mapsto 1\otimes b.
\end{equation*}
The maps $\tau $ (satisfying $\left( \ref{mu}\right) $ and $\left( \ref{un}%
\right) $) are called \emph{twisting maps}. Essentially, these are all of
TPP's, namely: any triple $\left( \mathbf{C},i_{A},i_{B}\right) $ is
isomorphic to a unique triple of the form $\left( \mathbf{A}\otimes _{\tau }%
\mathbf{B},\mathbf{i}_{A},\mathbf{i}_{B}\right) $, i.e. there exist a unique
twisting map $\tau :\mathbf{B}\otimes \mathbf{A\rightarrow A}\otimes \mathbf{%
B}$ such that 
\begin{equation*}
\left( \mathbf{A}\otimes _{\tau }\mathbf{B},\mathbf{i}_{A},\mathbf{i}%
_{B}\right) \backsimeq \left( \mathbf{C},i_{A},i_{B}\right)
\end{equation*}
(for a proof see \cite{cap}). The isomorphism is given by $\varphi :$ $%
\mathbf{A}\otimes \mathbf{B}\backsimeq \mathbf{C}$ and $\tau $ by the
equation 
\begin{equation}
\tau \left( b\otimes a\right) =\varphi ^{-1}\left( i_{B}\left( b\right)
\cdot i_{A}\left( a\right) \right) ,\;a\in \mathbf{A},\;b\in \mathbf{B}.
\label{7}
\end{equation}
In other words, the equivalence classes of TTP's are in one to one
correspondence to twisting maps.

We often refer to a TTP of $\mathbf{A}$ and $\mathbf{B}$ as an algebra $%
\mathbf{A}\otimes _{\tau }\mathbf{B}$ (omitting the canonical inclusion
maps), being $\tau $ some twisting map.

Any triple $\left( \mathbf{C},i_{\Bbbk },i_{A}\right) $ or $\left( \mathbf{C}%
,i_{A},i_{\Bbbk }\right) $ (remember that $\Bbbk $ is a unit object for the
monoidal category $\left( \Bbbk \mathrm{-Alg},\otimes \right) $), because of
the Eq. $\left( \ref{un}\right) $, is isomorphic to the algebra $\mathbf{A}$%
, i.e. $\mathbf{A}\otimes _{\tau }\Bbbk =\Bbbk \otimes _{\tau }\mathbf{A}=%
\mathbf{A}$ for all twisting map $\tau $.

We shall see that the canonical \emph{flipping map }$\tau _{o}$, $\tau
_{o}\left( b\otimes a\right) =a\otimes b$, is related to the TTP's with
commuting generic points, as it could be expect and is implicit in \cite{cap}%
. This class of TTP's belongs to the class of the usual tensor product of
algebras $\mathbf{A}\otimes \mathbf{B}$.\smallskip

\begin{proposition}
Given a pair of algebras $\mathbf{A}$ and $\mathbf{B}$, the set of twisted
tensor products $\left( \mathbf{C},i_{A},i_{B}\right) $ such that $\left[
i_{A},i_{B}\right] _{C}=0$, i.e. $i_{A}\ $and $i_{B}$ are commuting points,
form an equivalence class whose associated twisting map is the flipping map.
\end{proposition}

\begin{proof}
Let $\left( \mathbf{C},i_{A},i_{B}\right) $ be a TTP of $\mathbf{A}$ and $%
\mathbf{B}$. By the above theorem, there exist a unique twisting map $\tau
\left( b\otimes a\right) =\varphi ^{-1}\left( i_{B}\left( b\right) \cdot
i_{A}\left( a\right) \right) $, such that $\varphi $ defines the algebra
isomorphism $\mathbf{A}\otimes _{\tau }\mathbf{B\backsimeq C}$, satisfying 
\begin{equation*}
\varphi \,\,\mathbf{i}_{A}=i_{A}\;\;\;and\;\;\;\varphi \,\,\mathbf{i}%
_{B}=i_{B}.
\end{equation*}
Therefore, if $\left[ i_{A},i_{B}\right] _{C}=0$, 
\begin{eqnarray*}
\tau \left( b\otimes a\right) =\varphi ^{-1}\left( i_{B}\left( b\right)
\cdot i_{A}\left( a\right) \right) =\varphi ^{-1}\left( i_{A}\left( a\right)
\cdot i_{B}\left( b\right) \right) \\
=\varphi ^{-1}\left( i_{A}\left( a\right) \right) \cdot _{\tau }\varphi
^{-1}\left( i_{B}\left( b\right) \right) =\mathbf{i}_{A}\left( a\right)
\cdot _{\tau }\mathbf{i}_{B}\left( b\right) \\
=\left( a\otimes 1\right) \cdot _{\tau }\left( 1\otimes b\right) =a\otimes b
\end{eqnarray*}
where we are denoting $``\cdot "$ for the product in $\mathbf{C}$ and $%
``\cdot _{\tau }"$ for the one in $\mathbf{A}\otimes _{\tau }\mathbf{B}$.
The last equation means that $\tau $ is the canonical flipping map.
\end{proof}

Therefore, it can be stated that those twisting $\tau $\ different from
flipping map give rise to a (direct)\ product space where the points of each
factor do not commute among themselves.

For latter convenience, we will characterize the classes of TTP's with
bijective twisting maps.

\begin{definition}
Let $\left( \mathbf{C},i_{A},i_{B}\right) $ be a $\mathbf{TTP}$ for $\mathbf{%
A}$ and $\mathbf{B}$. It will be called \textbf{symmetric} twisted tensor
product $\left( \mathbf{STTP}\right) $ provided the linear map 
\begin{equation*}
\varphi ^{op}\doteq m_{C}\,\left( i_{B}\otimes i_{A}\right) :\mathbf{B}%
\otimes \mathbf{A\rightarrow C}
\end{equation*}
is a bijection.\ \ \ $\blacksquare $
\end{definition}

\begin{proposition}
The class of a \textbf{STTP} corresponds to a \textbf{bijective} twisting
map, and viceversa.
\end{proposition}

\begin{proof}
\emph{\ }The twisting map related to $\left( \mathbf{C},i_{A},i_{B}\right) $
is $\tau =\varphi ^{-1}\,\varphi ^{op}$, as it can see from $\left( \ref{7}%
\right) $. Then $\tau $ is bijective \emph{iff} $\varphi ^{op}$ is a
bijection.
\end{proof}

In these cases $\varphi ^{op}$ defines the algebra isomorphism $\mathbf{B}%
\otimes _{\tau ^{op}}\mathbf{A}\backsimeq \mathbf{C}$, with 
\begin{equation*}
\tau ^{op}=\left( \varphi ^{op}\right) ^{-1}\,\varphi =\left( \varphi
^{-1}\,\varphi ^{op}\right) ^{-1}=\tau ^{-1}.
\end{equation*}
In particular we have $\mathbf{B}\otimes _{\tau ^{-1}}\mathbf{A}\backsimeq 
\mathbf{A}\otimes _{\tau }\mathbf{B}$.

Now we shall rephrase the above definitions for quantum spaces.

\subsection{Twisted monoid on $\mathrm{FGA}$'s objects}

The monoid $\circ $ stems from the usual tensor product $\otimes $, being $%
\mathbf{A}\circ \mathbf{B}$ a subalgebra of $\mathbf{A}\otimes \mathbf{B}$.
This enables us to introduce a natural definition of \textbf{TTP} on quantum
spaces as follows:

\begin{definition}
Let $\mathcal{A}=\left( \mathbf{A}_{1},\mathbf{A}\right) $ and $\mathcal{B}%
=\left( \mathbf{B}_{1},\mathbf{B}\right) $ be two objects of $\mathrm{FGA}$.
A \textbf{TTP of} $\mathcal{A}$ \textbf{and} $\mathcal{B}$ is a quadruple 
\begin{equation*}
\left( \mathcal{C},\mathbf{C},i_{A},i_{B}\right) ,\;\;\;\mathcal{C}=\left( 
\mathbf{C}_{1},\widetilde{\mathbf{C}}\right) ,
\end{equation*}
with $\widetilde{\mathbf{C}}$ a subalgebra of $\mathbf{C}$, and $i_{A}$ and $%
i_{B}$ are generic $\mathbf{C}$-points of $\mathbf{A}$ and $\mathbf{B}$,
resp., such that the linear map $\varphi \doteq m_{C}\,\left( i_{A}\otimes
i_{B}\right) $ is an isomorphism, with the additional condition that its
restriction to $\mathbf{A}_{1}\otimes \mathbf{B}_{1}$ gives $\varphi _{1}:%
\mathbf{A}_{1}\otimes \mathbf{B}_{1}\backsimeq \mathbf{C}_{1}$. A morphism
of TTP's is a morphism of algebras $\rho :\mathbf{C}\rightarrow \mathbf{C}%
^{\prime }$, with $\rho \,i_{A}=i_{A}^{\prime }$, $\rho
\,i_{B}=i_{B}^{\prime }$, and $\rho \left( \mathbf{C}_{1}\right) \subset 
\mathbf{C}_{1}^{\prime }$, i.e. $\rho $ restricted to $\widetilde{\mathbf{C}}
$ defines a morphism $\mathcal{C}\rightarrow \mathcal{C}^{\prime }$%
.\smallskip\ \ \ $\blacksquare $
\end{definition}

\begin{proposition}
Any quadruple $\left( \mathcal{C},\mathbf{C},i_{A},i_{B}\right) $ is
isomorphic to a unique one 
\begin{equation*}
\left( \mathcal{A}\circ _{\tau }\mathcal{B},\mathbf{A}\otimes _{\tau }%
\mathbf{B},\mathbf{i}_{A},\mathbf{i}_{B}\right) ,
\end{equation*}
where $\tau $ is a twisting map, $\mathcal{A}\circ _{\tau }\mathcal{B}\doteq
\left( \mathbf{A}_{1}\otimes \mathbf{B}_{1}\mathbf{,\mathbf{A}}\circ _{\tau }%
\mathbf{\mathbf{B}}\right) $, and $\mathbf{A}\circ _{\tau }\mathbf{B}$ is
the subalgebra of $\mathbf{A}\otimes _{\tau }\mathbf{B}$ generated by $%
\mathbf{A}_{1}\otimes \mathbf{B}_{1}$. The isomorphism is given by $\varphi $%
, the associated isomorphism of quantum spaces is 
\begin{equation*}
\left. \varphi \right| _{\mathbf{A}\circ _{\tau }\mathbf{B}}:\mathcal{A}%
\circ _{\tau }\mathcal{B}\backsimeq \mathcal{C},
\end{equation*}
and again $\tau \left( b\otimes a\right) =\varphi ^{-1}\left( i_{B}\left(
b\right) \cdot i_{A}\left( a\right) \right) $, $a\in \mathbf{A}$, $b\in 
\mathbf{B}$.
\end{proposition}

\begin{proof}
It follows from the (general) algebra case.
\end{proof}

Since this is a one to one correspondence, we often refer to TTP's of $%
\mathcal{A}$ and $\mathcal{B}$ as a quantum space $\mathcal{A}\circ _{\tau }%
\mathcal{B}$, being $\tau $ some twisting map.

\begin{remark}
Analogously to the previous section, any quadruple $\left( \mathcal{C},%
\mathbf{C},i_{A},i_{\Bbbk }\right) $ is isomorphic to $\mathcal{A}$, because 
$\mathcal{A}\circ _{\tau }\mathcal{I}=\mathcal{I}\circ _{\tau }\mathcal{A}=$ 
$\mathcal{A}$ for any twisting map \emph{(}remember that $\mathcal{I}$ is a
unit object for $\left( \mathrm{FGA},\circ \right) $\emph{)}. The same is
not true for the unit element $\mathcal{K}$ of $\left( \mathrm{CA},\circ
\right) $ \emph{(}and every $\left( \mathrm{CA}^{m},\circ \right) $\emph{)},
because the underlying algebra of $\mathcal{K}=\left( \Bbbk ,\Bbbk \left[ e%
\right] \right) $ is not just generated by the unit element \emph{(}as the
object $\mathcal{I}=\left( \Bbbk ,\Bbbk \right) $\emph{)}, but for the set $%
\left\{ 1,e^{n}\right\} _{n\geq 1}$.\ \ \ $\blacksquare $
\end{remark}

Respect to the relationship between commuting points and the monoid $\circ $
of quantum spaces, we have:

\begin{proposition}
Given $\mathcal{A}$,$\mathcal{B}\in \mathrm{FGA}$, the set of TTP's $\left( 
\mathcal{C},\mathbf{C},i_{A},i_{B}\right) $ such that $\left[ i_{A},i_{B}%
\right] _{C}=0$, form an equivalence class with associated twisting map
equal to the flipping map. In other words, the above set is the equivalence
class of $\mathcal{A}\circ \mathcal{B}$.\ \ \ $\blacksquare $
\end{proposition}

It is worth mentioning that the underlying vector spaces of $\mathbf{A}\circ
_{\tau }\mathbf{B}$ and $\mathbf{A}\circ \mathbf{B}$ (subspaces of $\mathbf{A%
}\otimes \mathbf{B}$) do not coincide in general\footnote{%
It is enough to take, for $\mathbf{A}=\mathbf{A}_{1}^{\otimes }$ and $%
\mathbf{B}=\mathbf{B}_{1}^{\otimes }$, a twisting map $\tau :\mathbf{B}%
\otimes \mathbf{A\rightarrow A}\otimes \mathbf{B}$ satisfying $\left( \ref
{un}\right) $, and in some complement $\mathbf{l}^{c}$ of $\mathbf{l}=%
\mathbf{B}\otimes \Bbbk +\Bbbk \otimes \mathbf{A}$ ($\Bbbk =\mathbf{A}_{0},%
\mathbf{B}_{0}$), to put $\left. \tau \right| _{\mathbf{l}^{c}}=0$. Hence,
we will have $\mathbf{A}\circ _{\tau }\mathbf{B}=\Bbbk $, clearly distinct
from $\mathbf{A}\circ \mathbf{B}$.} (while for the algebra case we have the
equality of vector spaces $\mathbf{A}\otimes \mathbf{B}=\mathbf{A}\otimes
_{\tau }\mathbf{B}$). Nevertheless, let us consider the symmetric analogous
of TTP of algebras:

\begin{definition}
A \textbf{symmetric} twisted tensor product of $\mathcal{A}$ and $\mathcal{B}
$, is a twisted tensor product $\left( \mathcal{C},\mathbf{C}%
,i_{A},i_{B}\right) $ such that 
\begin{equation*}
\varphi ^{op}\doteq m_{C}\,\left( i_{B}\otimes i_{A}\right) :\mathbf{B}%
\otimes \mathbf{A\rightarrow C}
\end{equation*}
is a linear isomorphism which restricted to $\mathbf{B}_{1}\otimes \mathbf{A}%
_{1}$ gives $\mathbf{B}_{1}\otimes \mathbf{A}_{1}\backsimeq \mathbf{C}_{1}$%
.\ \ \ $\blacksquare $
\end{definition}

As in the algebra case we have $\mathcal{A}\circ _{\tau }\mathcal{B}%
\backsimeq \mathcal{B}\circ _{\tau ^{-1}}\mathcal{A}$.

A direct consequence of the above definition is the following.

\begin{proposition}
The class of an \textbf{STTP} corresponds to a \textbf{bijective} twisting
map such that 
\begin{equation}
\tau \left( \mathbf{B}_{1}\otimes \mathbf{A}_{1}\right) \subset \mathbf{A}%
_{1}\otimes \mathbf{B}_{1},  \label{9}
\end{equation}
and viceversa.
\end{proposition}

\begin{proof}
As in the previous section, the twisting map related to $\left( \mathcal{C},%
\mathbf{C},i_{A},i_{B}\right) $ can be written $\tau =\varphi ^{-1}\,\varphi
^{op}$, from which we see that $\tau $ is a bijection \emph{iff} $\varphi
^{op}$ is a bijection. On the other hand, 
\begin{equation*}
\varphi \,\tau \left( \mathbf{B}_{1}\otimes \mathbf{A}_{1}\right) =\varphi
^{op}\left( \mathbf{B}_{1}\otimes \mathbf{A}_{1}\right) \subset \mathbf{C}%
_{1}
\end{equation*}

\emph{iff} Eq. $\left( \ref{9}\right) $ holds.
\end{proof}

We call \emph{symmetric twisting maps} those $\tau $'s related to STTP's.
They are completely defined by isomorphisms $\left. \mathbb{\tau }\right| _{%
\mathbf{B}_{1}\otimes \mathbf{A}_{1}}:\mathbf{B}_{1}\otimes \mathbf{A}%
_{1}\backsimeq \mathbf{A}_{1}\otimes \mathbf{B}_{1}$, since properties $%
\left( \ref{mu}\right) $ and $\left( \ref{un}\right) $. For convenience, we
often write the above maps as 
\begin{equation*}
\widehat{\mathbb{\tau }}:\mathbf{B}_{1}\otimes \mathbf{A}_{1}\backsimeq 
\mathbf{B}_{1}\otimes \mathbf{A}_{1},
\end{equation*}
using the flipping map $\mathbb{\tau }_{o}$ to reorder the factors, i.e. $%
\mathbb{\,\widehat{\mathbb{\tau }}=\tau }_{o}\,\mathbb{\tau }.$ Therefore,
given a pair of basis $\left\{ a_{i}\right\} $ and $\left\{ b_{i}\right\} $
of $\mathbf{A}_{1}$ and $\mathbf{B}_{1}$, resp., we have $\tau $ defined by 
\begin{equation*}
\widehat{\mathbb{\tau }}\left( b_{i}\otimes a_{j}\right) =\widehat{\mathbb{%
\tau }}_{ij}^{lk}\;b_{l}\otimes a_{k},
\end{equation*}
where we are using the sum over repeated indices convention.

\medskip

\textbf{\textbf{Example.} }Consider a pair of quadratic quantum spaces $%
\mathcal{A}$ and $\mathcal{B}$, such that, in terms of basis $\left\{
a_{i}\right\} _{i=1}^{\frak{n}}$ and $\left\{ b_{i}\right\} _{i=1}^{\frak{m}%
} $ of $\mathbf{A}_{1}$ and $\mathbf{B}_{1}$, resp., the kernel of the
associated canonical maps $\mathbf{A}_{1}^{\otimes }\twoheadrightarrow 
\mathbf{A}$ and $\mathbf{B}_{1}^{\otimes }\twoheadrightarrow \mathbf{B}$ are
generated by subspaces 
\begin{equation*}
span\left[ A_{ij}^{kl}\,a_{k}\otimes a_{l}\right] _{i,j=1}^{\frak{n}}\subset 
\mathbf{A}_{1}^{\otimes 2}\;\;and\;\;span\left[ B_{ij}^{kl}\,b_{k}\otimes
b_{l}\right] _{i,j=1}^{\frak{m}}\subset \mathbf{B}_{1}^{\otimes 2}.
\end{equation*}
Taking an invertible matrix $\widehat{\mathbb{\tau }}$ such that there exist 
$\Lambda $ and $\Omega $ satisfying 
\begin{equation*}
A_{ij}^{kl}\,\widehat{\mathbb{\tau }}_{rk}^{sp}\,\widehat{\mathbb{\tau }}%
_{sl}^{tq}=\Lambda _{rij}^{tkl}\,A_{kl}^{pq}\;\;and\;\;B_{ij}^{kl}\,\widehat{%
\mathbb{\tau }}_{lr}^{qs}\,\widehat{\mathbb{\tau }}_{ks}^{pt}=\Omega
_{rij}^{tkl}\,B_{kl}^{pq},
\end{equation*}
we get a STTP $\mathcal{A}\circ _{\tau }\mathcal{B}$. (Such matrices $%
\Lambda $ and $\Omega $ insure $\widehat{\mathbb{\tau }}$ preserve the
ideals.) For example, calling $\mathbb{A}$ and $\mathbb{B}$ the
endomorphisms $\mathbf{A}_{1}^{\otimes 2}\rightarrow \mathbf{A}_{1}^{\otimes
2}$ and $\mathbf{B}_{1}^{\otimes 2}\rightarrow \mathbf{B}_{1}^{\otimes 2}$,
defined by the coefficients $A_{ij}^{kl}$ and $B_{ij}^{kl}$ (on above
basis), resp., any pair of linear isomorphisms $\alpha :\mathbf{A}%
_{1}\backsimeq \mathbf{A}_{1}$ and $\beta :\mathbf{B}_{1}\backsimeq \mathbf{B%
}_{1}$ such that $\left[ \mathbb{A},\alpha \otimes \alpha \right] =\left[ 
\mathbb{B},\beta \otimes \beta \right] =0$, give rise to a symmetric
twisting map $\mathbb{\tau }$ with $\widehat{\mathbb{\tau }}=\beta \otimes
\alpha $. In fact, for such a case we can take 
\begin{equation*}
\Lambda _{rij}^{tkl}=\widehat{\mathbb{\tau }}_{ri}^{sk}\,\widehat{\mathbb{%
\tau }}_{sj}^{tl}\;\;and\;\;\Omega _{rij}^{tkl}=\widehat{\mathbb{\tau }}%
_{jr}^{ls}\,\widehat{\mathbb{\tau }}_{is}^{kt}.\;\;\;\;\blacksquare
\end{equation*}

\medskip

Now, the main properties of STTP's.

\begin{proposition}
Given a symmetric twisting map $\mathbb{\tau }$ on $\mathcal{A}$ and $%
\mathcal{B}$, the underlying vector spaces to $\mathbf{A}\circ _{\tau }%
\mathbf{B}$ and $\mathbf{A}\circ \mathbf{B}$ are equal. Moreover, the
related filtrations of the algebras $\mathbf{A}\circ _{\tau }\mathbf{B}$ and 
$\mathbf{A}\circ \mathbf{B}$ coincide.
\end{proposition}

\begin{proof}
As in the example above, let us \thinspace indicate by $\left\{
a_{i}\right\} $ and $\left\{ b_{i}\right\} $ a pair of basis for $\mathbf{A}%
_{1}$ and $\mathbf{B}_{1}$, respectively. The underlying vector space of $%
\mathbf{A}\circ _{\tau }\mathbf{B}\subset \mathbf{A}\otimes \mathbf{B}$ and $%
\mathbf{A}\circ \mathbf{B}\subset \mathbf{A}\otimes \mathbf{B}$ are
generated by the words written with elements $\left\{ a_{i}\otimes
b_{j}\right\} $ under $m_{A\otimes B}^{\tau }$ and $m_{A\otimes B}$. If we
write $a_{i}\otimes b_{j}\doteq ab$, $m_{A\otimes B}\doteq \cdot $ and $%
m_{A\otimes B}^{\tau }\doteq \cdot _{\tau }$, those words are related by 
\begin{equation}
ab\cdot _{\tau }ab\cdot _{\tau }...\cdot _{\tau }ab=\Gamma \;ab\cdot ab\cdot
...\cdot ab,  \label{relu}
\end{equation}
where $\Gamma $ is a matrix constructed from $\widehat{\tau }$ by using
successively Eq. $\left( \ref{mu}\right) $ (in any order, because both
products are associative). Because $\widehat{\tau }$ is invertible, $\Gamma $
is also invertible. Then, we can write $ab\cdot ...\cdot ab=\Gamma
^{-1}\;ab\cdot _{\tau }...\cdot _{\tau }ab$. Hence, any element of $\mathbf{A%
}\circ _{\tau }\mathbf{B}$ is an element of $\mathbf{A}\circ \mathbf{B}$,
and viceversa.\ 

The proof of the last statement is immediate from Eq. $\left( \ref{relu}%
\right) $.
\end{proof}

\begin{proposition}
Given $\mathcal{A}$ and $\mathcal{B}$ in $\mathrm{CA}$, a STTP $\mathcal{A}%
\circ _{\tau }\mathcal{B}$ is in $\mathrm{CA}$. Moreover, the related
gradations of $\mathcal{A}\circ _{\tau }\mathcal{B}$ and $\mathcal{A}\circ 
\mathcal{B}$ are equals.
\end{proposition}

\begin{proof}
The symmetric property of a twisting map $\tau :\mathbf{B}\otimes \mathbf{%
A\rightarrow A}\otimes \mathbf{B}$, where $\mathbf{A}$ and $\mathbf{B}$ are
graded algebras, insures that $\tau \left( \mathbf{B}_{n}\otimes \mathbf{A}%
_{m}\right) \subset \mathbf{A}_{m}\otimes \mathbf{B}_{n}$ (which can be seen
by direct calculation), and so the product $m_{A\otimes B}^{\tau }$
satisfies 
\begin{eqnarray*}
m_{A\otimes B}^{\tau }\left( \mathbf{A}_{n}\otimes \mathbf{B}_{n}\otimes 
\mathbf{A}_{m}\otimes \mathbf{B}_{m}\right) &=&\left( m_{A}\otimes
m_{B}\right) \left( \mathbf{A}_{n}\otimes \tau \left( \mathbf{B}_{n}\otimes 
\mathbf{A}_{m}\right) \otimes \mathbf{B}_{m}\right) \\
&\subset &\mathbf{A}_{n+m}\otimes \mathbf{B}_{n+m}.
\end{eqnarray*}
Thus, as a vector space $\mathbf{A}\circ _{\tau }\mathbf{B}=\bigoplus_{n}%
\mathbf{A}_{n}\otimes \mathbf{B}_{n}$ and 
\begin{equation*}
\left( \mathbf{A}_{n}\otimes \mathbf{B}_{n}\right) \cdot _{\tau }\left( 
\mathbf{A}_{m}\otimes \mathbf{B}_{m}\right) \subset \mathbf{A}_{n+m}\otimes 
\mathbf{B}_{n+m},
\end{equation*}
as we wanted to show.
\end{proof}

To conclude this section, let us observe that, in view of last proposition,
to define a STTP of conic quantum spaces $\mathcal{A}$ and $\mathcal{B}$ is
the same as fixing a particular algebra structure on the graded vector space 
$\mathbf{A}\circ \mathbf{B}$, such that $\mathbf{A}_{1}\otimes \mathbf{B}%
_{1} $ is the generating subspace. That means, given a symmetric twisting
map $\tau $, the conic quantum space $\mathcal{A}\circ _{\tau }\mathcal{B}$
(and in particular $\mathcal{A}\circ \mathcal{B}$) is an object of $\mathrm{%
GrVct}_{\Bbbk }$ (the category of graded vector spaces and homogeneous
linear maps) with an additional structure. Hence, defining on $\mathrm{GrVct}%
_{\Bbbk }$ the monoidal product 
\begin{equation*}
\mathbf{V}\circ \mathbf{W}=\bigoplus\nolimits_{n\in \mathbb{N}_{0}}\mathbf{V}%
_{n}\otimes \mathbf{W}_{n},
\end{equation*}
for $\mathbf{V}=\bigoplus_{n\in \mathbb{N}_{0}}\mathbf{V}_{n}$ and $\mathbf{W%
}=\bigoplus_{n\in \mathbb{N}_{0}}\mathbf{W}_{n}$, the forgetful functor $%
\frak{H}:\mathrm{CA}\hookrightarrow \mathrm{GrVct}_{\Bbbk }$ turns into a
monoidal one, and the equalities 
\begin{equation*}
\frak{H}\left( \mathcal{A}\circ _{\tau }\mathcal{B}\right) =\frak{H}\left( 
\mathcal{A}\circ \mathcal{B}\right) =\frak{H}\mathcal{A}\circ \frak{H}%
\mathcal{B}=\mathbf{A}\circ \mathbf{B}
\end{equation*}
hold for all twisting maps $\tau $.

\section{Twisted quantum matrix spaces}

In this section we finally study the consequences of replacing $\circ $ by $%
\circ _{\tau }$ in maps $\mathcal{A}\rightarrow \mathcal{H}\circ \mathcal{B}$
(the diagrams of $\left( \mathcal{A}\downarrow \mathrm{CA}\circ \mathcal{B}%
\right) $), obtaining in this way the categories $\Upsilon ^{\mathcal{A},%
\mathcal{B}}$ mentioned in \S \textbf{2.2}. Because $\mathrm{CA}\circ _{\tau
}\mathcal{B}$ is no longer a functor on \textrm{$CA$, }this change cannot be
performed in the framework of comma categories. However, one is able to
define arrows $\mathcal{A}\rightarrow \mathcal{H}\circ _{\tau }\mathcal{B}$
(which essentially give the objects defining $\Upsilon ^{\mathcal{A},%
\mathcal{B}}$) in terms of the morphisms $\frak{H}\mathcal{A}\rightarrow 
\frak{H}\mathcal{H}\circ \frak{H}\mathcal{B}$ in $\mathrm{GrVct}_{\Bbbk }$.
Hence, the comma categories\footnote{%
Its objects are pairs $\left\langle \varphi ,\mathcal{H}\right\rangle $
where $\mathcal{H}\in \mathrm{FGA}_{\mathrm{G}}$ and $\varphi $ is an arrow
in $\mathrm{GrVct}_{\Bbbk }$, 
\begin{equation*}
\varphi :\frak{H}\mathcal{A}\rightarrow \frak{H}\left( \mathcal{H}\circ 
\mathcal{B}\right) =\frak{H}\mathcal{H}\circ \frak{H}\mathcal{B}.
\end{equation*}
} $\left( \frak{H}\mathcal{A}\downarrow \frak{H}\left( \mathrm{CA}\circ 
\mathcal{B}\right) \right) $, where $\frak{H}\left( \mathrm{CA}\circ 
\mathcal{B}\right) $ is the composition of the functors $\mathrm{CA}\circ 
\mathcal{B}$ and $\frak{H}$, will be the cornerstone in the following
construction, since we will define the categories $\Upsilon ^{\mathcal{A},%
\mathcal{B}}$ as full subcategories of them.

\subsection{The categories $\Upsilon ^{\cdot }$}

As we made at the end of \S \textbf{1.1}, for every $\left\langle \varphi ,%
\mathcal{H}\right\rangle $ in $\left( \frak{H}\mathcal{A}\downarrow \frak{H}%
\left( \mathrm{CA}\circ \mathcal{B}\right) \right) $ it can be defined the
surjection $\pi ^{\varphi }:\mathbf{B}_{1}^{\ast }\otimes \mathbf{A}%
_{1}\twoheadrightarrow \mathbf{H}_{1}^{\varphi }$, such that$\;b^{j}\otimes
a_{i}\mapsto h_{i}^{j}\in \mathbf{H}_{1}$, and a related functor $\frak{F}%
:\left( \frak{H}\mathcal{A}\downarrow \frak{H}\left( \mathrm{CA}\circ 
\mathcal{B}\right) \right) \rightarrow \mathrm{CA}$ such that 
\begin{equation}
\frak{F}\left\langle \varphi ,\mathcal{H}\right\rangle =\mathcal{H}^{\varphi
}=\left( \mathbf{H}_{1}^{\varphi },\mathbf{H}^{\varphi }\right)
\;\;;\;\;\;\;\alpha \mapsto \left. \alpha \right| _{\mathbf{H}^{\varphi }},
\label{f}
\end{equation}
where $\mathbf{H}^{\varphi }$ is the subalgebra of $\mathbf{H}$ generated by 
$\mathbf{H}_{1}^{\varphi }$. (Note that in general, $\mathbf{H}^{\varphi }$
has not relation with the image of $\varphi $.)

Now, consider a linear bijection 
\begin{equation}
\widehat{\tau }:\mathbf{B}_{1}\otimes Lin\left[ \mathbf{B}_{1},\mathbf{A}_{1}%
\right] \backsimeq \mathbf{B}_{1}\otimes Lin\left[ \mathbf{B}_{1},\mathbf{A}%
_{1}\right] ,  \label{lb}
\end{equation}
associated to $\mathcal{A}$ and $\mathcal{B}$. Eventually, it could be
defined on $\mathbf{B}_{1}\otimes \mathbf{H}_{1}^{\varphi }$ through $\pi
^{\varphi }:Lin\left[ \mathbf{B}_{1},\mathbf{A}_{1}\right]
\twoheadrightarrow \mathbf{H}_{1}^{\varphi }$ in such a way that the diagram 
\begin{equation*}
\begin{diagram} \QTR{bf}{B}_{1}\otimes Lin\left[
\QTR{bf}{B}_{1},\QTR{bf}{A}_{1}\right] & \rOnto &\QTR{bf}{B}_{1}\otimes
\QTR{bf}{H}_{1}^{\varphi }\\ \dCorr~{\backsimeq } & & \dTo \\
\QTR{bf}{B}_{1}\otimes Lin\left[ \QTR{bf}{B}_{1},\QTR{bf}{A}_{1}\right] &
\rOnto &\QTR{bf}{B}_{1}\otimes \QTR{bf}{H}_{1}^{\varphi }\\ \end{diagram}
\end{equation*}
be commutative, and extend it to all of $\mathbf{B}\otimes \mathbf{H}%
^{\varphi }$ as a symmetric twisting map using the properties $\left( \ref
{mu}\right) $ and $\left( \ref{un}\right) $.

\begin{definition}
For every pair $\mathcal{A},\mathcal{B}\in \mathrm{CA}$ and every linear
bijection $\widehat{\tau }$ \emph{(}as in Eq. $\left( \ref{lb}\right) $\emph{%
)}, we define $\Upsilon ^{\mathcal{A},\mathcal{B}}$ as the full subcategory
of $\left( \frak{H}\mathcal{A}\downarrow \frak{H\,}\left( \mathrm{CA}\circ 
\mathcal{B}\right) \right) $, associated to $\widehat{\tau }$, formed out by
diagrams $\left\langle \varphi ,\mathcal{H}\right\rangle $ such that $%
\widehat{\tau }$ defines a symmetric twisting map for $\frak{F}\left\langle
\varphi ,\mathcal{H}\right\rangle =\mathcal{H}^{\varphi }$ and $\mathcal{B}$%
, and the homogeneous linear map $\varphi $ is a morphism of quantum spaces $%
\mathcal{A}\rightarrow \mathcal{H}^{\varphi }\circ _{\tau }\mathcal{B}$.\ \
\ $\blacksquare $
\end{definition}

Given now a collection $\left\{ \widehat{\tau }_{\mathcal{A},\mathcal{B}%
}\right\} _{\mathcal{A},\mathcal{B}\in \mathrm{CA}}$ of linear bijections,
we name $\Upsilon ^{\cdot }$ the disjoint union of the categories $\Upsilon
^{\mathcal{A},\mathcal{B}}$ defined above. Clearly, $\mathrm{CA}^{\circ }$
is a category $\Upsilon ^{\cdot }$ with associated collection 
\begin{equation*}
\widehat{\tau }_{\mathcal{A},\mathcal{B}}=I_{B_{1}\otimes Lin\left[
B_{1},A_{1}\right] },\;\;\forall \mathcal{A},\mathcal{B}\in \mathrm{CA}.
\end{equation*}
Moreover, calling $\frak{H}\mathrm{CA}^{\circ }$ the disjoint union of $%
\left( \frak{H}\mathcal{A}\downarrow \frak{H}\left( \mathrm{CA}\circ 
\mathcal{B}\right) \right) $, it follows that any $\Upsilon ^{\cdot }$ (in
particular $\mathrm{CA}^{\circ }$) is a full subcategory of $\frak{H}\mathrm{%
CA}^{\circ }$. This sets $\mathrm{CA}^{\circ }$ and a generic $\Upsilon
^{\cdot }$ on an equal footing.

\begin{theorem}
Each category $\Upsilon ^{\mathcal{A},\mathcal{B}}\subset \Upsilon ^{\cdot }$
has initial object.
\end{theorem}

\begin{proof}
Given $\mathcal{A}$ and $\mathcal{B}$ generated by linear spaces $\mathbf{A}%
_{1}$ and $\mathbf{B}_{1}$ with $\dim \mathbf{A}_{1}=\frak{n}$ and $\dim 
\mathbf{B}_{1}=\frak{m}$, the initial object of $\Upsilon ^{\mathcal{A},%
\mathcal{B}}$ is a pair 
\begin{equation*}
\left\langle \delta _{_{\mathcal{A},\mathcal{B}}},\underline{hom}^{\Upsilon }%
\left[ \mathcal{B},\mathcal{A}\right] \right\rangle ;\;\underline{hom}%
^{\Upsilon }\left[ \mathcal{B},\mathcal{A}\right] =\left( \mathbf{D}_{1},%
\mathbf{D}\right) ;\;\;\mathbf{D}_{1}\doteq span\left[ z_{i}^{j}\right]
_{i,j=1}^{\frak{n},\frak{m}}.
\end{equation*}
Also, $\delta _{_{\mathcal{A},\mathcal{B}}}=\delta $ is defined by the
injective linear map $\delta _{1}:a_{i}\mapsto z_{i}^{j}\otimes b_{j}$. Let
us show it.

Because of $Lin\left[ \mathbf{B}_{1},\mathbf{A}_{1}\right] \backsimeq 
\mathbf{D}_{1}$, $\widehat{\tau }_{\mathcal{A},\mathcal{B}}=\widehat{\tau }$
can be extended to all of $\mathbf{B}_{1}^{\otimes }\otimes \mathbf{D}%
_{1}^{\otimes }$ as a symmetric twisting map $\tau ^{\otimes }:\mathbf{B}%
_{1}^{\otimes }\otimes \mathbf{D}_{1}^{\otimes }\rightarrow \mathbf{D}%
_{1}^{\otimes }\otimes \mathbf{B}_{1}^{\otimes }$, and $\delta _{1}$
extended to an algebra homomorphism $\delta _{1}^{\otimes }:\mathbf{A}%
_{1}^{\otimes }\rightarrow \mathbf{D}_{1}^{\otimes }\otimes _{\tau ^{\otimes
}}\mathbf{B}_{1}^{\otimes }$. We are going to build up an algebra $\mathbf{D}
$ as a quotient of $\mathbf{D}_{1}^{\otimes }$, such that the linear map $%
\tau :\mathbf{B}\otimes \mathbf{D}\rightarrow \mathbf{D}\otimes \mathbf{B}$
and the algebra homomorphism $\mathbf{A}\rightarrow \mathbf{D}\circ _{\tau }%
\mathbf{B}$, given as the quotient maps of $\tau ^{\otimes }$ and $\delta
_{1}^{\otimes }$, respectively, are well-defined.

For $\mathbf{J}=\ker \left[ \mathbf{B}_{1}^{\otimes }\twoheadrightarrow 
\mathbf{B}\right] $ consider some complement $\mathbf{J}^{c}$, such that $%
\mathbf{J}\oplus \mathbf{J}^{c}=\mathbf{B}_{1}^{\otimes }$. Let $\left\{
I_{\lambda }\right\} _{\lambda \in \Lambda }$ be a set of linear generators
for $\mathbf{I}=\ker \left[ \mathbf{A}_{1}^{\otimes }\twoheadrightarrow 
\mathbf{A}\right] $ and \thinspace indicate by $\left\{ J_{\omega
}^{c}\right\} _{\omega \in \Omega }$ a basis for $\mathbf{J}^{c}$. Then, we
can write 
\begin{equation}
\delta _{1}^{\otimes }\left( I_{\lambda }\right) =d_{\lambda }^{\omega
}\otimes J_{\omega }^{c}+...\,,  \label{N}
\end{equation}
where $d_{\lambda }^{\omega }\in \mathbf{D}_{1}^{\otimes }$ are linear
combinations of monomials in $z_{i}^{j}$, and $``..."$ denotes terms
contained in $\mathbf{D}_{1}^{\otimes }\otimes \mathbf{J}$. On the other
hand, given a set of generators $\left\{ J_{\theta }\right\} _{\theta \in
\Theta }$ for $\mathbf{J}$, and basis $\left\{ b_{S}\right\} _{S\in \mathbf{S%
}}$ and $\left\{ z_{R}\right\} _{R\in \mathbf{R}}$ for $\mathbf{B}%
_{1}^{\otimes }$ and $\mathbf{D}_{1}^{\otimes }$, resp. (given the latter,
for instance, by monomials on $\mathbf{B}_{1}$ and $\mathbf{D}_{1}$ basis
elements, being $S$ and $R$ the obvious multi-indices), let us write 
\begin{equation}
\tau ^{\otimes }\left( J_{\theta }\otimes z_{R}\right) =d_{\theta
,R}^{\omega _{1}}\otimes J_{\omega _{1}}^{c}+...;\;\,\qquad \tau ^{\otimes
}\left( b_{S_{1}}\otimes d_{\lambda }^{\omega _{1}}\right) =d_{S_{1},\lambda
}^{\omega _{1},\omega _{2}}\otimes J_{\omega _{2}}^{c}+....\;  \label{M}
\end{equation}
Evaluating $\tau ^{\otimes }$ on $b_{S_{1}}\otimes d_{\theta ,R}^{\omega
_{1}}$ and $b_{S_{2}}\otimes d_{S_{1},\lambda }^{\omega _{1},\omega _{2}}$,
then evaluate the part of the result contained in $\mathbf{D}_{1}^{\otimes
}\otimes \mathbf{J}^{c}$, and so on, one arrives to (repeating the process $%
m $ times) 
\begin{equation*}
\tau ^{\otimes }\left( b_{S_{m-1}}\otimes d_{S_{1},...,S_{m-2},\theta
,R}^{\omega _{1},...,\omega _{m-1}}\right) =d_{S_{1},...,S_{m-1},\theta
,R}^{\omega _{1},...,\omega _{m}}\otimes J_{\omega _{m}}^{c}+...
\end{equation*}
and 
\begin{equation*}
\tau ^{\otimes }\left( b_{S_{m}}\otimes d_{S_{1},...,S_{m-1},\lambda
}^{\omega _{1},...,\omega _{m}}\right) =d_{S_{1},...,S_{m},\lambda }^{\omega
_{1},...,\omega _{m+1}}\otimes J_{\omega _{m+1}}^{c}+....
\end{equation*}
Defining in $\mathbf{D}_{1}^{\otimes }$ the bilateral ideal $\mathbf{L}$
algebraically generated by the set 
\begin{equation}
\left\{ \left\{ d_{S_{1},...,S_{m-1},\theta ,R}^{\omega _{1},...,\omega
_{m}}\right\} _{\theta \in \Theta }^{R\in \mathbf{R}}\bigcup \left\{
d_{S_{1},...,S_{m-1},\lambda }^{\omega _{1},...,\omega _{m}}\right\}
_{\lambda \in \Lambda }\right\} _{m\in \mathbb{N},S\in \mathbf{S}}^{\omega
\in \Omega },  \label{mi}
\end{equation}
(where for $m=1$ the subindexes $S$'s disappear) we get for the algebra $%
\mathbf{D}\doteq \left. \mathbf{D}_{1}^{\otimes }\right/ \mathbf{L}$ the
symmetric twisting map $\tau :\mathbf{B}\otimes \mathbf{D}\backsimeq \mathbf{%
D}\otimes \mathbf{B}$ and the algebra homomorphism $\delta :\mathbf{A}%
\rightarrow \mathbf{D}\circ _{\tau }\mathbf{B}$, given by the linear maps $%
\widehat{\tau }$ and $a_{i}\mapsto z_{i}^{j}\otimes b_{j}$, respectively.

From the gradation of $\mathbf{I}$ and $\mathbf{J}$, the ideal $\mathbf{L}$
inherits a graded structure. Thus $\mathcal{D}_{\mathcal{B},\mathcal{A}}=%
\mathcal{D}=\left( \mathbf{D}_{1},\mathbf{D}\right) \in \mathrm{CA}$ and $%
\delta $ defines the morphism of quantum spaces $\mathcal{A}\rightarrow 
\mathcal{D}\circ _{\tau }\mathcal{B}$ (in particular, an homogeneous linear
map). It is obvious that $\mathcal{D}=\frak{F}\left\langle \delta ,\mathcal{D%
}\right\rangle =\mathcal{D}^{\delta }$, so $\left\langle \delta ,\mathcal{D}%
\right\rangle $ is a diagram of $\Upsilon ^{\mathcal{A},\mathcal{B}}$. The
initial character of $\left\langle \delta ,\mathcal{D}\right\rangle $ is
immediate.
\end{proof}

Observe that above result does not involve the gradation of the algebras, so
it can be extended to the general case. Consider the category of pairs $%
\left( \mathbf{A}_{1},\mathbf{A}\right) $, with $\mathbf{A}\in \mathrm{Vct}%
_{\Bbbk }$ and $\mathbf{A}_{1}\subset \mathbf{A}$ a linear subspace of $%
\mathbf{A}$, and the related forgetful functor $\frak{\mathsf{h}}:\mathcal{A}%
\mapsto \left( \mathbf{A}_{1},\mathbf{A}\right) $. Let us construct the
comma categories $\left( \frak{\mathsf{h}}\mathcal{A}\downarrow \frak{%
\mathsf{h}}\left( \mathrm{FGA}\circ \mathcal{B}\right) \right) $ and their
related disjoint union $\frak{\mathsf{h}}\mathrm{FGA}^{\circ }$. Its objects
are diagrams $\left\langle \varphi ,\mathcal{H}\right\rangle $, with $%
\mathcal{H}\in \mathrm{FGA}$ and $\varphi $ a linear map $\mathbf{A}%
\rightarrow \mathbf{H}\circ \mathbf{B}$ satisfying $\varphi \left( \mathbf{A}%
_{1}\right) \subset \mathbf{H}_{1}\otimes \mathbf{B}_{1}$. So, we can define
a functor $\mathsf{f}:\frak{\mathsf{h}}\mathrm{FGA}^{\circ }\rightarrow 
\mathrm{FGA}$, as in Eq. $\left( \ref{f}\right) $, with $\mathsf{f}%
\left\langle \varphi ,\mathcal{H}\right\rangle =\mathcal{H}^{\varphi }$, and
define full subcategories $\Upsilon ^{\cdot }\subset \frak{\mathsf{h}}%
\mathrm{FGA}^{\circ }$ associated to collections $\left\{ \widehat{\tau }%
\right\} $, whose objects are diagrams $\left\langle \varphi ,\mathcal{H}%
\right\rangle \in \frak{\mathsf{h}}\mathrm{FGA}^{\circ }$ such that $\varphi 
$ gives the arrow $\mathcal{A}\rightarrow \mathsf{f}\left\langle \varphi ,%
\mathcal{H}\right\rangle \circ _{\tau }\mathcal{B}$. It can be shown
(directly\emph{\ }from the proof of the previous theorem) that each $%
\Upsilon ^{\mathcal{A},\mathcal{B}}$ has initial objects. However it is not
clear for us how to endow this $\Upsilon ^{\cdot }$ with a semigroupoid
structure, as we will do for the conic case.

\subsubsection{The initial objects of $\Upsilon ^{\mathcal{A},\mathcal{B}}$
and the geometric interpretation}

As at the end of \S \textbf{1.1}, for every initial object $\mathcal{D}_{%
\mathcal{B},\mathcal{A}}$ there is an epi $\mathcal{D}_{\mathcal{B},\mathcal{%
A}}\twoheadrightarrow \frak{F}\left\langle \varphi ,\mathcal{H}\right\rangle 
$, $\forall \left\langle \varphi ,\mathcal{H}\right\rangle \in \Upsilon ^{%
\mathcal{A},\mathcal{B}}$. Thus, the opposite of any $\frak{F}\left\langle
\varphi ,\mathcal{H}\right\rangle $ can be regarded as a \emph{non
commutative subvariety }$\frak{F}\left\langle \varphi ,\mathcal{H}%
\right\rangle ^{op}\hookrightarrow \mathcal{D}_{\mathcal{B},\mathcal{A}%
}^{op} $.

Following \cite{man0}, we develope this idea in the terminology of $\mathbf{C%
}$-points. Fix a couple of quantum spaces $\mathcal{A},\mathcal{B}$ and
consider the initial object $\mathcal{D}_{\mathcal{B},\mathcal{A}}=\mathcal{D%
}$ of $\Upsilon ^{\mathcal{A},\mathcal{B}}$. Suppose an algebra $\mathbf{C}$
and a pair of $\mathbf{C}$-points $i_{B}$ and $i_{D}$, of $\mathbf{B}$ and $%
\mathbf{D}$, are given, in such a way that 
\begin{equation*}
i_{B}\left( b_{i}\right) \cdot i_{D}\left( z_{k}^{l}\right) =\widehat{\tau }%
_{ink}^{jlm}\;i_{D}\left( z_{m}^{n}\right) \cdot i_{B}\left( b_{j}\right) .
\end{equation*}
We are denoting by $``\cdot "$ the product of $\mathbf{C}$. If we write $%
i_{B}\left( b_{i}\right) \doteq \mathbf{b}_{i}$ and $i_{D}\left(
z_{k}^{l}\right) \doteq \mathbf{z}_{k}^{l}$, then the equation above says 
\begin{equation}
\mathbf{b}_{i}\cdot \mathbf{z}_{k}^{l}=\widehat{\tau }_{ink}^{jlm}\;\mathbf{z%
}_{m}^{n}\cdot \mathbf{b}_{j}.  \label{2i}
\end{equation}
That means $i_{B}$ and $i_{D}$ are non commuting $\mathbf{C}$-points, and
the non commutativity is given by a symmetric twisting map defining $%
\Upsilon ^{\mathcal{A},\mathcal{B}}$. Since the last theorem the map $%
i_{A}\doteq \left( i_{D}\cdot i_{B}\right) \,\delta $, given by 
\begin{equation*}
i_{A}\left( a_{j}\right) \doteq \mathbf{a}_{j}\doteq \mathbf{z}_{j}^{k}\cdot 
\mathbf{b}_{k}\in \mathbf{C},
\end{equation*}
defines a $\mathbf{C}$-point of $\mathbf{A}$. Moreover,

\begin{theorem}
Consider the quantum spaces $\left( \mathbf{A}_{1},\mathbf{A}\right) $ and $%
\left( \mathbf{H}_{1},\mathbf{H}\right) $, with 
\begin{equation*}
\mathbf{H}_{1}=span\left[ h_{j}^{k}\right] _{i,j=1}^{\frak{n},\frak{m}%
},\;h_{j}^{k}\in \mathbf{H},
\end{equation*}
an algebra $\mathbf{C}$, and $\mathbf{C}$-points $i_{B}:b_{j}\mapsto \mathbf{%
b}_{j}$ and $i_{H}:h_{j}^{k}\mapsto \mathbf{h}_{j}^{k}$ satisfying Eq. $%
\left( \ref{2i}\right) $ \emph{(}changing $\mathbf{z}$ by $\mathbf{h}$\emph{)%
}. Suppose in addition that $i_{B}$ is a generic point. Then, the map $%
a_{j}\in \mathbf{A}_{1}\mapsto \mathbf{h}_{j}^{k}\cdot \mathbf{b}_{k}\in 
\mathbf{C}$ can be extended to a $\mathbf{C}$-point of $\mathbf{A}$ \emph{iff%
} the assignment $z_{j}^{k}\mapsto \mathbf{h}_{j}^{k}$ defines a $\mathbf{C}$%
-point of $\mathbf{D}$.
\end{theorem}

\begin{proof}
Suppose $a_{j}\mapsto \mathbf{h}_{j}^{k}\cdot \mathbf{b}_{k}$ can be
extended to a $\mathbf{C}$-point of $\mathbf{A}$, and call $i_{A}$ such a
map. Let us come back to the proof of previous theorem. There, for instance,
we denote by $I_{\lambda }$ an element of the ideal $\mathbf{I}$ related to $%
\mathbf{A}$. Here, we will identify $I_{\lambda }$ with its projection (by $%
\Pi $) over $\mathbf{A}$, so $I_{\lambda }=0$. Applying an analogous
criterium to the other symbols, we have from Eq. $\left( \ref{N}\right) $
that 
\begin{equation*}
0=i_{A}\left( I_{\lambda }\right) =i_{H}\left( d_{\lambda }^{\omega }\right)
\cdot i_{B}\left( J_{\omega }^{c}\right) .
\end{equation*}
Now, the symbols $d_{\lambda }^{\omega }\in \mathbf{H}$ represent linear
combinations of monomial in $h_{j}^{k}$. They are formally identical to the
ones in $\mathbf{D}$ (we just have to replace $z$ by $h$). Since $i_{B}$ is
a generic $\mathbf{C}$-point, the elements $i_{B}\left( J_{\omega
}^{c}\right) $ are linearly independent, and consequently each $i_{H}\left(
d_{\lambda }^{\omega }\right) $ must be the null element. Similarly, from $%
\left( \ref{M}\right) $ follows that 
\begin{equation*}
0=i_{B}\left( J_{\theta }\right) \cdot i_{H}\left( h_{R}\right) =i_{H}\left(
d_{\theta ,R}^{\omega _{1}}\right) \cdot i_{B}\left( J_{\omega
_{1}}^{c}\right) ,
\end{equation*}
and for the same reason, $i_{H}\left( d_{\theta ,R}^{\omega _{1}}\right) =0$%
. The equation above can be seen as the result of passing $i_{H}\left(
h_{R}\right) $ through $i_{B}\left( J_{\theta }\right) $ using Eq. $\left( 
\ref{2i}\right) $ (replacing $\mathbf{z}$ by $\mathbf{h}$). Following that
process, we finally arrive at the equations (see Eq. $\left( \ref{mi}\right) 
$) 
\begin{equation}
i_{H}\left( d_{S_{1},...,S_{m-1},\theta ,R}^{\omega _{1},...,\omega
_{m}}\right) =0\;\;and\;\;i_{H}\left( d_{S_{1},...,S_{m-1},\lambda }^{\omega
_{1},...,\omega _{m}}\right) =0.  \label{mii}
\end{equation}
That means, the elements $i_{H}\left( h_{j}^{k}\right) =\mathbf{h}%
_{j}^{k}\in \mathbf{C}$ must satisfy unless the same relations that have to
satisfy the elements $z_{j}^{k}\in \mathbf{D}$. Thus, $z_{j}^{k}\mapsto 
\mathbf{h}_{j}^{k}$ extends to an algebra homomorphism.

Reciprocally, if $z_{j}^{k}\mapsto \mathbf{h}_{j}^{k}$ extends to an algebra
homomorphism, then the elements $\mathbf{h}_{j}^{k}$ must satisfy unless the
relations $\left( \ref{mii}\right) $, what implies that $a_{j}\mapsto 
\mathbf{h}_{j}^{k}\cdot \mathbf{b}_{k}$ is a well defined $\mathbf{C}$-point.
\end{proof}

Thus, this theorem is simply another way to express the initial property of $%
\mathcal{D}$ in the category $\Upsilon ^{\mathcal{A},\mathcal{B}}$. Note
that again, the gradation has not been used.

\subsubsection{Factorizable bijections}

We can see from Equation $\left( \ref{mi}\right) $ or $\left( \ref{mii}%
\right) $ that the algebras $\mathbf{D}$ could be very \emph{small}. Things
change by considering special classes of twisting maps.

We say that a collection $\left\{ \widehat{\tau }_{\mathcal{A},\mathcal{B}%
}\right\} _{\mathcal{A}\in \mathrm{CA}}$ is \emph{factorizable} if there
exists another collection 
\begin{equation*}
\left\{ \sigma _{\mathcal{A}}\right\} _{\mathcal{A}\in \mathrm{CA}%
},\;\;\sigma _{\mathcal{A}}\in Aut\left[ \mathbf{A}_{1}\right] ,
\end{equation*}
such that 
\begin{equation*}
\widehat{\tau }_{\mathcal{A},\mathcal{B}}=id\otimes \sigma _{\mathcal{B}%
}^{\ast -1}\otimes \sigma _{\mathcal{A}}:\mathbf{B}_{1}\otimes \mathbf{B}%
_{1}^{\ast }\otimes \mathbf{A}_{1}\backsimeq \mathbf{B}_{1}\otimes \mathbf{B}%
_{1}^{\ast }\otimes \mathbf{A}_{1}.
\end{equation*}
In that case we say $\widehat{\tau }_{\mathcal{A},\mathcal{B}}$ is a
factorizable bijection.

As an example, consider $\mathcal{A}$ and $\mathcal{B}$ with $\mathbf{A}%
\backsimeq \left. \mathbf{A}_{1}^{\otimes }\right/ \mathbf{I}$ and $\mathbf{B%
}\backsimeq \left. \mathbf{B}_{1}^{\otimes }\right/ \mathbf{J}$, and ideals
given by Eqs. $\left( \ref{ideals}\right) $ of \S \textbf{1.3}. Let us take $%
\Upsilon ^{\mathcal{A},\mathcal{B}}\mathrm{\ }$with associated factorizable
bijection $\widehat{\tau }_{\mathcal{A},\mathcal{B}}$ as above. The algebra $%
\mathbf{D}$ will be 
\begin{equation}
\mathbf{D}=\left. \left[ span\left[ z_{i}^{j}\right] _{i,j=1}^{\frak{n},%
\frak{m}}\right] ^{\otimes }\right/ \mathbf{L},  \label{A}
\end{equation}
being $\mathbf{L}$ the ideal algebraically generated by the elements
(compare to Eq. $\left( \ref{ideal}\right) $) 
\begin{equation}
R_{\lambda _{n}}^{k_{1}...k_{n}}\;\check{z}_{k_{1}}^{\left\langle
m-1\right\rangle j_{1}}\cdot \check{z}_{k_{2}}^{\left\langle m\right\rangle
j_{2}}\;...\;\check{z}_{k_{n}}^{\left\langle m+n-2\right\rangle
j_{n}}\;\left( S^{\bot }\right) _{j_{1}...j_{n}}^{\omega _{n}},\;\;\;\check{z%
}_{u}^{\left\langle r\right\rangle v}=\left[ \sigma _{\mathcal{A}}^{r}\right]
_{u}^{p}\,z_{p}^{q}\,\left[ \sigma _{\mathcal{B}}^{-r}\right] _{q}^{v},
\label{B}
\end{equation}
$\lambda _{n}\in \Lambda _{n}$, $\omega _{n}\in \Omega _{n}$, $m\in \mathbb{N%
}$, $n\in \mathbb{N}_{0}$. In addition, if $\sigma _{\mathcal{A}}$ and $%
\sigma _{\mathcal{B}}$ can be extended to $\mathbf{A}$ and $\mathbf{B}$ as
algebras automorphisms, the ideal defined by Eq. $\left( \ref{B}\right) $
reduces to the one defined by the elements 
\begin{equation*}
R_{\lambda _{n}}^{k_{1}...k_{n}}\;z_{k_{1}}^{j_{1}}\cdot \check{z}%
_{k_{2}}^{j_{2}}\cdot \check{z}_{k_{3}}^{\left\langle 2\right\rangle
j_{3}}\;...\;\check{z}_{k_{n}}^{\left\langle n-1\right\rangle j_{n}}\;\left(
S^{\bot }\right) _{j_{1}...j_{n}}^{\omega _{n}},
\end{equation*}
$\lambda _{n}\in \Lambda _{n}$, $\omega _{n}\in \Omega _{n}$, $n\in \mathbb{N%
}_{0}$, as we will show in \S \textbf{4.2.2}.

\subsection{The semigroupoid structure of $\Upsilon ^{\cdot }$}

We begin this section with the following observation: $\frak{H}\mathrm{CA}%
^{\circ }$ has a semigroupoid structure given by the partial product functor 
\begin{equation}
\left\langle \varphi ,\mathcal{H}\right\rangle \times \left\langle \chi ,%
\mathcal{G}\right\rangle \mapsto \left\langle \left( I_{H}\circ \chi \right)
\,\varphi ,\mathcal{H}\circ \mathcal{G}\right\rangle ;\;\;\alpha \times
\beta \mapsto \alpha \circ \beta ,  \label{ma}
\end{equation}
with domain 
\begin{equation*}
\bigvee_{_{\mathcal{A},\mathcal{B},\mathcal{C}\in \mathrm{CA}}}\left( \frak{H%
}\mathcal{A}\downarrow \frak{H}\left( \mathrm{CA}\circ \mathcal{C}\right)
\right) \times \left( \frak{H}\mathcal{C}\downarrow \frak{H}\left( \mathrm{CA%
}\circ \mathcal{B}\right) \right) ,
\end{equation*}
and codomain $\frak{H}\mathrm{CA}^{\circ }$, and $\mathrm{CA}^{\circ
}\subset \frak{H}\mathrm{CA}^{\circ }$ is a\textbf{\ }sub-semigroupoid. Its
associativity comes from that of $\circ $, and the unit elements are given
by the diagrams $\left\langle \ell _{\mathcal{A}},\mathcal{K}\right\rangle $%
, where $\ell _{\mathcal{A}}$ is the homogeneous isomorphism $\mathbf{A}%
\backsimeq \Bbbk \left[ e\right] \otimes \mathbf{A}$, such that $a\mapsto
e^{n}\otimes a$ if $a\in \mathbf{A}_{n}$. Finally, it is immediate from Eq. $%
\left( \ref{dok}\right) $ that $\mathrm{CA}^{\circ }$ is a sub-semigroupoid
of $\frak{H}\mathrm{CA}^{\circ }$.$\smallskip $

Nevertheless, for a generic collection of linear bijections, $\Upsilon
^{\cdot }$ fails to be a semigroupoid. To address this problem, we shall
consider particular collections of twisting maps.

Given a pair of linear endomorphisms $\alpha $ and $\beta $ on vector spaces 
$\otimes _{i\in I}\mathbf{V}_{i}$ and $\otimes _{i\in I}\mathbf{W}_{i}$ ($I$
a set of indices), respectively, we shall say that $\alpha $ is \emph{%
strongly congruent}\textbf{\ }to $\beta $, and note $\alpha \backsim \beta $%
, if there exist a family of isomorphisms $\zeta _{i}:\mathbf{V}%
_{i}\backsimeq \mathbf{W}_{i}$, $\forall i\in I$, such that 
\begin{equation*}
\alpha =\left( \otimes _{i\in I}\;\zeta _{i}^{-1}\right) \;\beta \;\left(
\otimes _{i\in I}\;\zeta _{i}\right) .
\end{equation*}
A collection $\left\{ \widehat{\tau }_{\mathcal{A},\mathcal{B}}\right\} _{%
\mathcal{A},\mathcal{B}\in \mathrm{CA}}$ of linear bijections $\left( \ref
{lb}\right) $ will be called \emph{global collection} when 
\begin{equation*}
\widehat{\tau }_{\mathcal{A},\mathcal{B}}\backsim \widehat{\tau }_{\mathcal{C%
},\mathcal{D}}\;\;\;\Leftrightarrow \;\;\;\mathbf{A}_{1}\backsimeq \mathbf{C}%
_{1}\;\;and\;\;\mathbf{B}_{1}\backsimeq \mathbf{D}_{1}.
\end{equation*}

Note that two $\widehat{\tau }$'s are strongly congruent \emph{iff} there
exists basis on the respective linear spaces such that the matrices defined
by the maps $\widehat{\tau }$'s on these basis coincide. Hence, the
equivalence classes can be represented by a collection $\left\{ \tau _{\frak{%
n},\frak{m}}\right\} _{\frak{n},\frak{m}\in \mathbb{N}_{0}}$ of invertible
matrices $\tau _{\frak{n},\frak{m}}\in GL\left( \frak{m}^{2}\cdot \frak{n}%
\right) $ or linear isomorphisms 
\begin{equation*}
\tau _{\frak{n},\frak{m}}:\Bbbk ^{\frak{m}}\otimes \Bbbk ^{\ast \frak{m}%
}\otimes \Bbbk ^{\frak{n}}\backsimeq \Bbbk ^{\frak{m}}\otimes \Bbbk ^{\ast 
\frak{m}}\otimes \Bbbk ^{\frak{n}}.
\end{equation*}
In particular, for a factorizable $\tau _{\frak{n},\frak{m}}$ we have $\tau
_{\frak{n},\frak{m}}=id_{\frak{m}}\otimes \sigma _{\frak{m}}^{\ast
-1}\otimes \sigma _{\frak{n}}$, with $\sigma _{\frak{n}}\in GL\left( \frak{n}%
\right) $ or $\sigma _{\frak{n}}:\Bbbk ^{\frak{n}}\backsimeq \Bbbk ^{\frak{n}%
}$; in indices 
\begin{equation*}
\left[ \widehat{\tau }_{\frak{n},\frak{m}}\right] _{i\alpha k}^{jl\beta
}=\delta _{i}^{j}\,\left[ \sigma _{\frak{m}}^{-1}\right] _{\alpha }^{l}\,%
\left[ \sigma _{\frak{n}}\right] _{k}^{\beta }.
\end{equation*}
In these cases the equivalence classes are given by a collection $\left\{
\sigma _{\frak{n}}\right\} _{\frak{n}\in \mathbb{N}_{0}}$.

In working with a global collection, it can be regarded as the replacement
of the category $\mathrm{FGA}$ by another one formed out by pairs $\left( 
\frak{n},\mathbf{I}\right) $, where $\frak{n}\in \mathbb{N}_{0}$ and $%
\mathbf{I}\ $is a bilateral ideal of $\left[ \Bbbk ^{\frak{n}}\right]
^{\otimes }$, with arrows $\left( \frak{n},\mathbf{I}\right) \rightarrow
\left( \frak{m},\mathbf{J}\right) $ given by linear maps $\alpha :\Bbbk ^{%
\frak{n}}\rightarrow \Bbbk ^{\frak{m}}$ defining a inclusion $\alpha
^{\otimes }\left( \mathbf{I}\right) \subset \mathbf{J}$.

\subsubsection{The $\mathrm{CA}$-cobased categories}

Now, let us present main result of this section.

\begin{theorem}
If the category $\Upsilon ^{\cdot }$ is associated to a global collection,
the following statements are equivalent:

a\emph{) }$\Upsilon ^{\cdot }\subset \frak{H}\mathrm{CA}^{\circ }$ is a
sub-semigroupoid.

b\emph{) }The collection that defines $\Upsilon ^{\cdot }$ is factorizable.
\end{theorem}

\begin{proof}
$\bullet \;a)\Rightarrow b)$ Consider $\mathcal{A}$ with $\dim \mathbf{A}%
_{1}=\frak{n}$. If every $\Upsilon ^{\mathcal{A}}$ is a submonoid of $\left( 
\frak{H}\mathcal{A}\downarrow \frak{H}\left( \mathrm{CA}\circ \mathcal{A}%
\right) \right) $, then $\left\langle \left( I\circ \delta \right) \,\delta ,%
\mathcal{D}\circ \mathcal{D}\right\rangle \in \Upsilon ^{\mathcal{A}}$ and
in consequence $\left( I\circ \delta \right) \,\delta $ defines the algebra
morphism 
\begin{equation}
\left( I\circ \delta \right) \,\delta :\mathbf{A}\rightarrow \left( \mathbf{D%
}\circ \mathbf{D}\right) \otimes _{\tau }\mathbf{A}  \label{eq}
\end{equation}
with $\;\left( I\circ \delta \right) \,\delta \left( a_{m}\right)
=z_{m}^{r}\otimes z_{r}^{l}\otimes a_{l}$. The twisting map $\tau $, given
by a matrix $\tau _{\frak{n},\frak{n}}=\tau _{\frak{n}}$, is 
\begin{equation}
\widehat{\tau }_{\mathcal{A}}\left[ a_{i}\otimes \left( z_{m}^{r}\otimes
z_{r}^{l}\right) \right] =\left[ \tau _{\frak{n}}\right] _{ij\beta
}^{kl\alpha }\;a_{k}\otimes \left( z_{\alpha }^{r}\otimes z_{r}^{j}\right) .
\label{eq1}
\end{equation}
Note that $\tau $ is only defined on $\frak{F}\left\langle \left( I\circ
\delta \right) \,\delta ,\mathcal{D}\circ \mathcal{D}\right\rangle $, given
the latter by the subalgebra of $\mathbf{D}\circ \mathbf{D}$ generated by 
\begin{equation}
span\left[ \sum\nolimits_{j=1}^{\frak{n}}z_{i}^{j}\otimes z_{j}^{k}\right]
_{i,k=1}^{\frak{n}}\subset \mathbf{D}_{1}\otimes \mathbf{D}_{1}.  \label{acl}
\end{equation}
Since Eq. $\left( \ref{eq}\right) $ and $\left( \ref{eq1}\right) $, $\left(
I\circ \delta \right) \,\delta $ applied to $a_{i}\otimes a_{j}$ gives 
\begin{equation*}
\left[ \left[ \tau _{\frak{n}}\right] _{kyj}^{slx}\,\left[ \tau _{\frak{n}}%
\right] _{\alpha pl}^{r\beta q}-\delta _{k}^{s}\delta _{y}^{q}\,\left[ \tau
_{\frak{n}}\right] _{\alpha pj}^{r\beta x}\right] \left( z_{i}^{k}\cdot
z_{x}^{y}\otimes z_{s}^{\alpha }\cdot z_{q}^{p}\otimes a_{r}\cdot a_{\beta
}\right) =0.
\end{equation*}
Because this eq. must hold for all$\frak{\ }\mathcal{A}\in \mathrm{CA}$, in
particular for those $\mathcal{A}$ defined by free algebras, 
\begin{equation}
\left[ \tau _{\frak{n}}\right] _{kyj}^{slx}\,\left[ \tau _{\frak{n}}\right]
_{\alpha pl}^{r\beta q}=\delta _{k}^{s}\delta _{y}^{q}\,\left[ \tau _{\frak{n%
}}\right] _{\alpha pj}^{r\beta x},\;\forall \frak{n}\in \mathbb{N}_{0},
\label{coas}
\end{equation}
is necessary from the global character of the considered collection. On the
other hand, due to $\Upsilon ^{\mathcal{A}}$ is a submonoid, the unit $%
\left\langle \ell _{\mathcal{A}},\mathcal{K}\right\rangle $ of $\left( \frak{%
H}\mathcal{A}\downarrow \frak{H}\left( \mathrm{CA}\circ \mathcal{A}\right)
\right) $, is also a unit for $\Upsilon ^{\mathcal{A}}$. In particular, $%
\ell _{\mathcal{A}}$ must define the isomorphisms $\mathcal{A}\backsimeq 
\mathcal{K}\circ \mathcal{A}$ and $\mathcal{A}\backsimeq \mathcal{K}\circ
_{\tau }\mathcal{A}$ (for all $\mathcal{A}$ with $\dim \mathbf{A}_{1}=\frak{n%
}$). Since the map $\widehat{\tau }_{\mathcal{A}}$ for $\left\langle \ell _{%
\mathcal{A}},\mathcal{K}\right\rangle $ is given by an isomorphism $\Bbbk
\otimes \mathbf{A}_{1}\backsimeq \Bbbk \otimes \mathbf{A}_{1}$ such that 
\begin{equation*}
a_{i}\otimes e\,\delta _{\beta }^{l}\mapsto \left[ \tau _{\frak{n}}\right]
_{ij\beta }^{kl\alpha }\;a_{k}\otimes e\,\delta _{\alpha }^{j}=\left[ \tau _{%
\frak{n}}\right] _{i\alpha \beta }^{kl\alpha }\;a_{k}\otimes e,
\end{equation*}
it is easy to see that $\left[ \tau _{\frak{n}}\right] _{i\alpha \beta
}^{kl\alpha }=\delta _{\beta }^{l}\tau _{i}^{k}$ with $\widehat{\tau }_{%
\mathcal{A}}\left( a_{i}\otimes e\right) =\tau _{i}^{k}\;a_{k}\otimes e$.
This immediately implies that $\ell _{\mathcal{A}}$ defines the above
isomorphisms \emph{iff} 
\begin{equation}
\left[ \tau _{\frak{n}}\right] _{i\alpha \beta }^{kl\alpha }=\delta
_{i}^{k}\delta _{\beta }^{l}.  \label{cou}
\end{equation}
The solutions to Eqs. $\left( \ref{coas}\right) $ and $\left( \ref{cou}%
\right) $ (see lemma at the end of this section) are 
\begin{equation*}
\tau _{\frak{n}}=id_{\frak{n}}\otimes \overline{\sigma }_{\frak{n}}\otimes
\sigma _{\frak{n}}\;\;with\;\;\sigma _{\frak{n}}\in GL\left( \frak{n}\right)
\diagup \mathcal{Z}_{GL\left( \frak{n}\right) },
\end{equation*}
or in indices $\left[ \tau _{\frak{n}}\right] _{ij\beta }^{kl\alpha }=\delta
_{i}^{k}\,\left[ \overline{\sigma }_{\frak{n}}\right] _{j}^{l}\,\left[
\sigma _{\frak{n}}\right] _{\beta }^{\alpha }$, being $\overline{\sigma }%
=\sigma ^{-1}$.

We now address the general case. Let $\mathcal{B}$ and $\mathcal{C}$ be
quantum spaces, and basis $\left\{ b_{i}\right\} _{i=1}^{\frak{m}}$ and $%
\left\{ c_{i}\right\} _{i=1}^{\frak{p}}$ of $\mathbf{B}_{1}$ and $\mathbf{C}%
_{1}$. We shall note $\left\langle \delta ,\mathcal{D}\right\rangle $ and $%
\left\langle \delta ^{\prime },\mathcal{D}^{\prime }\right\rangle $ the
corresponding initial objects of $\Upsilon ^{\mathcal{A},\mathcal{B}}$ and $%
\Upsilon ^{\mathcal{B},\mathcal{C}}$, respectively. (Again, the indices will
run on the range that correspond to their associated dimensions.) Because of 
\begin{equation*}
\left\langle \left( I\circ \delta ^{\prime }\right) \,\delta ,\mathcal{D}%
\circ \mathcal{D}^{\prime }\right\rangle \in \Upsilon ^{\mathcal{A},\mathcal{%
C}},
\end{equation*}
the map $\left( I\circ \delta ^{\prime }\right) \,\delta $ gives rise to the
algebra morphism 
\begin{equation}
\mathbf{A}\rightarrow \left( \mathbf{D}\circ \mathbf{D}^{\prime }\right)
\otimes _{\tau }\mathbf{C\;}:\;a_{m}\mapsto z_{m}^{r}\otimes \left(
z^{\prime }\right) _{r}^{l}\otimes c_{l},  \label{equ}
\end{equation}
with twisting map $\tau $ given by (it is valid the same observation that
leads to Eq. $\left( \ref{acl}\right) $) 
\begin{equation*}
\widehat{\tau }_{\mathcal{A},\mathcal{C}}\left[ c_{i}\otimes \left( z_{\beta
}^{r}\otimes \left( z^{\prime }\right) _{r}^{l}\right) \right] =\left[ \tau
_{\frak{n},\frak{p}}\right] _{ij\beta }^{kl\alpha }\;c_{k}\otimes \left(
z_{\alpha }^{r}\otimes \left( z^{\prime }\right) _{r}^{j}\right) .
\end{equation*}
To be $\left( I\circ \delta ^{\prime }\right) \,\delta $ an algebra morphism
defined by $\left( \ref{equ}\right) $%
\begin{equation*}
\left[ \left[ \tau _{\frak{n},\frak{m}}\right] _{kyj}^{slx}\,\left[ \tau _{%
\frak{m},\frak{p}}\right] _{\alpha pl}^{r\beta q}-\delta _{k}^{s}\delta
_{y}^{q}\,\left[ \tau _{\frak{n},\frak{p}}\right] _{\alpha pj}^{r\beta x}%
\right] \left( z_{i}^{k}\cdot \left( z^{\prime }\right) _{x}^{y}\otimes
z_{s}^{\alpha }\cdot \left( z^{\prime }\right) _{q}^{p}\otimes c_{r}\cdot
c_{\beta }\right) =0
\end{equation*}
is necessary. Asking that the above eq. holds for all $\mathcal{A},\mathcal{B%
},\mathcal{C}\in \mathrm{CA}$ with the associated dimensions $\frak{n},\frak{%
m}$ and $\frak{p}$, we need 
\begin{equation*}
\left[ \tau _{\frak{n},\frak{m}}\right] _{kyj}^{slx}\,\left[ \tau _{\frak{m},%
\frak{p}}\right] _{\alpha pl}^{r\beta q}=\delta _{k}^{s}\delta _{y}^{q}\,%
\left[ \tau _{\frak{n},\frak{p}}\right] _{\alpha pj}^{r\beta x},\;\forall 
\frak{n},\frak{m},\frak{p}\in \mathbb{N}_{0}.
\end{equation*}
We can solve these equations in terms of the solution of the case $\frak{n}=%
\frak{m}=\frak{p}$. Taking $\frak{n}=\frak{p}$, the last eq. can be written 
\begin{equation}
\left[ \tau _{\frak{n},\frak{m}}\right] _{kyj}^{slx}\,\left[ \tau _{\frak{m},%
\frak{n}}\right] _{\alpha pl}^{r\beta q}=\delta _{k}^{s}\delta _{y}^{q}\,%
\left[ \tau _{\frak{n}}\right] _{\alpha pj}^{r\beta x}=\delta _{k}^{s}\delta
_{y}^{q}\delta _{\alpha }^{r}\,\left[ \overline{\sigma }_{\frak{n}}\right]
_{p}^{\beta }\,\left[ \sigma _{\frak{n}}\right] _{j}^{x},  \label{inv}
\end{equation}
and contracting $x$ with $p$ and $\alpha $ with $s$, it reduces to 
\begin{equation}
\left[ \tau _{\frak{n},\frak{m}}\right] _{kyj}^{\alpha lp}\,\left[ \tau _{%
\frak{m},\frak{n}}\right] _{\alpha pl}^{r\beta q}=\delta _{k}^{r}\delta
_{y}^{q}\delta _{j}^{\beta }.  \label{inve}
\end{equation}
Changing $\frak{n}$ by $\frak{m}$ (and reordering indexes) it is obvious
that $\tau _{\frak{n},\frak{m}}$ and $\tau _{\frak{m},\frak{n}}$ are
mutually inverses. Therefore, multiplying $\left( \ref{inv}\right) $ to the
right by $\left[ \tau _{\frak{n},\frak{m}}\right] _{rq\beta }^{uwv}$ and
contracting $r,q,\beta $, 
\begin{equation*}
\left[ \tau _{\frak{n},\frak{m}}\right] _{kyj}^{swx}\delta _{\alpha
}^{u}\delta _{p}^{v}=\delta _{k}^{s}\,\left[ \sigma _{\frak{n}}\right]
_{j}^{x}\left[ \left[ \overline{\sigma }_{\frak{n}}\right] _{p}^{\beta }%
\left[ \tau _{\frak{n},\frak{m}}\right] _{\alpha y\beta }^{uwv}\right]
\end{equation*}
is fulfilled. Contracting $u$ with $\alpha $ and $p$ with $v$, and defining 
\begin{equation*}
\left[ \beta _{\frak{m}}\right] _{y}^{w}\doteq \frac{1}{\frak{m}^{2}}\left[ 
\overline{\sigma }_{\frak{n}}\right] _{p}^{\beta }\,\left[ \tau _{\frak{n},%
\frak{m}}\right] _{\alpha y\beta }^{\alpha wp},
\end{equation*}
it sees that $\left[ \tau _{\frak{n},\frak{m}}\right] _{kyj}^{swx}=\delta
_{k}^{s}\,\left[ \beta _{\frak{m}}\right] _{y}^{w}\,\left[ \sigma _{\frak{n}}%
\right] _{j}^{x}$; and changing $\frak{n}$ by $\frak{m}$, 
\begin{equation*}
\left[ \tau _{\frak{m},\frak{n}}\right] _{kyj}^{swx}=\delta _{k}^{s}\,\left[
\beta _{\frak{n}}\right] _{y}^{w}\,\left[ \sigma _{\frak{m}}\right] _{j}^{x}.
\end{equation*}
Now, substituting the last two expression on $\left( \ref{inve}\right) $, it
is straightforward to show that $\beta _{\frak{n}}=\overline{\sigma }_{\frak{%
n}}$ and $\beta _{\frak{m}}=\overline{\sigma }_{\frak{m}}$; therefore 
\begin{equation*}
\tau _{\frak{n},\frak{m}}=id_{\frak{m}}\otimes \overline{\sigma }_{\frak{m}%
}\otimes \sigma _{\frak{n}},\;\sigma _{\frak{n}}\in GL\left( \frak{n}\right)
,\;\sigma _{\frak{m}}\in GL\left( \frak{m}\right) .
\end{equation*}
$\bullet \;b)\Rightarrow a)$ Let us suppose that $\Upsilon ^{\cdot }$ is
defined by a factorizable collection of bijections. Consider again the
quantum spaces $\mathcal{A}$, $\mathcal{B}$ and $\mathcal{C}$, and diagrams $%
\left\langle \varphi ,\mathcal{H}\right\rangle \in \Upsilon ^{\mathcal{A},%
\mathcal{B}}$ and $\left\langle \psi ,\mathcal{G}\right\rangle \in \Upsilon
^{\mathcal{B},\mathcal{C}}$, with associated linear spaces (\emph{via }the
functor\emph{\ }$\frak{F}$) 
\begin{equation*}
\mathbf{H}_{1}^{\varphi }=span\left[ h_{i}^{j}\right] _{i,j=1}^{\frak{n},%
\frak{m}}\;\;and\;\;\mathbf{G}_{1}^{\psi }=span\left[ g_{i}^{j}\right]
_{i,j=1}^{\frak{m},\frak{p}}.
\end{equation*}
Denoting by $\mathbf{h}$ and $\mathbf{g}$ the matrices with entries $\left\{
h_{i}^{j}\right\} \subset \mathbf{H}_{1}$ and $\left\{ g_{i}^{j}\right\}
\subset \mathbf{G}_{1}$, respectively, and by $\mathbf{a}$, $\mathbf{b}$ and 
$\mathbf{c}$ the vectors whose components are $\left\{ a_{i}\right\} \subset 
\mathbf{A}_{1}$, $\left\{ b_{i}\right\} \subset \mathbf{B}_{1}$ and $\left\{
c_{i}\right\} \subset \mathbf{C}_{1}$. Then, we can write 
\begin{equation*}
\widehat{\tau }_{\mathcal{A},\mathcal{B}}:\mathbf{b}\otimes \mathbf{h}%
\mapsto \mathbf{b}\otimes \mathbf{\hat{h}}\;\;and\;\;\widehat{\tau }_{%
\mathcal{B},\mathcal{C}}:\mathbf{c}\otimes \mathbf{g}\mapsto \mathbf{c}%
\otimes \mathbf{\hat{g}}
\end{equation*}
being $\mathbf{\hat{h}}=\sigma _{\frak{n}}\cdot \mathbf{h}\cdot \sigma _{%
\frak{m}}^{-1}$ and $\mathbf{\hat{g}}=\sigma _{\frak{m}}\cdot \mathbf{g}%
\cdot \sigma _{\frak{p}}^{-1}$, where $``\cdot "$ \thinspace indicates
matrix multiplication (respects to the chosen basis of $\mathbf{A}_{1}$, $%
\mathbf{B}_{1}$ and $\mathbf{C}_{1}$). There is not contraction between the
above matrices and vector, unless we include the symbols $``\cdot "$ or $``%
\overset{.}{\otimes }"$ in expressions where they appear.

It must be shown that $\left\langle \left( I_{H}\circ \chi \right) \,\varphi
,\mathcal{H}\circ \mathcal{G}\right\rangle $ is an object of $\Upsilon ^{%
\mathcal{A},\mathcal{C}}$, which is true if and only if $\left( I_{H}\circ
\chi \right) \,\varphi $ define the arrow $\mathcal{A}\rightarrow \left( 
\mathcal{H}\circ \mathcal{G}\right) \circ _{\tau }\mathcal{C}$, with $\tau $
given by $\widehat{\tau }_{\mathcal{A},\mathcal{C}}$, which in the basis
above is defined by the assignment 
\begin{equation*}
\widehat{\tau }_{\mathcal{A},\mathcal{C}}:\mathbf{c}\otimes \left( \mathbf{h}%
\overset{.}{\otimes }\mathbf{g}\right) \mapsto \mathbf{c}\otimes \left( 
\mathbf{\hat{h}}\overset{.}{\otimes }\mathbf{\hat{g}}\right) .
\end{equation*}
Let us see firstly that $\widehat{\tau }_{\mathcal{A},\mathcal{C}}$ can be
effectively extended to a twisting map $\tau $ from $\mathbf{C}\otimes
\left( \mathbf{H}\circ \mathbf{G}\right) $ to $\left( \mathbf{H}\circ 
\mathbf{G}\right) \otimes \mathbf{C}$, and then that $\left( I_{H}\circ \chi
\right) \,\varphi $ defines the arrow above. Such an extension would be
given by 
\begin{equation}
\tau :\mathbf{c}_{1}...\mathbf{c}_{r}\otimes \mathbf{h}_{1}...\mathbf{h}_{s}%
\overset{.}{\otimes }\mathbf{g}_{1}...\mathbf{g}_{s}\mapsto \mathbf{\hat{h}}%
_{1}^{\left\langle r\right\rangle }...\;\mathbf{\hat{h}}_{s}^{\left\langle
r\right\rangle }\overset{.}{\otimes }\mathbf{\hat{g}}_{1}^{\left\langle
r\right\rangle }...\;\mathbf{\hat{g}}_{s}^{\left\langle r\right\rangle
}\otimes \mathbf{c}_{1}...\mathbf{c}_{r}  \label{exp}
\end{equation}
where the $\mathbf{h}_{k}$'s are copies of $\mathbf{h}$, the contractions
are done between $\mathbf{h}_{k}$'s and $\mathbf{g}_{l}$'s with $k=l$, and $%
\mathbf{\hat{h}}^{\left\langle r\right\rangle }=\left( \sigma _{\frak{n}%
}\right) ^{r}\cdot \mathbf{h}\cdot \left( \sigma _{\frak{m}}\right) ^{-r}\;$%
(idem $\mathbf{g}$). From expression $\left( \ref{exp}\right) $, it can be
checked directly that $\tau $ is well-defined as a linear map (to this end
we can further extend $\tau $ as a map on $\mathbf{C}_{1}^{\otimes }\otimes
\left( \mathbf{H}_{1}^{\otimes }\otimes \mathbf{G}_{1}^{\otimes }\right) $,
changing in $\left( \ref{exp}\right) $ the symbol $``\overset{.}{\otimes }"$
by $``\otimes "$, and show that $\tau $ respects the involved ideals). At
this point the conic character of the spaces we are considering is crucial.
For general quantum spaces the map $\tau $ given by $\left( \ref{exp}\right) 
$ may be not well-defined (see \textbf{Appendix}). Now, it follows from
straightforward calculations that $\left( I_{H}\circ \chi \right) \,\varphi $
(an homogeneous linear map $\mathbf{A}\rightarrow \mathbf{H}\otimes \mathbf{G%
}\otimes \mathbf{C}$) is a morphism $\mathcal{A}\rightarrow \left( \mathcal{H%
}\circ \mathcal{G}\right) \circ _{\tau }\mathcal{C}$.

Finally, we observe that units $\left\langle \ell _{\mathcal{A}},\mathcal{K}%
\right\rangle $\textsf{\ }are objects of $\Upsilon ^{\mathcal{A}}$. So, the
theorem has been proved.
\end{proof}

\begin{remark}
We see from proof above that, for $\Upsilon ^{\cdot }$ to be a
sub-semigroupoid of $\frak{H}\mathrm{CA}^{\circ }$, factorization property
is a \textbf{generic }condition we must impose on the collection $\left\{ 
\widehat{\tau }\right\} $. Thus, for global collections the mentioned
generic condition is needed in order to have the inclusion $\Upsilon ^{\cdot
}\subset \frak{H}\mathrm{CA}^{\circ }$ of semigroupoids.\ \ $\blacksquare $
\end{remark}

Observe that, in proving $b)\Rightarrow a)$ the global character of the
associated collection that defines $\Upsilon ^{\cdot }$ is not involved. So
we can enunciate the following corollaries (without proof).

\begin{corollary}
Any category $\Upsilon ^{\cdot }$ associated to a factorizable collection is
a sub-semigroupoid of $\frak{H}\mathrm{CA}^{\circ }$.\ \ $\blacksquare $
\end{corollary}

\medskip We call $\circ _{\Upsilon }$ the partial product associated to such
a semigroupoid. Because the embedding $\frak{P}^{\Upsilon }:\Upsilon ^{\cdot
}\hookrightarrow \mathrm{CA}$, such that $\left\langle \varphi ,\mathcal{H}%
\right\rangle \mapsto \mathcal{H}$, satisfies $\frak{P}^{\Upsilon }\,\circ
_{\Upsilon }=\circ \,\left( \frak{P}^{\Upsilon }\times \frak{P}^{\Upsilon
}\right) $ and $\frak{P}^{\Upsilon }\left\langle \ell _{\mathcal{A}},%
\mathcal{K}\right\rangle =\mathcal{K}$ (see \S \textbf{2.2}), it follows:

\begin{corollary}
For every factorizable collection, the function 
\begin{equation*}
\left( \mathcal{A},\mathcal{B}\right) \mapsto \underline{hom}^{\Upsilon }%
\left[ \mathcal{B},\mathcal{A}\right] \backsimeq \mathcal{D}_{\mathcal{B},%
\mathcal{A}}\in \mathrm{CA}
\end{equation*}
defines an $\mathrm{CA}$-cobased or $\mathsf{QLS}$-based category with
arrows \emph{(}opposite to\emph{)} 
\begin{equation*}
\underline{hom}^{\Upsilon }\left[ \mathcal{C},\mathcal{A}\right] \rightarrow 
\underline{hom}^{\Upsilon }\left[ \mathcal{B},\mathcal{A}\right] \circ 
\underline{hom}^{\Upsilon }\left[ \mathcal{C},\mathcal{B}\right]
\end{equation*}
the cocomposition, and for $\underline{end}^{\Upsilon }\left[ \mathcal{A}%
\right] $ the counit epimorphism 
\begin{equation*}
\underline{end}^{\Upsilon }\left[ \mathcal{A}\right] \twoheadrightarrow 
\mathcal{K}\;\;\;/\;\;\;z_{i}^{j}\mapsto \delta _{i}^{j}\;e,
\end{equation*}
and the monomorphic comultiplication 
\begin{equation*}
\underline{end}^{\Upsilon }\left[ \mathcal{A}\right] \hookrightarrow 
\underline{end}^{\Upsilon }\left[ \mathcal{A}\right] \circ \underline{end}%
^{\Upsilon }\left[ \mathcal{A}\right] \;\;\;/\;\;\;z_{i}^{j}\mapsto
z_{i}^{k}\otimes z_{k}^{j}.\;\;\;\blacksquare
\end{equation*}
\end{corollary}

Now, the annunciated lemma to conclude proof of above theorem.

\begin{lemma}
The solutions to Eqs. $\left( \ref{coas}\right) $ and $\left( \ref{cou}%
\right) $ are given by 
\begin{equation*}
\left[ \tau _{\frak{n}}\right] _{ij\beta }^{kl\alpha }=\delta _{i}^{k}\,%
\left[ \overline{\sigma }_{\frak{n}}\right] _{j}^{l}\,\left[ \sigma _{\frak{n%
}}\right] _{\beta }^{\alpha },\;\sigma _{\frak{n}}\in \in GL\left( \frak{n}%
\right) \diagup \mathcal{Z}_{GL\left( \frak{n}\right) }.
\end{equation*}
\end{lemma}

\begin{proof}
We want to solve equations 
\begin{equation}
\left[ \tau _{\frak{n}}\right] _{kyj}^{slx}\,\left[ \tau _{\frak{n}}\right]
_{\alpha pl}^{r\beta q}=\delta _{k}^{s}\delta _{y}^{q}\,\left[ \tau _{\frak{n%
}}\right] _{\alpha pj}^{r\beta x},  \label{ea}
\end{equation}
and 
\begin{equation}
\left[ \tau _{\frak{n}}\right] _{i\alpha \beta }^{kl\alpha }=\delta
_{i}^{k}\delta _{\beta }^{l}.  \label{ebb}
\end{equation}
The solutions will be of the form 
\begin{equation}
\left[ \tau _{\frak{n}}\right] _{kyj}^{slx}=\delta _{k}^{s}\theta _{jy}^{xl}.
\label{c}
\end{equation}
In fact, contracting $p$-$q$ in $\left( \ref{ea}\right) $ we have, using $%
\left( \ref{ebb}\right) $, 
\begin{equation*}
\left[ \tau _{\frak{n}}\right] _{kyj}^{slx}\delta _{\alpha }^{r}=\delta
_{k}^{s}\left[ \tau _{\frak{n}}\right] _{\alpha yj}^{rlx}.
\end{equation*}
Then, taking $\alpha =r$, we arrive at the equality $\left[ \tau _{\frak{n}}%
\right] _{kyj}^{slx}=\delta _{k}^{s}\left[ \tau _{\frak{n}}\right]
_{ryj}^{rlx}$, $1\leq r\leq \frak{n}$ (without contraction in $r$); hence,
we can define $\theta _{jy}^{xl}\doteq \left[ \tau _{\frak{n}}\right]
_{ryj}^{rlx}$ and obtain $\left( \ref{c}\right) $.

Inserting $\left( \ref{c}\right) $ on $\left( \ref{ea}\right) $ and $\left( 
\ref{ebb}\right) $, the latter reduce to 
\begin{equation*}
\theta _{jy}^{lx}\theta _{lp}^{\beta q}\,=\delta _{y}^{q}\,\theta
_{jp}^{\beta x};\;\;\theta _{\beta \alpha }^{l\alpha }=\delta _{\beta }^{l}.
\end{equation*}
If we define $\eta _{jx}^{ly}\doteq \theta _{jy}^{lx}/\frak{n}$, and
contract $q$-$y$ we have $\eta _{jx}^{ly}\eta _{ly}^{\beta p}\,=\eta
_{jx}^{\beta p}$, $\eta _{\beta \alpha }^{l\alpha }=\delta _{\beta }^{l}/%
\frak{n}$. Now, contracting $l$-$\beta $ in the latter, both equations
reduce to $\eta ^{2}=\eta $ and $tr\,\eta =1$, regarding $\eta $ as an $%
\frak{n}^{2}\times \frak{n}^{2}$ matrix. The condition $\eta ^{2}=\eta $
says that $\eta $ is diagonalizable and has eigenvalues $1$ and $0$. If it
has $n_{\lambda }$ eigenvalues equals to $\lambda $ , $\lambda =0,1$, then $%
tr\eta =n_{0}\cdot 0+n_{1}\cdot 1=n_{1}$. But $tr\,\eta =1$, so $n_{1}=1$
and $n_{0}=\frak{n}^{2}-1$. In particular, $\eta $ can be diagonalized to a
matrix $d=diag_{\frak{n}^{2}}(1,0,...,0)$ through some invertible matrix $%
\varkappa $, i.e. 
\begin{equation*}
\eta _{ij}^{kl}=\varkappa _{ij}^{xy}\,d_{xy}^{vw}\,\left( \varkappa
^{-1}\right) _{vw}^{kl}=\varkappa _{ij}^{11}\,\left( \varkappa ^{-1}\right)
_{11}^{kl}.
\end{equation*}
Coming back to $\theta $ we can write $\theta _{il}^{kj}=\frak{n}\,\varkappa
_{ij}^{11}\,\left( \varkappa ^{-1}\right) _{11}^{kl}$, and defining $\left[
\sigma _{\frak{n}}\right] _{i}^{j}\doteq \frak{n}\,\varkappa _{ij}^{11}$ and 
$\left[ \overline{\sigma }_{\frak{n}}\right] _{l}^{k}\doteq \left( \varkappa
^{-1}\right) _{11}^{kl}$, we will have 
\begin{equation*}
\theta _{il}^{kj}=\left[ \sigma _{\frak{n}}\right] _{i}^{j}\,\left[ 
\overline{\sigma }_{\frak{n}}\right] _{l}^{k}.
\end{equation*}
But $\theta _{il}^{kl}=\delta _{i}^{k}$, so $\left[ \sigma _{\frak{n}}\right]
_{i}^{j}\,\left[ \overline{\sigma }_{\frak{n}}\right] _{j}^{k}=\left[ 
\overline{\sigma }_{\frak{n}}\right] _{i}^{j}\,\left[ \sigma _{\frak{n}}%
\right] _{j}^{k}=\delta _{i}^{k}$, i.e. $\overline{\sigma }_{\frak{n}%
}=\sigma _{\frak{n}}^{-1}$. Thus, $\tau _{\frak{n}}$ is a solution of $%
\left( \ref{ea}\right) $ and $\left( \ref{ebb}\right) $ only if 
\begin{equation*}
\left[ \tau _{\frak{n}}\right] _{kyj}^{slx}=\delta _{k}^{s}\,\left[ 
\overline{\sigma }_{\frak{n}}\right] _{y}^{l}\,\left[ \sigma _{\frak{n}}%
\right] _{j}^{x},
\end{equation*}
for some $\sigma _{\frak{n}}\in GL\left( \frak{n}\right) $. Reciprocally, if 
$\sigma _{\frak{n}}\in GL\left( \frak{n}\right) $, it is easy to see that
last expression for $\tau _{\frak{n}}$ gives a solution for $\left( \ref
{coas}\right) $ and $\left( \ref{cou}\right) $.
\end{proof}

\subsubsection{Factorizable collections given by a family of automorphisms}

We are going to study the case in which a given factorizable map 
\begin{equation*}
\widehat{\tau }_{\mathcal{A},\mathcal{B}}=id\otimes \rho \otimes \phi
,\qquad \;\phi =\sigma _{\mathcal{A}},\;\rho ^{-1}=\sigma _{\mathcal{B}},
\end{equation*}
associated to $\Upsilon ^{\mathcal{A},\mathcal{B}}$, is such that $\sigma _{%
\mathcal{A}}:\mathbf{A}_{1}\backsimeq \mathbf{A}_{1}$ and $\sigma _{\mathcal{%
B}}:\mathbf{B}_{1}\backsimeq \mathbf{B}_{1}$ can be extended to quantum
space automorphisms $\mathcal{A}\backsimeq \mathcal{A}$ and $\mathcal{B}%
\backsimeq \mathcal{B}$.

Consider $\mathcal{A}$ and $\mathcal{B}$ with related ideals $\mathbf{I}$
and $\mathbf{J}$ given in Eq. $\left( \ref{ideals}\right) $, and define
another ideals, namely $\mathbf{I}_{\sigma }$ and $\mathbf{J}_{\sigma }$,
linearly generated by 
\begin{equation*}
\left\{ \left\{ _{_{\sigma }}R_{\lambda
_{n}}^{k_{1}...k_{n}}\;a_{k_{1}}...a_{k_{n}}\right\} _{\lambda _{n}\in
\Lambda _{n}}\right\} _{n\in \mathbb{N}_{0}},\;\;\left\{ \left\{ _{_{\sigma
}}S_{\mu _{n}}^{k_{1}...k_{n}}\;b_{k_{1}}...b_{k_{n}}\right\} _{\mu _{n}\in
\Phi _{n}}\right\} _{n\in \mathbb{N}_{0}},
\end{equation*}
respectively, being 
\begin{eqnarray}
_{_{\sigma }}R_{\lambda _{n}}^{k_{1}...k_{n}} &\doteq &R_{\lambda
_{n}}^{k_{1}j_{2}...j_{n}}\;\phi _{j_{2}}^{k_{2}}\left( \phi ^{2}\right)
_{j_{3}}^{k_{3}}...\left( \phi ^{n-1}\right) _{j_{n}}^{k_{n}},  \label{fed}
\\
&&  \notag \\
\;\;_{_{\sigma }}S_{\mu _{n}}^{k_{1}...k_{n}} &\doteq &S_{\mu
_{n}}^{k_{1}j_{2}...j_{n}}\;\left( \rho ^{-1}\right) _{j_{2}}^{k_{2}}\left(
\rho ^{-2}\right) _{j_{3}}^{k_{3}}...\left( \rho ^{1-n}\right)
_{j_{n}}^{k_{n}}.
\end{eqnarray}
Of course, $\mathbf{I}_{\sigma }\backsimeq \mathbf{I}$ and $\mathbf{J}%
_{\sigma }\backsimeq \mathbf{J}$. Therefore, $\left( \mathbf{I}_{\sigma
}\right) _{n}^{\perp }\backsimeq \mathbf{I}_{n}^{\perp }$ and $\left( 
\mathbf{J}_{\sigma }\right) _{n}^{\perp }\backsimeq \mathbf{J}_{n}^{\perp }$
. In particular, the spaces $\left( \mathbf{J}_{\sigma }\right) _{n}^{\perp
} $ will be generated by the set 
\begin{equation*}
\left\{ b^{j_{1}}b^{k_{2}}...b^{k_{n}}\;\rho _{k_{2}}^{j_{2}}\left( \rho
^{2}\right) _{k_{3}}^{j_{2}}...\left( \rho ^{n-1}\right)
_{k_{n}}^{j_{n}}\;\left( S^{\bot }\right) _{j_{1}...j_{n}}^{\omega
_{n}}\right\} _{\omega _{n}\in \Omega _{n}},
\end{equation*}
so we can define 
\begin{equation}
\left( _{_{\sigma }}S^{\bot }\right) _{j_{1}...j_{n}}^{\omega _{n}}\doteq
\rho _{k_{2}}^{j_{2}}\left( \rho ^{2}\right) _{k_{3}}^{j_{2}}...\left( \rho
^{n-1}\right) _{k_{n}}^{j_{n}}\;\left( S^{\bot }\right)
_{j_{1}...j_{n}}^{\omega _{n}}.  \label{fedd}
\end{equation}

\begin{proposition}
If $\sigma _{\mathcal{A}}:\mathbf{A}_{1}\backsimeq \mathbf{A}_{1}$ and $%
\sigma _{\mathcal{B}}:\mathbf{B}_{1}\backsimeq \mathbf{B}_{1}$ can be
extended to quantum space automorphisms $\mathcal{A}\backsimeq \mathcal{A}$
and $\mathcal{B}\backsimeq \mathcal{B}$, we can define 
\begin{equation*}
\begin{array}{c}
\underline{hom}^{\Upsilon }\left[ \mathcal{B},\mathcal{A}\right] \doteq 
\mathcal{B}^{\Upsilon }\triangleright \mathcal{A}^{\Upsilon }, \\ 
\\ 
\mathcal{A}^{\Upsilon }\doteq \left( \mathbf{A}_{1},\left. \mathbf{A}%
_{1}^{\otimes }\right/ \mathbf{I}_{\sigma }\right) ,\;\;\mathcal{B}%
^{\Upsilon }\doteq \left( \mathbf{B}_{1},\left. \mathbf{B}_{1}^{\otimes
}\right/ \mathbf{J}_{\sigma }\right) ,
\end{array}
\end{equation*}
and for $\mathcal{A}$,$\mathcal{B}\in \mathrm{QA}$ \emph{(}or any $\mathrm{CA%
}^{m}$, $m\in \mathbb{N}$\emph{)}, see Eq. $\left( \ref{compa}\right) $, 
\begin{equation*}
\underline{hom}^{\Upsilon }\left[ \mathcal{B},\mathcal{A}\right] \doteq
\left( \mathcal{B}^{\Upsilon }\right) ^{!}\bullet \mathcal{A}^{\Upsilon }.
\end{equation*}
In particular, $\underline{hom}^{\Upsilon }\left[ \mathcal{K},\mathcal{A}%
\right] =\mathcal{A}^{\Upsilon }$.
\end{proposition}

\begin{proof}
The maps $\sigma _{\mathcal{A}}$ and $\sigma _{\mathcal{B}}$ can be extended
to algebra isomorphisms \emph{iff} for every $n\in \mathbb{N}$ there exists
numbers $C_{\lambda _{n}}^{\lambda _{n}^{\prime }}$'s and $D_{\mu _{n}}^{\mu
_{n}^{\prime }}$'s (defining invertible matrices in $GL\left( \Lambda
_{n}\right) $ and $GL\left( \Phi _{n}\right) $, resp.) such that 
\begin{eqnarray}
R_{\lambda _{n}}^{k_{1}...k_{n}}\;\phi _{k_{1}}^{j_{1}}...\phi
_{k_{n}}^{j_{n}} &=&C_{\lambda _{n}}^{\lambda _{n}^{\prime }}\;R_{\lambda
_{n}^{\prime }}^{j_{1}...j_{n}},  \label{conc} \\
&&  \notag \\
S_{\mu _{n}}^{k_{1}...k_{n}}\;\left( \rho ^{-1}\right)
_{k_{1}}^{j_{1}}...\left( \rho ^{-1}\right) _{k_{n}}^{j_{n}} &=&D_{\mu
_{n}}^{\mu _{n}^{\prime }}\;S_{\mu _{n}^{\prime }}^{j_{1}...j_{n}}.
\end{eqnarray}
In such a case, because $\sigma _{\mathcal{A}}^{\ast }:\mathbf{A}_{1}^{\ast
}\backsimeq \mathbf{A}_{1}^{\ast }$ and $\sigma _{\mathcal{B}}^{\ast }:%
\mathbf{B}_{1}^{\ast }\backsimeq \mathbf{B}_{1}^{\ast }$ could be extended
to algebra homomorphisms $\mathbf{A}^{!}\backsimeq \mathbf{A}^{!}$ and $%
\mathbf{B}^{!}\backsimeq \mathbf{B}^{!}$, we will have, in particular (for $%
\mathbf{B}$), numbers $E_{\omega _{n}}^{\omega _{n}^{\prime }}$'s such that 
\begin{equation*}
\rho _{k_{2}}^{j_{2}}\rho _{k_{3}}^{j_{2}}...\rho _{k_{n}}^{j_{n}}\;\left(
S^{\bot }\right) _{k_{1}...j_{n}}^{\omega _{n}}=\left( S^{\bot }\right)
_{j_{1}...j_{n}}^{\omega _{n}}\;E_{\omega _{n}}^{\omega _{n}^{\prime }}.
\end{equation*}
Remember that for a factorizable map the underlying algebras of the initial
objects of $\Upsilon ^{\mathcal{A},\mathcal{B}}$ are isomorphic to the one
given by Eqs. $\left( \ref{A}\right) $ and $\left( \ref{B}\right) $. Hence,
from $\left( \ref{conc}\right) $ and the last eq., it is obvious that (as we
affirm in \S \textbf{4.1.2}) Eq. $\left( \ref{B}\right) $ reduces to 
\begin{equation}
\left\{ \left\{ R_{\lambda _{n}}^{k_{1}...k_{n}}\;z_{k_{1}}^{j_{1}}\cdot 
\check{z}_{k_{2}}^{j_{2}}\;...\;\check{z}_{k_{n}}^{\left\langle
n-1\right\rangle j_{n}}\;\left( S^{\bot }\right) _{j_{1}...j_{n}}^{\omega
_{n}}\right\} _{\omega _{n}\in \Omega _{n}}^{\lambda _{n}\in \Lambda
_{n}}\right\} _{n\in \mathbb{N}_{0}},  \label{C}
\end{equation}
where we have to read $z_{i}^{j}=b^{j}\otimes a_{i}$. This means that $%
\mathbf{L}$ is generated by 
\begin{equation*}
\left\{ \left\{ _{\sigma }R_{\lambda
_{n}}^{k_{1}...k_{n}}\;z_{k_{1}}^{j_{1}}\cdot
z_{k_{2}}^{j_{2}}\;...\;z_{k_{n}}^{j_{n}}\;\left( _{\sigma }S^{\bot }\right)
_{j_{1}...j_{n}}^{\omega _{n}}\right\} _{\omega _{n}\in \Omega
_{n}}^{\lambda _{n}\in \Lambda _{n}}\right\} _{n\in \mathbb{N}_{0}},
\end{equation*}
and immediately that $\mathbf{D}\backsimeq \left( \left. \mathbf{B}%
_{1}^{\otimes }\right/ \mathbf{J}_{\sigma }\right) \triangleright \left(
\left. \mathbf{A}_{1}^{\otimes }\right/ \mathbf{I}_{\sigma }\right) $.

To show the equality $\mathcal{A}^{\Upsilon }=\underline{hom}^{\Upsilon }%
\left[ \mathcal{K},\mathcal{A}\right] $, it is enough to realize that $%
\sigma _{\mathcal{K}}:\Bbbk \backsimeq \Bbbk $ always can be extended to an
automorphism $\mathcal{K}\backsimeq \mathcal{K}$ and always $\mathcal{K}%
^{\Upsilon }=\mathcal{K}$. So, $\mathcal{K}^{\Upsilon }\triangleright 
\mathcal{A}^{\Upsilon }=\mathcal{A}^{\Upsilon }$. In particular, $\underline{%
hom}^{\Upsilon }\left[ \mathcal{K},\mathcal{K}\right] =\mathcal{K}$.
\end{proof}

As an example, let us take $\mathcal{A}=\mathcal{B}=\mathcal{A}_{R}$, i.e.
the quadratic quantum space generated by $\mathbf{A}_{1}=span\left[ a_{i}%
\right] _{i=1}^{n}$ and restricted to relations 
\begin{equation}
\mathbf{I}_{R}=span\left[ R_{ij}^{kl}\;a_{k}\otimes a_{l}-a_{j}\otimes a_{i}%
\right] _{i,j=1}^{n},  \label{idr}
\end{equation}
being $R=R_{su\left( n\right) }$, the $R$-matrix related to the Lie algebra $%
su\left( n\right) $ and the quantum groups $GL_{q}\left( n\right) $ and $%
SL_{q}\left( n\right) $ \cite{frt}. Recall that $GL_{q}\left( n\right) $ is
a quadratic bialgebra generated by symbols $T_{i}^{j}$ satisfying 
\begin{equation}
R_{ij}^{kl}\;T_{k}^{n}T_{l}^{m}-T_{j}^{l}T_{i}^{k}\;R_{kl}^{nm}.
\label{ifrt}
\end{equation}
It coacts upon $\mathcal{A}_{R}$ through the map $a_{i}\mapsto
T_{i}^{k}\otimes a_{k}$, and consequently, $GL_{q}\left( n\right) $ defines
a diagram in the category $\left( \mathcal{A}\downarrow \mathrm{QA}\circ 
\mathcal{A}\right) $. It follows from the initiality of $\underline{end}%
\left[ \mathcal{A}_{R}\right] $ that there exists an epimorphism $\underline{%
end}\left[ \mathcal{A}_{R}\right] \twoheadrightarrow GL_{q}\left( n\right) $%
. Now, for the bijection $\widehat{\tau }_{\mathcal{A}_{R}}=id\otimes \phi
^{-1}\otimes \phi $, with $\phi =\sigma _{\mathcal{A}_{R}}:\mathbf{A}%
_{1}\rightarrow \mathbf{A}_{1}$ given by an invertible diagonal matrix, we
have that $\left[ R,\phi \otimes \phi \right] =0$. Then, Eq. $\left( \ref
{conc}\right) $ holds, and 
\begin{equation*}
\underline{end}^{\Upsilon }\left[ \mathcal{A}_{R}\right] =\left( \left( 
\mathcal{A}_{R}\right) ^{\Upsilon }\right) ^{!}\bullet \left( \mathcal{A}%
_{R}\right) ^{\Upsilon }=\mathcal{A}_{_{\sigma }R}^{!}\bullet \mathcal{A}%
_{_{\sigma }R}.
\end{equation*}
Moreover, since straightforwardly $\mathcal{A}_{_{\sigma }R}=\mathcal{A}%
_{R^{\phi }}$, being 
\begin{equation*}
R^{\phi }=R_{su\left( n\right) }^{\phi }=\left( \phi ^{-1}\otimes id\right)
\,R_{su\left( n\right) }\left( id\otimes \phi \right) ,
\end{equation*}
it follows that 
\begin{equation*}
\underline{end}^{\Upsilon }\left[ \mathcal{A}_{R}\right] =\mathcal{A}%
_{R^{\phi }}^{!}\bullet \mathcal{A}_{R^{\phi }}=\underline{end}\left[ 
\mathcal{A}_{R^{\phi }}\right] .
\end{equation*}
$R^{\phi }$ defines precisely a (multiparametric) quantum group $GL_{q,\phi
}\left( n\right) $ \cite{demindov}\cite{man2}, so we have, as before, an
epimorphic map $\underline{end}^{\Upsilon }\left[ \mathcal{A}_{R}\right]
\twoheadrightarrow GL_{q,\phi }\left( n\right) $. In particular, we can say
that $\left[ GL_{q,\phi }\left( n\right) \right] ^{op}$ is a \emph{quantum
subspace} of $\underline{End}^{\Upsilon }\left[ \mathcal{A}_{R}^{op}\right] $%
.

\subsection{Twisted internal coEnd objects}

In this section, we shall connect $\underline{end}$'s and $\underline{end}%
^{\Upsilon }$'s by \emph{2-cocycle twisting }of bialgebras \cite{drin}, as $%
GL_{q}\left( n\right) $ and its multiparametric versions. We actually use a
dual approach with respect to Drinfeld one, following the book of Majid \cite
{maj} and references therein.

Let us restrict ourself to collections $\left\{ \sigma _{\mathcal{A}%
}\right\} _{\mathcal{A}\in \mathrm{CA}}$ that define automorphisms 
\begin{equation*}
\left\{ \sigma :\mathcal{A}\backsimeq \mathcal{A}\right\} _{\mathcal{A}\in 
\mathrm{CA}}.
\end{equation*}
We shall see each object $\underline{end}^{\Upsilon }\left[ \mathcal{A}%
\right] =\mathcal{A}^{\Upsilon }\triangleright \mathcal{A}^{\Upsilon }$,
endowed with the bialgebra structure 
\begin{equation*}
\Delta :\underline{end}^{\Upsilon }\left[ \mathcal{A}\right] \hookrightarrow 
\underline{end}^{\Upsilon }\left[ \mathcal{A}\right] \circ \underline{end}%
^{\Upsilon }\left[ \mathcal{A}\right] ;\;\;\;\varepsilon :\underline{end}%
^{\Upsilon }\left[ \mathcal{A}\right] \twoheadrightarrow \mathcal{K},
\end{equation*}
is a counital 2-cocycle twisting of the proper coHom object $\underline{end}%
\left[ \mathcal{A}\right] $; and that the same twisting induces one on $%
\mathcal{A}$ by the coaction\textbf{\ }$\mathcal{A}\hookrightarrow 
\underline{end}\left[ \mathcal{A}\right] \circ \mathcal{A}$, giving the
quantum space $\mathcal{A}^{\Upsilon }$ (c.f. \cite{maj}, page 54).

\bigskip

Consider a pair $\left( \mathbf{A}_{1},\mathbf{A}\right) $ with related
automorphism $\sigma _{\mathcal{A}}$ given by $a_{i}\mapsto \sigma
_{i}^{j}\,a_{j}$. Let us write $z_{i}^{j}=a^{j}\otimes a_{i}$ and define on $%
\underline{end}\left[ \mathcal{A}\right] =\mathcal{A}\triangleright \mathcal{%
A}$ the linear 2-form $\,\chi :\left( \mathbf{A}\triangleright \mathbf{A}%
\right) \otimes \left( \mathbf{A}\triangleright \mathbf{A}\right)
\rightarrow \Bbbk $ as\footnote{%
It is well-defined because of Eq. $\left( \ref{conc}\right) $ for $\phi
=\sigma $.} 
\begin{equation}
\begin{array}{l}
\,\chi \left( \mathbf{z}_{1}...\mathbf{z}_{r}\otimes \mathbf{z}_{r+1}...%
\mathbf{z}_{r+s}\right) =\left( \sigma ^{-r}\right) _{r+1}...\left( \sigma
^{-r}\right) _{r+s}, \\ 
\\ 
\,\chi \left( 1\otimes x\right) =\,\chi \left( x\otimes 1\right)
=\varepsilon \left( x\right) \in \Bbbk ,\;\forall x\in \mathbf{A}%
\triangleright \mathbf{A}.
\end{array}
\label{defcoc}
\end{equation}
The first eq. must be understood as, for $r=2$ and $s=3$, 
\begin{equation*}
\,\chi \left( z_{i_{1}}^{j_{1}}\cdot z_{i_{2}}^{j_{2}}\otimes
z_{i_{3}}^{j_{3}}\cdot z_{i_{4}}^{j_{4}}\cdot z_{i_{5}}^{j_{5}}\right)
=\left( \sigma ^{-2}\right) _{i_{3}}^{j_{3}}\,\left( \sigma ^{-2}\right)
_{i_{4}}^{j_{4}}\,\left( \sigma ^{-2}\right) _{i_{5}}^{j_{5}}\in \Bbbk .
\end{equation*}
The last eq. says that $\,\chi $ is counital. By $\ast $ we indicate the
convolution product of linear forms on a bialgebra.

\begin{proposition}
If $\left( \mu ,\eta ,\Delta ,\varepsilon \right) $ is the bialgebra
structure of $\mathbf{A}\triangleright \mathbf{A}$, we can define another
structure $\left( \mu _{\,\chi },\eta ,\Delta ,\varepsilon \right) $, where $%
\mu _{\,\chi }\doteq \,\chi \ast \mu \ast \,\chi ^{-1}$.
\end{proposition}

This new bialgebra is denoted $\left( \mathbf{A}\triangleright \mathbf{A}%
\right) _{\,\chi }$ and called the twisting of $\mathbf{A}\triangleright 
\mathbf{A}$ by $\,\chi $.

\begin{proof}
The counital property for $\,\chi $ insures that $\eta $ is a unit for $\mu
_{\,\chi }$. We just must prove $\mu _{\,\chi }$ is associative. This
follows from the fact that $\,\chi $ is a 2-cocycle, i.e. $\,\chi _{12}\ast
\,\chi \,\left( \mu \otimes I\right) =\,\chi _{23}\ast \,\chi \,\left(
I\otimes \mu \right) $; where, $\chi _{12}=\,\chi \otimes \varepsilon $, $%
\,\chi _{23}=\varepsilon \otimes \,\chi $ and $\,\chi _{13}=\,\chi
_{12}\,\left( I\otimes \tau _{o}\right) $. In fact, acting with $\chi
_{12}\ast \,\chi \,\left( \mu \otimes I\right) $ upon the element 
\begin{equation*}
\mathbf{z}_{1}...\mathbf{z}_{r}\otimes \mathbf{z}_{r+1}...\mathbf{z}%
_{r+s}\otimes \mathbf{z}_{r+s+1}...\mathbf{z}_{r+s+t}\in \left( \mathbf{A}%
\triangleright \mathbf{A}\right) ^{\otimes 3}
\end{equation*}
we obtain the matrix 
\begin{equation*}
\left( \sigma ^{-r}\right) _{r+1}...\left( \sigma ^{-r}\right) _{r+s}\left(
\sigma ^{-r-s}\right) _{r+s+1}...\left( \sigma ^{-r-s}\right) _{r+s+t};
\end{equation*}
and acting with $\,\,\chi _{23}\ast \,\chi \,\left( I\otimes \mu \right) $
we have 
\begin{equation*}
\begin{array}{l}
\left( \sigma ^{-s}\right) _{r+s+1}...\left( \sigma ^{-s}\right)
_{r+s+t}\cdot \left( \sigma ^{-r}\right) _{r+1}...\left( \sigma ^{-r}\right)
_{r+s+t}= \\ 
\\ 
=\left( \sigma ^{-r}\right) _{r+1}...\left( \sigma ^{-r}\right) _{r+s}\left(
\sigma ^{-r-s}\right) _{r+s+1}...\left( \sigma ^{-r-s}\right) _{r+s+t},
\end{array}
\end{equation*}
where $\cdot $ \thinspace indicates matrix multiplication. Thus $\,\chi $ is
a 2-cocycle. Therefore, $\left( \mu _{\,\chi },\eta \right) $ define a
unital algebra structure.

$\Delta $ and $\varepsilon $ are algebra morphisms for $\left( \mu _{\,\chi
},\eta \right) $ (see \cite{maj} for a proof), and consequently $\left( \mu
_{\,\chi },\eta ,\Delta ,\varepsilon \right) $ is a bialgebra structure on
the linear space $\mathbf{A}\triangleright \mathbf{A}$, as we wanted to show.
\end{proof}

Now, we are going to show that $\left( \mathbf{A}\triangleright \mathbf{A}%
\right) _{\,\chi }\ $is a quantum space, which is isomorphic to $\mathcal{A}%
^{\Upsilon }\triangleright \mathcal{A}^{\Upsilon }$ as quantum spaces and as
bialgebras. By definition, the underlying vector space of $\left( \mathbf{A}%
\triangleright \mathbf{A}\right) _{\,\chi }$ is $\mathbf{A}\triangleright 
\mathbf{A}$, so is given by monomials in $z_{i}^{j}=a^{j}\otimes a_{i}$
under the product $\mu $, satisfying (see Eq. $\left( \ref{ideal}\right) $
of \S \textbf{1.3}) 
\begin{equation}
R_{\lambda _{n}}^{k_{1}...k_{n}}\;z_{k_{1}}^{j_{1}}\cdot
z_{k_{2}}^{j_{2}}\cdot ...\cdot z_{k_{n}}^{j_{n}}\;\left( R^{\bot }\right)
_{j_{1}...j_{n}}^{\omega _{n}}=0,  \label{x}
\end{equation}
which can be written 
\begin{equation}
R_{\lambda _{n}}^{k_{1}...k_{n}}\;z_{k_{1}}^{j_{1}}\cdot _{\,\chi }\check{z}%
_{k_{2}}^{j_{2}}\cdot _{\,\chi }...\cdot _{\,\chi }\check{z}%
_{k_{n}}^{\left\langle n-1\right\rangle j_{n}}\;\left( R^{\bot }\right)
_{j_{1}...j_{n}}^{\omega _{n}}=0  \label{y}
\end{equation}
or 
\begin{equation}
_{\sigma }R_{\lambda _{n}}^{k_{1}...k_{n}}\;z_{k_{1}}^{j_{1}}\cdot _{\,\chi
}z_{k_{2}}^{j_{2}}\cdot _{\,\chi }...\cdot _{\,\chi
}z_{k_{n}}^{j_{n}}\;\left( _{\sigma }R^{\bot }\right)
_{j_{1}...j_{n}}^{\omega _{n}}=0,  \label{z}
\end{equation}
where $\cdot _{\,\chi }\doteq \mu _{\,\chi }$.\footnote{%
To see the equivalence of Eqs. $\left( \ref{x}\right) $ and $\left( \ref{y}%
\right) $ note that $z_{i}^{j}\cdot _{\,\chi }z_{k}^{l}=z_{i}^{j}\cdot 
\check{z}_{k}^{l}$ (and use $\left( \ref{defcoc}\right) $ for higher order
monomials); and the equivalence of those equations with $\left( \ref{z}%
\right) $ follows from definition of $_{\sigma }R$ and $_{\sigma }R^{\bot }$
(see Eqs. $\left( \ref{fed}\right) $ and $\left( \ref{fedd}\right) $).}
Clearly, the monomials in $z_{i}^{j}$ under the product $\mu _{\,\chi }$
also span $\left( \mathbf{A}\triangleright \mathbf{A}\right) _{\,\chi }$.
Then, we have a quantum space 
\begin{equation*}
\left( \mathcal{A}\triangleright \mathcal{A}\right) _{\,\chi }\doteq \left( 
\mathbf{A}_{1}^{\ast }\otimes \mathbf{A}_{1},\left( \mathbf{A}\triangleright 
\mathbf{A}\right) _{\,\chi }\right)
\end{equation*}
with algebra structure $\left( \mu _{\,\chi },\eta \right) $ and related
ideal generated by the elements given in Eq. $\left( \ref{z}\right) $, with $%
\lambda _{n}\in \Lambda _{n}$, $\omega _{n}\in \Omega _{n}$, $n\in \mathbb{N}
$. Because the elements that define the ideal of $\mathcal{A}^{\Upsilon
}\triangleright \mathcal{A}^{\Upsilon }$ are 
\begin{equation*}
_{\sigma }R_{\lambda _{n}}^{k_{1}...k_{n}}\;z_{k_{1}}^{j_{1}}\cdot
z_{k_{2}}^{j_{2}}\cdot ...\cdot z_{k_{n}}^{j_{n}}\;\left( _{\sigma }R^{\bot
}\right) _{j_{1}...j_{n}}^{\omega _{n}},
\end{equation*}
it is obvious that $\mathcal{A}^{\Upsilon }\triangleright \mathcal{A}%
^{\Upsilon }$ and $\left( \mathcal{A}\triangleright \mathcal{A}\right)
_{\,\chi }$ are isomorphic quantum spaces. On the other hand, the coalgebra
structure are identical, hence they are isomorphic as bialgebras as well.
The following theorem resumes these results.

\begin{theorem}
The twisting $\left( \mathbf{A}\triangleright \mathbf{A}\right) _{\,\chi }$ 
\emph{(}of the bialgebra $\mathbf{A}\triangleright \mathbf{A}$ by $\,\chi $%
\emph{) }defines a pair 
\begin{equation*}
\underline{end}\left[ \mathcal{A}\right] _{\,\chi }=\left( \mathcal{A}%
\triangleright \mathcal{A}\right) _{\,\chi }\doteq \left( \mathbf{A}%
_{1}^{\ast }\otimes \mathbf{A}_{1},\left( \mathbf{A}\triangleright \mathbf{A}%
\right) _{\,\chi }\right)
\end{equation*}
isomorphic to $\mathcal{A}^{\Upsilon }\triangleright \mathcal{A}^{\Upsilon }=%
\underline{end}^{\Upsilon }\left[ \mathcal{A}\right] $; in other terms,
there exists a counital 2-cocycle $\chi $ such that 
\begin{equation*}
\underline{end}^{\Upsilon }\left[ \mathcal{A}\right] \backsimeq \underline{%
end}\left[ \mathcal{A}\right] _{\,\chi }.\;\;\;\blacksquare
\end{equation*}
\end{theorem}

Finally, consider in $\mathbf{A}$ the product 
\begin{equation*}
m_{\,\chi }=\left( \,\chi \otimes m\right) \,\left( I\otimes \tau
_{o}\otimes I\right) \,\left( \delta \otimes \delta \right) ,
\end{equation*}
where $\delta $ is the proper coevaluation $\mathcal{A}\rightarrow 
\underline{end}\left[ \mathcal{A}\right] \circ \mathcal{A}$ and $\tau _{o}$
is the flipping map. In particular, 
\begin{equation*}
a_{i}\cdot _{\chi }a_{j}\doteq m_{\,\chi }\left( a_{i}\otimes a_{j}\right)
=\left( \sigma ^{-1}\right) _{j}^{k}\;a_{i}\cdot a_{k}.
\end{equation*}
Since $\chi $ is a counital 2-cocycle, $m_{\,\chi }$ is associative with the
same unit as $m$, and $\mathbf{A}_{1}$ generates $\mathbf{A}$ under $%
m_{\,\chi }$. Then, the algebra $\left( \mathbf{A},m_{\,\chi }\right) $
(with the same unit as $\left( \mathbf{A},m\right) $) is a quantum space,
namely $\mathcal{A}_{\,\chi }$. Using $m_{\,\chi }$, the relations 
\begin{equation*}
R_{\lambda _{n}}^{k_{1}...k_{n}}\;a_{k_{1}}\cdot ...\cdot a_{k_{n}}=0
\end{equation*}
can be expressed as 
\begin{equation}
_{\sigma }R_{\lambda _{n}}^{k_{1}...k_{n}}\;a_{k_{1}}\cdot _{\,\chi
}...\cdot _{\,\chi }a_{k_{n}}=0,  \label{tws}
\end{equation}
then,

\begin{theorem}
The quantum space $\mathcal{A}^{\Upsilon }=\underline{hom}^{\Upsilon }\left[ 
\mathcal{K},\mathcal{A}\right] $ is isomorphic to the twisting $\mathcal{A}%
_{\,\chi }$ of $\mathcal{A}$ induced by the coaction $\mathbf{A}\rightarrow
\left( \mathbf{A}\triangleright \mathbf{A}\right) \circ \mathbf{A}$ and the
2-form $\,\chi $ \emph{(}given in $\left( \ref{defcoc}\right) $\emph{)}; in
other words, $\mathcal{A}^{\Upsilon }\backsimeq \mathcal{A}_{\,\chi
}.\;\;\;\;\blacksquare $
\end{theorem}

\begin{corollary}
Since $\underline{hom}^{\Upsilon }\left[ \mathcal{B},\mathcal{A}\right] =%
\mathcal{B}^{\Upsilon }\triangleright \mathcal{A}^{\Upsilon }$, it follows
that 
\begin{equation*}
\underline{hom}^{\Upsilon }\left[ \mathcal{B},\mathcal{A}\right] \backsimeq 
\mathcal{B}_{\,\chi }\triangleright \mathcal{A}_{\,\chi }=\underline{hom}%
\left[ \mathcal{B}_{\,\chi },\mathcal{A}_{\,\chi }\right] .
\end{equation*}
For coEnd objects we have in addition 
\begin{equation*}
\underline{end}^{\Upsilon }\left[ \mathcal{A}\right] \backsimeq \mathcal{A}%
_{\,\chi }\triangleright \mathcal{A}_{\,\chi }=\underline{end}\left[ 
\mathcal{A}_{\,\chi }\right] \backsimeq \underline{end}\left[ \mathcal{A}%
\right] _{\,\chi }.\;\;\;\;\blacksquare
\end{equation*}
\end{corollary}

\section*{Conclusions}

We carried out the construction of $\underline{hom}^{\Upsilon }\left[ 
\mathcal{B},\mathcal{A}\right] $ looking for the quantum spaces giving, in a
universal way, a notion of `quantum space of transformations' when one
allows for some non commutativity among its \emph{points} and the points of
the space to be transformed. In doing so, we defined a twisted tensor
product between objects of the category $\mathrm{FGA}$ of (general) quantum
spaces (in an analogous way that the one defined by \v{C}ap and his
collaborators for unital algebras). Making use of the Manin algebraic
geometric terminology, we were able to connect the idea of twisted tensor
product structures with the mentioned non commutativity. Then, among certain
subclasses of maps $\mathcal{A}\rightarrow \mathcal{H}\circ _{\tau }\mathcal{%
B}$, we find universal elements, defined by objects $\underline{hom}%
^{\Upsilon }\left[ \mathcal{B},\mathcal{A}\right] \in \mathrm{FGA}$, giving
us for each subclass the notion of space of morphisms we are looking for.
Moreover, on a subcategory $\mathrm{CA}\subset \mathrm{FGA}$ of quantum
spaces, the opposite objects to $\underline{hom}^{\Upsilon }\left[ \mathcal{B%
},\mathcal{A}\right] $ define a $\mathsf{QLS}$-based category.

We showed that under certain circumstances, the bialgebras $\underline{end}%
^{\Upsilon }\left[ \mathcal{A}\right] $ are 2-cocycle twisting of $%
\underline{end}\left[ \mathcal{A}\right] $ (the proper internal coEnd
objects of $\mathrm{CA}$). In a forthcoming paper \cite{gm} we shall define
twisting transformations of quantum spaces, partially controlled by a
(multiplicative) cochain quasicomplex. In these terms, we will see that all
objects $\underline{hom}^{\Upsilon }\left[ \mathcal{B},\mathcal{A}\right] $
(under similar circumstances) are twisting by 2-cocycles of the proper coHom
ones. Then, we can say that both objects are \emph{twist} or \emph{gauge
equivalent}, in an analogous sense to the gauge equivalence that Drinfeld
defined on quasi-Hopf algebras. On the other hand, symmetric twisted tensor
products $\mathcal{A}\circ _{\tau }\mathcal{B}$ can be seen as particular
2-cocycle twisting of the quantum space $\mathcal{A}\circ \mathcal{B}$,
which enable us to develope a generalization of the concept of twisted coHom
objects in the setting of twisting of quantum spaces, as we shall discuss in 
\cite{gg}.

\section*{Acknowledgments}

H.M. thanks to CONICET and S.G. thanks to CNEA, Argentina, for financial
support. The authors also thank to L. Bruschi for useful discussions.

\section*{Appendix}

Consider the general quantum spaces $\mathcal{A}=\left( \mathbf{A}_{1},%
\mathbf{A}_{1}^{\otimes }\right) $, $\mathcal{B}=\left( \mathbf{B}_{1},%
\mathbf{B}_{1}^{\otimes }\right) $ and $\mathcal{C}=\left( \mathbf{C}_{1},%
\mathbf{C}\right) $, and the basis $\left\{ a_{i}\right\} \subset \mathbf{A}%
_{1}$, $\left\{ b_{i}\right\} \subset \mathbf{B}_{1}$ and $\left\{
c_{i}\right\} \subset \mathbf{C}_{1}$. The algebra of $\mathcal{C}$ is given
by $\mathbf{C}\doteq \mathbf{C}_{1}^{\otimes }\left/ I\left[
D^{k}\,c_{k}+\lambda \,1\right] \right. $,$\;$with $\lambda \neq 0$. So, $%
\mathbf{C}$ is not a graded algebra. Let $\left\langle \varphi ,\mathcal{H}%
\right\rangle $ and $\left\langle \psi ,\mathcal{G}\right\rangle $ be
initial diagrams in $\Upsilon ^{\mathcal{A},\mathcal{B}}$ and $\Upsilon ^{%
\mathcal{B},\mathcal{C}}$, respectively, where $\Upsilon ^{\cdot }$ is built
up, just for simplicity, from a factorizable collection. Because $\mathcal{A}
$ and $\mathcal{B}$ are given by free algebras, we can take $\mathcal{H}%
=\left( \mathbf{B}_{1}^{\ast }\otimes \mathbf{A}_{1},\left[ \mathbf{B}%
_{1}^{\ast }\otimes \mathbf{A}_{1}\right] ^{\otimes }\right) $, and $\varphi
_{1}=\left. \varphi \right| _{\mathbf{A}_{1}}:a_{i}\mapsto b^{j}\otimes
a_{i}\otimes b_{j}$. On the other hand, the algebra $\mathbf{G}$ of $%
\mathcal{G}=\left( \mathbf{G}_{1},\mathbf{G}\right) $ can be taken as a
quotient of $\left[ \mathbf{C}_{1}^{\ast }\otimes \mathbf{B}_{1}\right]
^{\otimes }$, where such a quotient is needed to make the corresponding $%
\tau $ of $\Upsilon ^{\mathcal{B},\mathcal{C}}$ well-defined. In particular,
calling $g_{j}^{l}$ the generators of $\mathbf{G}_{1}$, we must have 
\begin{equation*}
\widehat{\tau }_{\mathcal{B},\mathcal{C}}\left( \left( D^{k}\,c_{k}+\lambda
\,1\right) \otimes g_{j}^{l}\right) =D^{k}\,c_{k}\otimes \hat{g}%
_{j}^{l}+\lambda \,1\otimes g_{j}^{l}=\lambda \,1\otimes \left( g_{j}^{l}-%
\hat{g}_{j}^{l}\right) =0,
\end{equation*}
what implies that $\hat{g}_{j}^{l}-g_{j}^{l}=0$. Then, the symbols $%
g_{j}^{l} $ are not linearly independent. Let $\left\{ y_{\mu }\right\} $ be
a basis for $\mathbf{G}_{1}$, such that $g_{j}^{l}=\alpha _{j}^{l\mu
}\,y_{\mu }$, and let us denote by $h_{n}^{j}$ the basis of $\mathbf{B}%
_{1}^{\ast }\otimes \mathbf{A}_{1}$. If we want to define $\widehat{\tau }_{%
\mathcal{A},\mathcal{C}}$ on $\mathbf{H}\circ \mathbf{G}\otimes \mathbf{C}$,
we need that 
\begin{equation*}
\begin{array}{l}
\widehat{\tau }_{\mathcal{A},\mathcal{C}}\left( \left( D^{k}\,c_{k}+\lambda
\,1\right) \otimes \left( h_{n}^{j}\otimes g_{j}^{l}\right) \right) =\lambda
\,1\otimes \left( h_{n}^{j}\otimes g_{j}^{l}-\hat{h}_{n}^{j}\otimes \hat{g}%
_{j}^{l}\right) \\ 
\\ 
=\lambda \,1\otimes \left( h_{n}^{j}-\hat{h}_{n}^{j}\right) \otimes
g_{j}^{l}=\lambda \,1\otimes \left( h_{n}^{j}-\hat{h}_{n}^{j}\right) \alpha
_{j}^{l\mu }\otimes y_{\mu }=0,
\end{array}
\end{equation*}
which is true if and only if $\left( h_{n}^{j}-\hat{h}_{n}^{j}\right)
\,\alpha _{j}^{l\mu }=0$. But the symbols $h_{n}^{j}$ are linearly
independent, hence, in general, $\widehat{\tau }_{\mathcal{A},\mathcal{C}}$
is not well-defined on $\mathbf{H}\circ \mathbf{G}\otimes \mathbf{C}$.

\end{document}